\documentclass[12pt,a4]{article}
\usepackage{amssymb,amsmath,amsthm,amsfonts,amscd,amstext}
\usepackage[mathscr]{eucal}
\usepackage{array}
\usepackage{graphicx}
\newtheorem{lem}{Lemma}[section]
\newtheorem{prop}[lem]{Proposition}
\newtheorem{thm}[lem]{Theorem}
\newtheorem{cor}[lem]{Corollary}
\newtheorem{df}{Definition}[section]

\title{\bf Automata in groups and  dynamics and
induced systems of PDE in tropical geometry }
\author{\bf Tsuyoshi Kato}
\date{}

\begin{document}
\maketitle

\section{Introduction}
Automata appear in many branches of mathematics.
From the view point of dynamical systems, they produce actions on spaces.
Of particular interest for us is {\em theory of automata groups} which studies
discrete group actions on trees.

In this paper 
we develop a new construction  of dynamical scaling transform 
by use of {\em tropical geometry}, and apply 
 automata groups to  its framework in order to
  analyze global behaviour of both rational dynamical systems and
some classes of systems of
non linear partial differential equations.

Tropical geometry is a kind of scale transform between dynamical systems, which
provides with a correspondence between automata and real rational dynamics.
It allows us to study two dynamical systems at the same time, 
whose dynamical natures are very 
different from each other, by taking scaling limits of parameters.

 There has been  the extended development of partial differential equations
by use of  approximations by real rational dynamics.
Combination of these two aspects creates a new connection between geometric group theory 
and partial differential equations. This is the main aspects which we will focus on in this paper.

A discrete group is called an automata group, if it is generated by  actions on the rooted trees,
whose rules are represented by an automaton. From the view point of dynamical scaling limits,
automata can be regarded as {\em frame-dynamics} which play the role of underlying mechanisms 
for rational dynamical systems and systems of partial differential equations.

In this paper we give a construction of dynamical scale transform from automata groups to partial differential 
equations. In particular we develop a method of dynamical scale transform from 
Mealy automata groups 
to first order hyperbolic systems of PDE.

If rational dynamical systems or  systems of PDE contain  automata groups as their underlying dynamics,
then the geometric group  structures will reflect to the dynamical structures of solutions to the systems.
So a basic and general question is whenever we focus on some geometric properties of discrete groups,
whether one can find some  systems with global solutions which 
satisfy structural similarity with these properties.

Theory of automata groups provide plenty of examples which possess quite characteristic 
properties. Here we focus on the  finitely generated 
 {\em infintie torsion groups}
whose existence is a question called the {\em Burnside problem}.
The first example was discovered  by Ajijan-Novikov ([AN]), and then
by Aleshin as an automata group ([Al]). 
In terms of dynamical systems, one can restate its property as 
infinite recursivity
 that any actions have finite orders, and the number of the minimal orders  are infinite.
 In this paper we show that Aleshin's automaton produces the
rational dynamical systems which satisfy {\em infinite quasi-recursivity}.

In this paper we  include basic analysis of the hyperbolic systems 
of first order PDE of $2$-variables,
which arise from Mealy automata. In particular 
it includes existence, uniqueness, energy estimates with explicit estimates on the constants, 
and so on. It allows us to produce solutions to the PDE systems 
which can be applied to the asymptotic comparison theorem which we will describe below.

Even though the topics of this paper include two different fields,
where one is discrete group theory and the other is PDE analysis,
still  the contents are written in a self-contained way.

\vspace{2mm}

Now let us introduce a class of dynamical systems which we treat.
Let $Q$ be a set and $X$ be a space.
A  {\em state dynamics} on $X$  is the dynamical system given by elements of $Q$,
such that
every $q^0 \in Q$ determines an action on  $X$.

 An
 {\em automaton} is consisted by finite rules which  can create quite complicated state dynamics 
 on the sequences of alphabets.
Let $S$ be another set, and consider all the set of infinite sequences:
$$X_S=\{(s_0,s_1, \dots): s_i \in S\}.$$
Let us consider an automaton ${\bf A}$ which is given by  a pair of functions:
\begin{align*}
& \psi: Q \times S^{\alpha+1} \to S, \\
& \phi: Q \times S^{\beta+1} \to Q
\end{align*}
where $\alpha, \beta \geq 0$.
It gives rise to the state dynamics on $X_S$ as follows.
Let us choose any $q\in Q$ and $\bar{s}=(s_0,s_1, \dots) \in X_S$.
Then:
$${ \bf A}_q: X_S \to X_S$$
 ${\bf A}_q (\bar{s})=(s_0',s_1', \dots)$ is determined 
 inductively by:
  $$s_i'= \psi(q_i,s_i, \dots,s_{i+\alpha}), \quad
 q_{i+1} = \phi(q_i, s_i, \dots, s_{i+\beta}) \ \ (q_0=q).$$
Besides the dynamics over $X_S$, 
the change of the state sets play important roles in the hidden dynamics.

 Any sequences
 $\bar{q}^j =(q^0,\dots, q^j) \in X_Q^{j+1}$
 give dynamics by compositions:
 $${\bf A}_{\bar{q}^j} = {\bf A}_{q^j} \circ \dots \circ {\bf A}_{q^0} : X_S \to X_S.$$
It can happen  that  different automata give the same state dynamics.
In such case, the dynamics ${\bf A}_q$ are the same, but  dynamics of change of the state sets 
can be very different.
Such two automata are called {\em equivalent}.

 The above type of the state dynamics contains two 
 important cases of discrete dynamics which appear in tropical geometry:
 
 {\bf (A)} {\em Integrable systems of cell automaton:}
 One arises from the scaling limits of the integrable systems.
Let us just present  one of the typical case, the
 ultra-discrete Lotka Volterra cell automaton, whose 
 dynamics is given by the rule ([HT], [TTMS]):
 $$s_i'= \varphi(s_{i-1}',s_i, s_{i+1})= s_{i-1}' +\max(L_0, s_i) - \max(L_0, s_{i+1})$$ 
Let us rewrite this  by a state dynamics. 
For  $S=Q$ with the initial state $q_0=q$,
consider the state dynamics with $\psi : Q \times S  \to S,
\phi: Q \times S^3 \to Q$ by:
\begin{align*}
& \psi(q,s) =q, 
& \phi(q,s_0, s_1,s_2) = q + \max(L_0, s_1) - \max(L_0, s_2). \\
\end{align*}
One can assign $s_i' = q_i $ with $s_0'=q$, which describes the above automaton.
 
 There are many other cases of the integrable systems given by cell automata of 
 the above types (see [K4]).

{\bf (B)} {\em  Mealy automaton:}
Suppose both $\alpha=\beta=0$ and so: 
$$\psi: Q \times S \to S, \quad \phi: Q \times S \to Q$$
where
$\psi: (q, \quad) : S \cong S$ are isomorphic for all $q \in Q$.
If we identify $X_S$ with the rooted trees, then the state dynamics give the group actions on the trees,
since the actions can be restricted  level-setwisely.
 The groups generated by these states are called the {\em automata groups}
given by the automata $(\psi,\phi)$.

A general state dynamics will not give actions on the trees, but still 
it may be possible to hold isomorphisms ${\bf A}_q : X_S \cong X_S$,
which are identified with the actions on the boundary of the trees.
In $4.A.3$, we give an example:
\begin{lem}
There is a non Mealy automaton which induce the isomorphism actions on the boundary of the trees.
\end{lem}

Automata groups contain several  discrete groups
which are  quite characteristic in geometric group theory. 
Let us list some of known results (see also [Z]):

(1) Automata groups with intermediate growth ([G]),
more generally with fluctuations ([Br]),

(2) Automata groups which are infinite torsion ([AN],[Al]),

(3) Automata groups with non-uniformly exponential growth ([W]),

(4) A criterion of amenability ([BKN]),

(5) Some classificatin of $2$ state automata groups ([GNS]).
\vspace{3mm} \\
Our basic idea of reasearch direction  is to study how such geometric or analytic properties
reflect to the structures of the  rational dynamics and PDE systems, if they contain
such groups as their frame-dynamics. 
 As an applciation,
we will verify existence of  {\em infinite quasi-recursivity} 
for some rational dynamical systems.

Let us explian more details of the structure of this paper.
\vspace{3mm} \\
{\bf Tropical geometry:}
A relative $(\max, +)$-function $\varphi$ 
is a  piecewise linear function equipped with its presentation of the form: 
$$\varphi(\bar{x})=   
 \max(\alpha_1 + \bar{a}_1 \bar{x}, \dots , \alpha_m+ \bar{a}_m \bar{x})
 - \max(\beta_1 + \bar{b}_1 \bar{x}, \dots , \beta_l+ \bar{b}_l \bar{x}).$$
Tropical geometry associates  the parametrized  rational function  (see [Mi]),
which we call a {\em relatively elementary function} :
$$f_t(\bar{z}) \equiv \frac{k_t(\bar{z})}{h_t(\bar{z})}=
\frac{\Sigma_{k=1}^m t^{\alpha_k} \bar{z}^{\bar{a}_k}}
{\Sigma_{k=1}^l t^{\beta_k} \bar{z}^{\bar{b}_k}}.$$
They admit one to one correspondence of their prsentations to $\varphi$, and
 take positive real numbers if the inputs are also positive.

A cruitial property is that  $f_t$ converge to $\varphi$ by letting $t \to \infty$ in some sense.
Among various ways of scaling limits, tropical geometry behaves
quite nicely, which allows us to obtain several uniform estimates
by comparisons of both dynamical behaviours at the same time.

It is quite characteristic of  $(\max, +)$-functions that different presentations can give the same functions.
For example:
 $$\varphi(y,x) =\max(x,-x)-y, \quad \psi(x,y)=\max(\varphi(x,y), -y)$$ 
are the same functions but have different presentations.
The corresponding rational functions are mutually
$f_t(w,z)=w^{-1}(z+z^{-1})$ and $g_t(w,z)=w^{-1}(z+z^{-1}+1)$,
which are different even as functions.
This motivates us to introduce a notion of {\em tropical equivalence}
between such $f_t$ and $g_t$.

Let us take finite sets $S,Q \subset {\mathbb Z}$, and consider an automaton ${\bf A}$ given by:
$$\varphi:Q \times  S^{\alpha+1} \to S,  \quad
\psi: Q\times S^{\beta+1} \to Q.$$
In general these  can be extended 
over the real numbers  by 
 piecewise linear  functions which  
admit the presentations by the relative  $(\max,+)$ functions:
\begin{align*}
& \psi(\bar{r})=   
 \max(\alpha_1 + \bar{a}_1 \bar{r} , \dots, \alpha_{\alpha+2}+ \bar{a}_{\alpha+2} \bar{r})
 - \max(\beta_1 + \bar{b}_1 \bar{r},\dots ,  \beta_{\alpha+2}+ \bar{b}_{\alpha+2} \bar{r}), \\
& \phi(\bar{l})=   
 \max(\gamma_1 + \bar{c}_1 \bar{l},  \dots , \gamma_{\beta+2}+ \bar{c}_{\beta+2} \bar{l})
 - \max(\delta_1 + \bar{d}_1 \bar{l}, \dots , \delta_{\beta+2}+ \bar{d}_{\beta+2} \bar{l}).
 \end{align*}
\quad
\vspace{3mm} \\
{\bf State systems of the rational dynamics:}
  {\em  The state dynamics} with respect to the pair
 $(\psi, \phi)$ 
 are given by  the discrete dynamics 
 inductively defined by the iterations:
$$ x_i^{j+1}= \psi(y_i^j, x_i^j, \dots, x_{i+\alpha}^j), \quad
 y_{i+1}^j= \phi(y_i^j, x_i^j, \dots, x_{i+\beta}^j)$$
where $x_i^0=x_i$ and $y_0^j=y^j$ are the initial values for $i,j \geq 0$.

Let us visualize this dynamics as follows.
$y=y^0$ determines the map:
$$ {\bf A}_y : {\mathbb R}^{\mathbb N} \to {\mathbb R}^{\mathbb N}, \quad
 {\bf A}_y(x_0,x_1, \dots)= (x_0^1,x_1^1,x_2^1, \dots).$$
So by composition, finite  sequences $(y^0,y^1, \dots, y^l)$ give the maps:
\begin{align*}
&  {\bf A}_{(y^0,y^1, \dots, y^l)} =  {\bf A}_{y^l} \circ \dots \circ {\bf A}_{y^0}:
{\mathbb R}^{\mathbb N} \to {\mathbb R}^{\mathbb N}, \\
&  {\bf A}_{(y^0,y^1, \dots, y^j)} (x_0,x_1, \dots)= (x_0^{j+1},x_1^{j+1},x_2^{j+1}, \dots).
\end{align*}

If $\alpha= \beta=0$, then the Meay automaton gives the automata groups, and 
the state dynamics above exactly extend the group actions over the trees.

Let us consider the corresponding   parametrized rational functions
$f_t$ and $g_t$ with respect to $\varphi$ and $\psi$ respectively.

{\em The state system of the rational dynamics} is given by the 
corresponding rational dynamics:
$$ z_i^{j+1}  = f_t(w_i^j, z_i^j, \dots, z_{i+\alpha}^j), \quad
 w_{i+1}^j= g_t(w_i^j,  z_i^j,\dots, z_{i+\beta}^j )$$
   where  the initial values  are given by $z_i^0= z_i, w_0^j = w^j >0$.

Now let 
$(\psi^1, \phi^1)$ and $(\psi^2, \phi^2)$ be paris of $(\max,+)$-functions,
and $(f_t^1,g_t^1)$ and $(f_t^2,g_t^2)$ be the  corresponding relatively elementary functions.
Assume that both pairs are tropically equivalent:
$$ (f_t^1,g_t^1) \sim (f_t^2,g_t^2).$$

For $l=1,2$, 
let $\{w^j(l)\}_{j \geq 0}$ and $\{z_i(l)\}_{i \geq 0}$ be the initial sequences by positive numbers, 
and denote the  solutions by $(z_i^j(l),w_i^j(l) )$
 to the state systems of the rational dynamics:
 \begin{align*}
& z_i^{j+1}(l)  = f^l_t(w_i^j(l), z_i^j(l), \dots, z_{i+\alpha}^j(l)), \\
& w_{i+1}^j(l)= g^l_t(w_i^j(l),  z_i^j(l),\dots, z_{i+\beta}^j(l) ).
\end{align*}
with the  initial values 
$z_i^0(l)=z_i(l)$ and $w_0^j(l)=w^j(l)$.

Let $P_i(c) = \frac{c^{i+1}-1}{c-1}$.
Our basic  tool of the analysis for the rational dynamical systems are given by the following:
\begin{prop} There are constants $M\geq 1$ and $c\geq 0$
independently of the initial values  so that
the asymptotic uniform-estimates hold:
\begin{align*}
\max(  & \frac{z_i^{j+1}(1)}{z_i^{j+1}(2)} , \  \frac{z_i^{j+1}(2)}{z_i^{j+1}(1)}, \
 \frac{w_{i+1}^j(1)}{w_{i+1}^j(2)}, \   \frac{w_{i+1}^j(2)}{w_{i+1}^j(1)}) \\
 & \\
& \leq      M^{2P_{i+j(\gamma+1)}(c)} [ \{(z_i(l), w^j(l)) \}_{l=1}^2] ^{\tilde{c}^{i+1+j(\gamma+1)}} \\
\end{align*}
where $\gamma = \max(\alpha, \beta)$,  $\tilde{c} = \max(c,1)$,
 and
 the initial rates are given by:
 $$[ \{(z_i(l), w^j(l)) \}_{l=1}^2] 
= \sup_{i,j} \max( \frac{z_i(1)}{z_i(2)} , \frac{z_i(2)}{z_i(1)} , \frac{w^j(1)}{w^j(2)},  \frac{w^j(2)}{w^j(1)} ).$$
\end{prop}
In particular if $c< 1$ and if the initial values are the same, then 
one obtains the uniform estimates by some constants $M^C$.
Both $M$ and $c$ can be known  immediately once the presentations of the defining functions
are given.
In fact 
$c$ is  the Lipschitz constant of the pair $(\psi, \phi)$.

For our later purpose of the application to PDE analysis, 
we verify the uniform estimates when the orbits fluctuate
under controll of torpically equivalent functions as below.

Let $(\{z_i^j(l)\}_{i,j}, \{w_i^j(l)\}_{i,j})$ be as above.
\begin{thm}
 Suppose  $(f_t^1, g_t^1)  \sim (f_t^2, g_t^2)$ are pairwisely tropically equivalent. 

If another sequences $\{w_i^j\}_{i,j}$ and $\{z_i^j\}_{i,j}$
 satisfy the dynamical inequalities:
\begin{align*}
& f^1_t(w_i^j, z_i^j, \dots, z_{i+\alpha}^j)
 \leq z_i^{j+1} \leq  f^2_t(w_i^j, z_i^j, \dots, z_{i+\alpha}^j) , \\
& g^1_t(w_i^j,   z_i^j, \dots,z_{i+\beta}^j) \leq w_{i+1}^j \leq g^2_t(w_i^j,   z_i^j, \dots , z_{i+\beta}^j)
\end{align*}
then the uniform estimates hold for $l=1,2$:
$$( \frac{z_i^j(l)}{z_i^{j+1}} )^{ \pm 1} ,   \quad
 (\frac{w_i^j(l)}{w_{i+1}^j})^{\pm 1}   \\
 \leq  \ \  M^{4P_{i+j(\gamma+1)}(c)} 
[\sup_{i,j} \max \{(\frac{z^0_i(l)}{z_i^0} )^{\pm 1} 
 , (\frac{w^j_0(l)}{w^j_0})^{\pm 1}  \}] ^{\tilde{c}^{i+1+j(\gamma+1)}} $$
\end{thm}
We will apply these estimates for
$f_t^1= \frac{1}{N} f_t$ and $ f_t^2= Nf_t$
 to the analysis of  PDE 
 for some integer $N \geq 1$.

Notice that we have removed the condition of monotonicity  of functions 
which were assumed in [K3].
\vspace{3mm} \\
{\bf Automata groups and tropical geometry:}
Let us take finite sets $S,Q \subset {\mathbb Z}$, and consider an automaton ${\bf A}$ given by
the pair $\varphi:Q \times  S^{\alpha+1} \to S $ and
$\psi: Q\times S^{\beta+1} \to Q$ with their 
extensions  by  relative  $(\max,+)$ functions $( \psi,\phi)$. 

We say that the extensions are {\em stable} over $(Q,S)$, if there are some $0 < \delta <1$ and 
$0 \leq \mu <1$ so that  the Lipschitz constants of the pair $( \psi,\phi)$ 
is bounded  by $\mu $ on $\delta$ neighbourhoods   
of $Q, S \subset  {\mathbb R} $ respectively.
\begin{lem}
For Mealy automata, stable extensions always exist.
\end{lem}
It would be possible to construct stable extensiopns for general automata,
by performing little bit more complicated constructions.

If two automata are equivalent, then the corresponding  state systems of the 
rational dynamics also show mutual structural similarity on the large scale as:  

\begin{thm} Suppose ${\bf A}_1$ and ${\bf A}_2$ are equivalent,
and choose stable extensions of them.

Then for any $C \geq 1$, there exists $t_0 >1$ so that   for all  $t \geq t_0$
and any  initial values 
  $\bar{z}(l) $ and $ \bar{w}(l)$ with the bounds:
$$C^{-1} k_i \leq z_i(l) \leq Ck_i, \quad   C^{-1} q^j \leq w^j(l) \leq C q^j $$ 
for  $\bar{k} \in X_S$ and $\bar{q} \in X_Q$, then 
 the uniform estimates hold:
$$\max \{  \frac{z_i^j(1)}{z_i^j(2)} , \  \frac{z_i^j(2)}{z_i^j(1)} \}
 <  C^4.$$
 \end{thm}

There are many cases when automata groups are finite (see [GNS]).
In such case, $({\bf A}_{\bar{q}^m})^p=\text{ id }: X_S\to X_S$ hold for some $p$ and
 all $\bar{q}^m =(q^0, \dots, q^m) \in X_Q^{m+1}$.

This situation is restated by the state dynamics.
Let ${\bf A}$ be an automaton (not necessarily Mealy), and
consider the dynamics
${\bf A}_{\bar{q}^m}: X_S \to X_S$.
Let us denote the periodic sequence of $\bar{q}^m$ by:
$$\bar{q}^m_{per} \equiv (q^0, \dots, q^m, q^0, \dots ,q^m, \dots,  q^0, \dots, q^m, \dots) \in X_Q.$$

Let us choose a stable extension of ${ \bf A}$ and 
 consider the state dynamics:
 $$ z_i^{j+1}= f_t(w_i^j, z_i^j, \dots, z_{i+\alpha}), \quad
 w^j_{i+1}= g_t(w_i^j   ,z_i^j, \dots, z_{i+\beta}^j)$$ 
with the  initial values
  $C^{-1} t^{k_i} \leq z_i \leq Ct^{k_i}  , \quad  
 C^{-1} t^{q^j}\leq w^j \leq Ct^{q^j} $,
where $ \bar{k} =(k_0,k_1, \dots ) $ are some elements in $X_S$ and
$\bar{q}^m_{per} = (q^0,q^1, \dots) $.

\begin{prop}
Under  the above situation, suppose 
${\bf A}_{\bar{q}^m}: X_S \to X_S$ is of finite order with period $p $.

Then 
for any $C >0$, there exists $t_0>1$ and $D$ independent of the initial values
so that
 the uniform bounds:
$$(   \frac{z_i^j}{z_i^{j+p(m+1)l}} )^{\pm 1} \ \leq    D$$ 
hold for all $t \geq t_0$ and
 all $i,j,l =0,1,2, \dots$
\end{prop}

The automata group by Aleshin is generated by $2$ states and  {\em infinite torsion},
which gives  a solution to the Burnside problem.  
Let us transform such phenomena to the corresponding state systems of the rational dynamics,
which we call the {\em rational Burnside problem}:
\begin{df}
The state system of the parametrized rational dynamics by $(f_t,g_t)$
is quasi-recursive with respect to $(X_t,Y_t)$,
if  for any  $C, C'  \geq 1$, there exists $t_0 >1$  
so that for all $t \geq t_0$
and  any $\{w^j_0\}_j \in Y_t$, there exist some $p\in {\mathbb N}$
such  that:

(1)
any solutions $(\{z_i^j\}_{i,j}, \{w_i^j\}_{i,j})$ with 
 $\{z_i^0\}_i \in X_t$ 
 satisfy the uniform bounds:
$$ (\frac{z_i^{j+pl}}{z_i^j} )^{\pm 1} \leq C$$
  for all  $i,j, l =0,1,2, \dots  $, and 

(2) for any $1 \leq p' \leq p-1 $,
there are some  $\{z_i^0\}_i \in X_t$ so that the solutions 
 $(\{z_i^j\}_{i,j}, \{w_i^j\}_{i,j})$ satisfy the uniform lower bounds:
$$ (\frac{z_i^{j+p'}}{z_i^j} )^{\pm 1} \geq C' $$
 for all  $ j=0,1,2,  \dots$
 and  some  $ i$.

It is  infinitely quasi-recursive, if  infinitely many such $p$ exist.
 \end{df}

By use of stable extensions of the Aleshin's automaton, we have the following:
\begin{thm}
There exists a pair of relatively elementary functions
$(f_t,g_t)$ so that the state systems of the rational dynamics is infinitely quasi-recursive.
\end{thm}
\quad
\vspace{3mm} \\
{\bf Approximations of systems of PDE by rational dynamics:}
It would be of interest to study how globally analytic properties of automata group actions 
on trees effect on the associated dynamics of the PDE systems.

Let us state a general procedure to induce PDE systems from automata.
Let $0< \epsilon \leq 1$ be    constants.
Let us 
consider a function
$u: [0, \infty) \times [0, \infty) \to {\mathbb R}$ of class  $C^{\mu+1}$,
 and  take the Taylor expansions
up to order $\mu+1$:
 \begin{align*}
  u & (x  + i \epsilon   , s+ j \epsilon)
   = u+i  \epsilon u_x +  j \epsilon u_s+
    \frac{(i\epsilon)^2}{2} u_{2x} +  \frac{(j \epsilon)^2}{2} u_{2s}  +i j \epsilon^2  u_{xs} \\
& +
    \dots   + \frac{(j \epsilon)^{\mu}}{\mu !} u_{\mu s} 
  + \frac{(i\epsilon)^{(\mu +1)}}{(\mu +1)!} u_{(\mu +1)x}(\xi_{ij}) +
  \dots
 + \frac{(j\epsilon)^{(\mu +1)}}{(\mu +1)!} u_{(\mu +1)s}(\xi_{ij}) .\\
 \end{align*}

Let us describe a general way to approximate 
solutions to the systems of PDE by the state systems of rational dynamics.
 
 Let $f_t  =  \frac{a_t}{b_t}$ and $ g_t  =  \frac{c_t}{d_t}$ be 
  relatively elementary functions, and
 consider the state systems of the rational dynamics:
$$ z_i^{j+1}= f_t(w_i^j, z_i^j, \dots, z_{i+\alpha}), \quad
 w^j_{i+1}= g_t(w_i^j   ,z_i^j, \dots, z_{i+\beta}^j).$$

Let us  introduce the change of variables by:
$$ i = \frac{x}{\epsilon}, \quad j= \frac{s}{\epsilon}, \quad u(x,s) = z_i^j
\quad v(x,s) =w_i^j$$
Then we take  the difference, and   insert the Taylor expansions:
\begin{align*}
& z_i^{j+1}- f_t(w_i^j,z_i^j, \dots, z_{i+\alpha}^j) 
  =  u(  x,s+\epsilon) 
-   f_t(  v(x,s),   u(x,s), \dots,  u(x+\alpha \epsilon,s)) \\
& = 
\frac{ P_1(\epsilon,t, u,v, u_s, u_x ,\dots, u_{\mu x}) + R_1(\epsilon, t, u,v, \dots, u_{(\mu+1)x}(\xi))}
{ b_t(  v(x,s),    u(x,s), \dots, u(x+\alpha \epsilon,s))}  \\ 
& \equiv  {\bf L}_1(\epsilon, t,u,v,u_s, \dots, u_{\mu x}) 
 + \epsilon^{\mu+1} {\bf E}_1(\epsilon, t,u,v, \dots,
\{ u_{\bar{a}}(\xi_{ij}) \}_{\bar{a},i,j} ) .
\end{align*}
where $P_1$ and $R_1$ are polynomials, and 
each monomial in $R_1$ contains derivatives of $u$ of order $\mu+1$.

Similarly we have the expansions:
\begin{align*}
& w_{i+1}^j- g_t(w_i^j,z_i^j, \dots, z_{i+\beta}^j) 
  =  v(  x+\epsilon,s) 
-   g_t( v(x,s),   u(x,s), \dots,  u(x+\beta \epsilon,s)) \\
& =
\frac{ P_2(\epsilon,t, u,v,u_x, v_x, \dots, u_{\mu x}) +  R_2(\epsilon, t, u,v, \dots, u_{(\mu+1)x}(\xi'))}
{ d_t( v(x,s),    u(x,s), \dots, u(x+ \beta \epsilon,s))}  \\ 
& \equiv  {\bf L}_2(\epsilon, t,u,v,u_x, v_x, \dots, u_{\mu x}) \\
& \qquad \qquad + \epsilon^{\mu+1} {\bf E}_2(\epsilon, t,u,v,  \dots,
\{ u_{\bar{a}}(\xi'_{ij}) \}_{\bar{a},i,j} , \{ v_{\bar{a}}(\xi'_{ij}) \}_{\bar{a},i,j} ) 
\end{align*}
where each monomial in $R_1$ contains derivatives of $u$ or  $v$ of order $\mu+1$.

We say that ${\bf L}_i$ and ${\bf E}_i$ are  the {\em leading} and {\em error} terms  respectively, 
and call the  parametrized systems of PDE of  order $\mu$:
 $$ P_1(\epsilon, t,  u, v,u_s,  \dots, u_{\mu s}, u_{\mu x})=0 ,\quad
 P_2(\epsilon, t,  u,v, v_x, \dots, u_{\mu x}, v_{\mu x})=0 $$
   the {\em induced  systems of partial differential equations}
with respect to  $(\psi, \phi)$.

Let us take four relatively elementary functions $f^1,f^2$ and $g^1,g^2$ so that
$f^1 \sim f^2$ and $g^1 \sim g^2$ are tropically equivalent mutually.
For $l=1,2$, let:
 $$ P_1^l(\epsilon, t,  u, v,u_s,  \dots, u_{\mu s}, u_{\mu x})=0 ,\quad
 P_2^l(\epsilon, t,  u,v, v_x, \dots, u_{\mu x}, v_{\mu x})=0 $$
be the induced  systems of PDEs of order $\mu$.

   For functions of $C^{\mu+1}$ class, let us introduce 
the  {\em higher distorsion} which are the
relatively uniform norms of their $\mu+1$-differentiations  ($3.B$):
$$ 
K(u,v) \equiv  \sup_{(x,s) \in [0, \infty)^2 }   
\max[\frac{||(u,v)||_{ \mu, \alpha}^1}{  u(x,s+ \epsilon)}, \ \
\frac{||(u,v)||_{ \mu, \beta}^2}{  v(x+ \epsilon,s)}].$$

\begin{thm}  For $l=1,2$, let $C$ be the bigger one of their  error constants. 

Let  $(u^l,v^l): [0, \infty) \times [0, \infty) \to (0, \infty)$
be positive solutions to  the above systems 
respectively, so that the estimates:
 $$0  \leq   CK(u^l,v^l) \leq (1-  \delta) \epsilon^{-1}$$
are satisfied for some positive $ \delta >0$.

 Then
 they satisfy the asymptotic estimates
 for all $(x,s) \in [0,  \infty) \times [0, \infty)$:
 \begin{align*}
 (\frac{u^1}{u^2})^{\pm 1}(x,s), & \  \ (\frac{v^1}{v^2})^{\pm 1}(x,s) \\
& \leq (N_0 M)^{6P_{\epsilon^{-1}( x+s (\gamma+1))}(c)}  \ \
 ([ (u^1,v^1):(u^2,v^2)]_{\epsilon})^{\tilde{c}^{\epsilon^{-1}(x+s(\gamma+1))+1}}  
 \end{align*}
where $N_0$  is any integer  with $N_0 \geq  \max(\delta^{-1}, 2-  \delta)$.
 \end{thm}
We say that  a solution to the above PDE systems is {\em admissible},
if the higher distorsions satisfy the above estimates.

\vspace{2mm}

In practice, in order to apply the above result,
there are two types of questions where:

\vspace{2mm}

(1) existence of positive solutions, \ \
(2)
 admissibility of such solutions.
\vspace{2mm} \\
Among all the induced PDE systems,
 the hyperbolic Mealy systems which we treat below
are the one which arise from automata groups.
For such class, existence follows always and we construct the admissible systems of 
PDE quite concretely.
\vspace{3mm} \\
{\bf Hyperbolic systems of PDE in tropical geometry:}
In the case of   Mealy dynamics given by:
$$ x(i,j+1)= \psi(y(i,j), x(i,j)), \quad 
 y(i+1,j)= \phi(y(i,j), x(i,j))$$
 the induced  first order systems of the equations:
$$\epsilon \ u_s = f_t(v,u) -u, \quad
   \epsilon \ v_s  =g_t(v,u) -v$$
are called  the {\em hyperbolic Mealy systems}.

In this paper we develop  basic analysis of  the hyperbolic Mealy systems.
In particular we verify existence of positive solutions, uniqueness and 
explicit energy estimates. 
 Construction of  admissible solutions to the 
hyperbolic Mealy systems 
involve the interplay of estimates between piecewise linear and differentiable dynamics.

Let us treat a case when  the pair $(f_t,g_t)$ restricts to a self-dynamics over $[r,R]$ ($5.C$).
Let us put the initial domain:
$$I_0 = [0, \infty) \times \{0\} \  \cup   \ \{0\} \times [0, \infty) .$$

\begin{thm}
Suppose that  the pair $(f_t, g_t)$ restricts to a self-dynamics over $[r,R]$, 
and give the positive initial values:
$$u , \ v :   \ I_0  \to [r+q, R- q].$$
Then:
(1)  there exists a positive solution  
 $$u,v : [0, \infty) \times [0, \infty) \to (0, \infty)$$
 with the uniform bounds 
 $r +q \leq u(x,s), \ v(x,s) \leq R-q$.
  
 (2) The solution is unique.
 \end{thm}

Let us induce the energy estimates which are well known for the hyperbolic systems, but 
here we also estimate the constants explicitly.
Let us  introduce the numbers:
 \begin{align*}
& D=  \sup_{(u,v)   \in [r,R]^2}  \{ \   |f_t(v,u)-u|,   \ |g_t(v,u)-v| \  \},  \\
& B = \max(  ||(f_t)_u -1||C^0 , ||(f_t)_v||C^0, \ ||(g_t)_v-1||C^0 , ||(g_t)_u||C^0) , \\
& b =  \sup | (f_t)_v (g_t-v)|, \quad d=  \sup |(g_t)_u(f_t-u)|.
 \end{align*}

\begin{prop}
Suppose that  the pair $(f_t, g_t)$ restricts to a self-dynamics over $[r,R]$.

(1) Let us give the  initial values:
$$u(\quad ,0) , \ v(0, \quad) : [0,  \infty)       \to [r+q, R- q]$$
with uniformly bounded $C^1$norms:
$$||u_x||C^0([0,  \infty)  \times  \{0 \}) ,  \quad  ||v_s||C^0(   \{0 \} \times [0,  \infty)) \leq A  <  \infty.$$

Then  there is a constant $C$ so that solutions 
 $u,v : [0, \infty) \times [0, \infty) \to (0, \infty)$
 have the  asymptotic $C^1$ bounds:
  \begin{align*}
  &||\frac{\partial u}{\partial x}||C^0( [0, \infty)\times \{m \}) , \  \
  ||\frac{\partial v}{\partial s}||C^1(  \{m\} \times [0, \infty) ) \  \leq \ 2^{\tau^{-1}m } (A  + 2D), \\
&   ||\frac{\partial u}{\partial s}||C^0 , \quad  ||\frac{\partial v}{\partial x}||C^0 \ \ \leq  \ \ D
\end{align*}
where
$\quad  \tau \text{ Lip}_{\bar{f}_t,  \bar{g}_t} \leq \frac{1}{2}, 
\quad \tau \leq D^{-1}q , \quad \delta \equiv \tau B \leq \frac{1}{4}$.

\vspace{3mm}

(2)
Assume negativities:
$$-a   \ \ \leq (f_t)_u-1 , \quad   (g_t)_v -1 \leq  \  \ -c$$
for some $0 <a,b$. 
Then the uniform estimates hold:
\begin{align*}
& \frac{b}{a} + (u_x(x,0) -\frac{b}{a})\exp(-as) \leq u_x(x,s) \leq
\frac{b}{c} + (u_x(x,0) -\frac{b}{c})\exp(-as),  \\
& \frac{d}{a} + (v_s(0,s) -\frac{d}{a})\exp(-ax) \leq v_s(x,s) \leq
\frac{d}{c} + (v_s(0,s) -\frac{d}{c})\exp(-as)
\end{align*}

In particular $|u_x|$ and $|v_s|$ are both uniformly bounded.
\end{prop}

In $5.D$, we construct admissible systems of PDE from piecewisely linear functions.
For this purpose we introduce some modification of automata.

\begin{df} $(\bar{\psi}, \bar{\phi})$ is an
$\epsilon$-refinement of the pair $(\psi, \phi)$, if
 there is   a posiive number $N$ 
so that for any $\bar{q}\in X_Q$
and $\bar{k}\in X_S$,
there are paths $y: \{0, 1, \dots\} \to {\mathbb R}$
and $x: \{0,1, \dots \} \to {\mathbb R}$
with:
\begin{align*}
& y(jN)= q^j , \quad |y(j+1)-y(j)|\leq \epsilon, \\
& x(iN) =k_i, \quad |x(i+1) - x(i)| \leq \epsilon
\end{align*}
 for all $i , j \in \{ 0,1, \dots \}$, such that
the corresponding orbits:
$$(\{x_i^j\}, \{y_i^j\}), \quad (\{s_i^j\}, \{q_i^j\})$$
with respect to $(\bar{\psi}, \bar{\phi})$ and $(\psi, \phi)$ respectively, satisfy 
 the equalities:
$$x_{iN}^{jN}= k_i^j, \quad   y_{iN}^{jN} = q_i^j.$$
$(\psi, \phi)$ is refinable, if there is an $\epsilon$-refinement for any small $\epsilon >0$.
\end{df}

Let ${\bf A}$ be a Mealy automaton with $2$ alphabets, 
equipped with a representative $(\psi, \phi)$  by relatively $(\max,+)$-functions. 
For any $(s_0,s_1,\dots) \in X_2$ and $(q^0,q^1, \dots) \in X_{m+1}$,
let:
$${\bf A}_{(q^0, \dots,q^{l-1})} (s_0,s_1, \dots) = ((s_0^l, s_1^l, \dots ) \in X_2$$
be the orbits of the automata group actions.

For its refinement $(\bar{\psi}, \bar{\phi})$ and 
their tropical correspondences $( \bar{f}_t, \bar{g}_t)$, let us consider the hyperbolic Mealy systems:
$$u_s=  \bar{f}_t(v,u)-u, \quad v_x=  \bar{g}_t(v,u)-v.$$

\begin{thm}   For any $C>0$ and any $t \geq t(C) >1$, 
there are refinements 
 $(\bar{\psi}, \bar{\phi})$ of $( \psi,  \phi)$ with the   pairs of
 tropical correspondences $(\bar{f}_t,\bar{g}_t)$
so that:

(1) $(\bar{f}_t,\bar{g}_t)$ admits admissible solutions, 

(2)
for any another pairs $(f_t,g_t)$ toropically equivalent to 
$(\bar{f}_t,\bar{g}_t)$, 
 any admissible solutions to the equations:
$$u_s=f_t(v,u)-u, \quad v_x= g_t(v,u)-v$$
whose  initial values satisfy the inclusions for all $k =0,1,2, \dots$:
$$d(u(N k, 0) , S) ,\quad d(v(0,Nk ), Q) \leq C$$
then  they satisfy the asymptotic estimates  for some $M, c \geq 1$:
 \begin{align*}
 (\frac{u(N i, N j) }{s_i^j})^{\pm 1} & \  \ \leq M^{ P_{N(i+j) +1 } (c)}
 \end{align*}
 \end{thm}
This is a consequence of admissibility of the refinements as:
\begin{prop}
Let ${\bf A}$ be a Mealy automaton with $2$ alphabets.
There is an refinement of ${\bf A}$ with the pair of functions $(\bar{\phi}, \bar{\psi})$
so that the corresponding relatively elementary functions $(\tilde{f}_t, \tilde{g}_t)$
are admissible with the estimates:
 \begin{align*}
& [ \ |  ( \tilde{f}_t(v,u) -u) ( (\tilde{f}_t)_u(v,u)-1)|+
  |( \tilde{f}_t)_v(u,v)||v_s| \ ]
(x,s+ \alpha)   <2u(x,s+1),  \\
& [ \ |  ( \tilde{g}_t(v,u) -v) ( (\tilde{g}_t)_v(v,u)-1)|
+   |( \tilde{g}_t)_u(u,v)||u_x| \ ]
(x+ \alpha,s)   <2v(x+1,s).
 \end{align*}
for any solutions $(u,v)$ and all $0  \leq \alpha  \leq 1$.
 \end{prop}

It would be of interest to apply group-theoretic results to the analysis of PDE.
Compared with the dynamical Burnside problem,
we would like to propose the following.

Let $X \subset [0, \infty)$ be a net and   $Y \subset [0, \infty)$ be a periodic subset.
Let us take parametrized functional subsets
 $A_X, B_Y  \subset C[0, \infty)$.
  \vspace{3mm} \\
{ \em Conjecture 1.1:} There exists a Meay automaton  so that the corresponding hyperbolic system 
is  infinitely quasi-recursive  over $(A_X , B_Y)$ in the following sense;
 there exists $D\geq 1$ so that:

(1) There are solutions 
$u,v : [0, \infty)^2 \to (0, \infty)$ for  any initial values
 $u(x, 0)\in A_X$ and $ v(0,s) \in  B_Y$.

(2)   For  each $ y \in B_Y$, there is a minimal $p$ such that
any solutions $(u,v)$ with $u(\quad, 0) \in A_X$ and $v(0, \quad) =y$
satisfy  quasi periodicity:
$$(\frac{u_d(x,s)}{u_d(x, s+k p)})^{\pm 1} \leq D$$
for any $(x,s), (x, s+k p)
\in  X \times Y$.

(3) 
 $\lim_{t \to \infty}  [\sup_{(x,s) \in X \times Y} u(x,s) ][\inf_{(x,s) \in X \times Y } u(x,s) ]^{-1}=   \infty$.

  (4) Infinitely many such $p$ exist for each $t>>1$.

 \newpage

\section{Tropical geometry}
{\bf 2.A Tropical transform:}
A {\em relative} $(\max, +)$-{\em function} $\varphi$ 
is a  piecewise linear function equipped with its presentation of the form: 
$$\varphi(\bar{x})=   
 \max(\alpha_1 + \bar{a}_1 \bar{x}, \dots , \alpha_m+ \bar{a}_m \bar{x})
 - \max(\beta_1 + \bar{b}_1 \bar{x}, \dots , \beta_l+ \bar{b}_l \bar{x})$$
 where
 $\bar{a}_l \bar{x}= \Sigma_{i=1}^n a_l^i x_i$, 
 $\bar{x} =(x_1, \dots,x_n) \in {\mathbb R}^n$,  
$ \bar{a}_l =(a_l^1, \dots, a_l^n) , \bar{b}_l \in {\mathbb Q}^n$
 and $\alpha_i , \beta_i  \in {\mathbb R}$. 

We  say that the multiple integer $M \equiv ml$
 is the {\em number of the components}.

 $\varphi$  is Lipschitz  since it  is piecewise linear.
This is the most important property for our analytic estimates later.

 Throughout this paper, we equip the metric on ${\mathbb R}^n$ by:
 $$|(x_1, \dots, x_n)|= \sup_{1 \leq i \leq n} |x_i|.$$

Corresponding to $\varphi$, 
tropical geometry associates  the parametrized  rational function given by:
$$f_t(\bar{z}) \equiv \frac{k_t(\bar{z})}{h_t(\bar{z})}=
\frac{\Sigma_{k=1}^m t^{\alpha_k} \bar{z}^{\bar{a}_k}}
{\Sigma_{k=1}^l t^{\beta_k} \bar{z}^{\bar{b}_k}}$$
 where
$  \bar{z}^{\bar{a}}= \Pi_{i=1}^n z_i^{a^i}$, 
$\bar{z}=(z_1, \dots, z_n) \in {\mathbb R}^n_{>0}$  ([LM], [V], [Mi]).

We say that $f_t$  is a {\em relatively elementary function}, and
both terms $h_t(\bar{z})= \Sigma_{k=1}^l t^{\beta_k} \bar{z}^{\bar{b}_k}$ 
and $k_t(\bar{z})= \Sigma_{k=1}^m t^{\alpha_k} \bar{z}^{\bar{a}_k}$ are
just elementary functions.

These two functions $\varphi$ and $f_t$
admit one to one correspondence between their presentations.
In fact they  are connected passing through some intermediate
functions $\varphi_t$, which
Maslov inroduced as  dequantization of the real line ${\mathbb R}$.

Let us briefly explain the aspects of scaling limit in tropical geometry.
For $t>1$, there is a family of semi-rings $R_t $ which are all 
the real number ${\mathbb R}$ as sets. 
The multiplications and the additions are respectively  given by:
$$x \oplus_t y = \log_t (t^x + t^y), \quad  x \otimes_t y = x+y.$$
As $t \to \infty$ one obtains  the equality:
$$x \oplus_{\infty} y = \max (x,y).$$
By use of these semi-ring structure, one has
$R_t$-{\em polynomials} of the form:
$$\varphi_t(\bar{x})=   
(\alpha_1 + \bar{a}_1 \bar{x}) \oplus_t \dots \oplus_t ( \alpha_m+ \bar{a}_m \bar{x})
 - (\beta_1 + \bar{b}_1 \bar{x}) \oplus_t \dots  \oplus_t (\beta_l+ \bar{b}_l \bar{x})$$
\quad
\vspace{3mm} \\
{\bf 2.A.2 Basic properties:}
So far we have seen three different types of functions,
$\varphi, f_t$ and $\varphi_t$.
Let us list some basic properties they satisfy.

\vspace{2mm}

({\bf A})
$\varphi$ and $\varphi_t$ are connected as:
 $$\lim_{t \to \infty} \varphi_t =\varphi$$
 and the limit
  satisfies the relative  $(\max ,+)$ equation
 $ \varphi_{\infty}(\bar{x}) = \varphi(\bar{x})$ as in $2.A$.
 In fact we have the uniform estimates as below.
 Let $M$ be the number of the components for $\varphi$.

\begin{lem}[K2]
The uniform estimates hold:
$$\sup_{\bar{x}  \in {\mathbb R}^n} |\varphi_t(\bar{x}) - \varphi(\bar{x})| \leq \log_t M.$$
\end{lem}
{\em Proof:}
For convenience we include the proof for the simple case.
We verify the estimate 
$|x \oplus_t y  - \max(x, y)| \leq \log_t 2 $.
Assume  $x =  \max(x, y)$. 
Then: 
$$x \oplus_t y  = \log_t (t^x+ t^y) 
 = \log_t (t^x(1+ t^{y-x}  )) 
= x+ \log_t( 1+ t^{y-x}).$$
Since $y-x \leq 0$ are  non positive, 
 the estimates $\log_t(1+ t^{y-x} ) \leq \log_t 2 $ hold.
The general case easily follows from this.
This completes the proof.
\vspace{3mm} \\
{\em Remark 2.1:}
At a glance tropical geometry seems a special kind  among various  scale transforms.
However we would point out the following fact, which suggests 
that  tropical geometry posesses some universality 
from the arithmetic view points:
\begin{lem}
Let $f: { \mathbb R} \to  {\mathbb R}  $ be a continuous map 
with the property:
$$f(x+y) =f(x)f(y), \qquad x,y \in  {\mathbb R} .$$
 Then there is some $t >0$ so that 
one of $f(x)=t^x$ or $f(x)=0$ holds. 
\end{lem}
{\em Proof:}
Firstly $f(x) \geq 0$ hold, since $f(x) = f(\frac{x}{2})^2$.

$f(a)= (f(1))^a$ hold for any rational $a= \frac{m}{n} \in {\mathbb Q}_{>0}$.
This follows from the equalities  $f(1)= f(n \frac{1}{n})= f(\frac{1}{n})^n$ and
$f(a)= f(\frac{1}{n})^m = f(1)^{\frac{m}{n}}=f(1)^a$.
By continuity, the same formula holds for any $a \in {\mathbb R}$.

$f  \equiv 0$ hold iff $f(b)=0$ holds for some $b \in {\mathbb R}  $, since
the equalities  $f(a) = f(a-b)f(b)$ hold.

$f(0)=0$ or $f(0)=1$ must hold by the equality
$f(0)=f(0)f(0)$.
So $f\equiv 0$ hold  iff $f(0)=0$.

Conversely suppose $f(0)=1$ holds. Then
$f(a) > 0$ hold for any $a  \in {\mathbb R}$, and so
 $f(a)=t^a$ hold  with $t \equiv f(1)>0$.
 This completes the proof.

\vspace{3mm}
 
({\bf B})
Three  types of the  functions equipped with their presentations:
$$\varphi, \quad \varphi_t, \quad f_t$$ 
have one to one correspondences with each other.
Namely one can obtain the presentations of all these functions 
at the same time,
once  the coefficients 
$\alpha_i, \bar{a}_i, \beta_j, \bar{b}_j$ are determined.

Notice that  as functions, 
the  presentations of relative  $(\max, +)$ functions
$\varphi$ 
are not uniquely determined  in general, unlike to the case of rational functions.

For example $\varphi(x,y)= \max(-y,y)-x$ and $\psi(x,y)=\max(\varphi(x,y),-x)$
 are the same functions, but have the different presentations.
Correspondingly the rational functions have the presentations as: 
$$f(z,w)=z^{-1}(w^{-1}+w), \quad g(z,w)= z^{-1}(w^{-1}+w+1)$$
which are mutually different even as functions.

This leads us to  the following notion.
 Let $\varphi^1$ and $\varphi^2$ be two relative $(\max, +)$-functions with $n$ variables.
Then $\varphi^2$ is {\em equivalent} to $\varphi^1$ ($\varphi^1 \sim \varphi^2$), 
if they are the same as functions, so
$\varphi^1(x_1, \dots, x_n)=\varphi^2(x_1, \dots,x_n)$ hold for all $(x_1, \dots,x_n) \in {\mathbb R}^n$
but possibly they may have different presentations.
\begin{df}[K2]
Let $f^1_t$ and $f^2_t$ be  two relatively elementary functions.
$f_t^1$ and  $f^2_t$ are mutually  tropically equivalent, 
if the corresponding  $(\max, +)$-functions $\varphi^1$ and  $\varphi^2$ are equivalent.
\end{df}

For example if $f_t$ corresponds to $\varphi$,
 then $3f_t$ corresponds to $\max(\varphi,\varphi,\varphi)$, and so on.
 So $f_t$ and $N f_t$ are tropically equivalent for all $N=1,2,\dots$
 More generally $f_t$ and $a f_t$ are tropically equivalent for any $a  \in {\mathbb Q}_{>0}$.

\vspace{3mm}

For $N,M  \geq 1$ and $f_t$, let us put:
$$\tilde{f}_t(z_0, \dots,z_{n-1}) \equiv f_t(\frac{M}{N}z_0, \dots, \frac{M}{N}z_{n-1}).$$

\begin{lem}
$f_t$ and $\tilde{f}_t$ are tropically equivalent.
\end{lem}
{\em Proof:}
The function  $z \to \frac{M}{N}z$ corresponds tropically  to: 
$$x \to \mu(x) \equiv  \max(x, \dots, x)  - \max(0,\dots,0) =x.$$

Let $\varphi$ be the corresponding relative $(\max, +)$-function to $f_t$.
Then the another relative  $(\max,+)$-function:
$$\tilde{\varphi}(x_0, \dots,x_{n-1})\equiv 
\varphi( \mu(x_0), \dots, \mu(x_{n-1}))=
\varphi(x_0, \dots, x_{n-1})$$
corresponds to $\tilde{f}$.
This completes the proof.
\vspace{3mm}

\vspace{3mm}

({\bf C})
Let us relate $\varphi_t$ with $f_t$.
Let $\text{Log}_t: {\mathbb R}^n_{>0} \to {\mathbb R}^n$ be given by:
$$ (x_0, \dots, x_{n-1}) =\text{Log}_t(z_0, \dots, z_{n-1})\equiv    (\log_t z_0, \dots, \log_t z_{n-1}).$$

\begin{prop}[LM,V]
$f_t \equiv (\log_t)^{-1} \circ \varphi_t \circ \text{Log}_t : {\mathbb R}_{>0}^n \to (0, \infty)$
is a parametrized rational function
$f_t(\bar{z}) \equiv \frac{k_t(\bar{z})}{h_t(\bar{z})}=
\frac{\Sigma_{k=1}^m t^{\alpha_k} \bar{z}^{\bar{a}_k}}
{\Sigma_{k=1}^l t^{\beta_k} \bar{z}^{\bar{b}_k}}$.
\end{prop}
One can check this equality by  direct calculations.
Even though  verification is quite easy,  it plays an important role
in tropical geometry.

In particular $f_t>0$ take positive values on $ {\mathbb R}_{>0}^n=(0, \infty)^n$.
Notice that $f_t(\bar{0})=0$ may occur.

\vspace{3mm} 

Let $(\varphi, \varphi_t, f_t)$ be the triplet as above,
and $M$ be the number of the components.
Let us induce $C^0$ comparisons:
\begin{lem}
Suppose $\varphi$ is uniformly  bounded from both  above and below
as $a \leq \varphi(\bar{x})  \leq b$.
Then  $f_t$ admits uniform bounds:
$$ t^a M^{-1} \leq  f_t(\bar{z}) \
\leq t^b M. $$
\end{lem}
{\em Proof:}
For any $\bar{z} \in {\mathbb R}_{>0}^n$, let us denote $\bar{x}=\text{Log}_t(\bar{z})$.
By use of proposition $2.4$, 
we have the equalities:
$$ f_t(\bar{z}) = t^{\varphi_t(\text{Log}_t(\bar{z})))} \
= t^{\varphi_t(\bar{x})} \
 = t^{\varphi(\bar{x})} \ t^{\varphi_t   (\bar{x}) - \varphi(\bar{x})}.$$

 Then By  lemma $2.1$ we have the estimates:
$$ t^a M^{-1} \leq  f_t(\bar{z}) \
\leq t^b M. $$
This completes the proof.
\vspace{3mm} 

Next let us induce $C^1$ comparisons.
Let $\psi$ and $\phi$ be two $(\max, +)$-functions
so that the equality $\varphi = \psi - \phi$ holds.

\begin{lem}
Suppose (1)  $\varphi$ is uniformly  bounded  from both above and below,
and (2) $\phi$ is bounded from below.
Then 
any derivatives of $f_t$ are  uniformly bounded  from above,
after change   by troipcally  equivalent one, if necessarily.
\end{lem}
{\em Proof:}
The  bounds
$a \leq  \varphi = \psi - \phi \leq  b$ hold
by  the assumption.

Let us write $f_t =\frac{h_t}{k_t}$, where both $h_t$ and $k_t$ are polynomials
parametrized by $t>1$, which correspond to  $\psi$ and $\phi$ respectively.

By the assumption, there is some $N \geq 0$ so that 
the equalitiy $\phi' = \max( -N, \phi)$ holds as functions.
Then 
 we may assume  positivity 
 $k_t(0) > 0$, by change of $\phi$ by $\phi'$ if necessarily,
which corresponds to $k_t+ t^{-N}$.

The uniform bounds of the range of $f_t$ implies that 
degrees of $h_t$ and $k_t$ coincide each other.
Let us consider  the derivative $f_t' = \frac{h_t'}{k_t} - \frac{h_t k_t'}{k_t^2}$.
The both terms have the property that the  degree of the denominators 
are strictly larger than that of the numerators. 

These  imply that the derivatives of $f_t$ are also
uniformly bounded from above.
This completes the proof.
\vspace{3mm}

Notice that 
the higher derivateives can be considered similarly.
In section $5$, we induce more detailed estimates.
\vspace{3mm}\\
{\bf 2.B State dynamics:}
In $2.B$ we study analysis of dynamical systems by relatively elementary functions.
Let:
  $$\varphi: {\mathbb R}^{\alpha+2} \to {\mathbb R}, \quad
 \psi: {\mathbb R}^{\beta+2} \to {\mathbb R}$$
be two piecewise-linear  functions which  
admit their presentations by the  relative  $(\max,+)$ functions:
\begin{align*}
& \psi(\bar{r})=   
 \max(\alpha_1 + \bar{a}_1 \bar{r} , \dots, \alpha_{\alpha+2}+ \bar{a}_{\alpha+2} \bar{r})
 - \max(\beta_1 + \bar{b}_1 \bar{r},\dots ,  \beta_{\alpha+2}+ \bar{b}_{\alpha+2} \bar{r}), \\
& \phi(\bar{l})=   
 \max(\gamma_1 + \bar{c}_1 \bar{l},  \dots , \gamma_{\beta+2}+ \bar{c}_{\beta+2} \bar{l})
 - \max(\delta_1 + \bar{d}_1 \bar{l}, \dots , \delta_{\beta+2}+ \bar{d}_{\beta+2} \bar{l})
 \end{align*}
where $\alpha, \beta,\gamma, \delta \in {\mathbb R}$,
$\bar{a},  \bar{b} \in {\mathbb R}^{\alpha+2} $ and 
  $\bar{c}, \bar{d} \in {\mathbb R}^{\beta+2}$ 
are all constants.

Later on we denote by $M= \max(M_{\psi}, M_{\phi})$ and $c= \max(c_{\psi},c_{\phi})$ 
as the bigger ones  of the numbers of the components and
the Lipschitz constants for $\psi$ and $\phi$ respectively.

\vspace{3mm}

Let us take initial sequences in ${\mathbb R}$:
 $$\{x_i\}_{i \geq 0}, \qquad \{y^j\}_{j \geq 0}.$$ 
  \begin{df}
 The state dynamics with respect to the pair
 $(\psi, \phi)$ 
 are given by  the discrete dynamics 
 inductively defined by the iterations:
\begin{align*}
& x_i^{j+1}= \psi(y_i^j, x_i^j, \dots, x_{i+\alpha}^j), \\
& y_{i+1}^j= \phi(y_i^j, x_i^j, \dots, x_{i+\beta}^j)
\end{align*}
where the initial values are given by $x_i^0=x_i$ and $y_0^j=y^j$ for $i,j =0,1,2, \dots$
\end{df}
\quad
\vspace{3mm} \\
{\bf 2.B.2 Comparisons between iterated dynamics:}
Let $(\psi_t, \phi_t)$ be the tropical correspondences of the pair $(\psi, \phi)$.
One can also consider another  state dynamics by use of $(\psi_t, \phi_t)$ instead of the pair.
In fact we have uniform estimates of their orbits, which we describe below.

Recall the dynamics of automata in the introduction.
There are dynamics of the states $y_i^j$ behind the  actions $ {\bf A}_{(y^0,y^1, \dots, y^j)}$.
Firstly let us start analyzing the  orbits $\{y_0^j\}_j$, since their treatment
is relatively simple compared with  $\{x_i^j\}_{i,j}$.

Let $\varphi: {\mathbb R}^n \to {\mathbb R}$ be a relative $(\max,+)$-function of $n$ variables,
with the number of the components $M$ and  the Lipschitz constant $c$. Let
$\varphi_t$ be the corresponding function.

Now we consider the dynamics of the states, and
 choose any initial data  $(x_0,x_1, \dots) \in {\mathbb R}^{\mathbb N}$ 
and $y_0 \in {\mathbb R}$. Then one considers two discrete dynamics 
defined inductively by the iterations for $i \geq 0$:
\begin{align*}
& y_{i+1}= \varphi(y_i,  x_i, \dots, x_{i+\beta}), \\
& y_{i+1}'= \varphi_t(y_i', x_i, \dots, x_{i+\beta})
\end{align*}
where $y_0'=y_0$.
We  denote
 $\bar{x}_i =(x_i, \dots, x_{i+\beta})$
for simplicity of the notation. 

Let us put the polynomials of degree $i$:
$$P_i(c)= \frac{c^{i+1}-1}{c-1}.$$ 
Notice the equality:
$$cP_i(c)+1=P_{i+1}(c).$$

\begin{lem}
The uniform estimates:
$$|y_i-y_i'| \leq P_{i-1}(c) \log_t M$$
hold for all $i \geq 0$.
\end{lem}
{\em Proof:}
Firstly  the estimates
$|y_1' - y_1| = |\varphi_t(y_0,\bar{x}_0) -\varphi(y_0,\bar{x}_0)|
 \leq  \log_t M$ hold by lemma $2.1$.

Next we have  the estimates:
\begin{align*}
|y_2-y_2'| & =  |\varphi(y_1,\bar{x}_1) - \varphi_t(y_1',\bar{x}_1)| \\
& \leq  |\varphi(y_1,\bar{x}_1) - \varphi(y_1',\bar{x}_1)| +
 |\varphi(y_1',\bar{x}_1) - \varphi_t(y_1',\bar{x}_1)| \\
& \leq c|y_1-y_1'| + \log_t M
\leq (c+1) \log_t M
\end{align*}
Similarly  we have the following estimates:
 \begin{align*}
 |y_3& -y_3'|  =
  |\varphi(y_2,\bar{x}_2) - \varphi_t(y_2',\bar{x}_2)|  \\
& \leq 
  |\varphi(y_2,\bar{x}_2) - \varphi(y_2',\bar{x}_2)|  +
   |\varphi(y_2',\bar{x}_2) - \varphi_t(y_2',\bar{x}_2)|  
 \leq [c(c+1) +1] \log_t M
 \end{align*}
By iterating the same estimates,
one obtains the conclusion.

This completes the proof.
\vspace{3mm} \\
{\bf 2.C Basic estimates for orbits:}
The estimates for  the dynamics of  $\{x_i^j\}_{i,j}$ involve more complicated analysis.
As preliminaries, we verify some general estimates which will be used later.

For   four  sequences
$\{p_i\}_i, \{q^j\}_{j}, \{x_i\}_{i}, \{y^j\}_{j}$ by real numbers, let us
introduce the numbers:
$$|\{p_i,x_i\}_i; \{q^j,y^j\}_j|  \ \ \equiv  \ \
\sup_{i,j} \{ |p_i-x_i|, \ \ |y^j -q^j|\}.$$
Let us take  four sequences:
$$\{p_i^j\}_{i,j\geq 0}, \{q_i^j\}_{i,j\geq 0}, \{x_i^j\}_{i,j\geq 0}, \{y_i^j\}_{i,j\geq 0}.$$

For the Lipschitz constants $c$, let us put:
$$\tilde{c} = \max(c,1).$$

The following type of the estimates are applied when  we consider   Mealy automata.
The general cases are treated after this version.

\begin{lem}
Suppose these sequences  satisfy the following estimates:
$$ |p_i^{j+1} - x_i^{j+1}| , \ \ |q_{i+1}^j - y_{i+1}^j| 
 \leq c \max( | q_i^j-y_i^j|,  \ |p_i^j-x_i^j|) + T $$
 for some $T \in {\mathbb R}$ and $c \geq 0$.
 
Then they satisfy the estimates for all $i,j \geq 0$:
$$  |p_i^{j+1} - x_i^{j+1}| , \ \    |q_{i+1}^j - y_{i+1}^j|
\leq P_{i+j}(c) T + \tilde{c}^{i+j}|\{p_i^0,x_i^0\}_i ; \{q^j_0, y^j_0\}_j|.$$
\end{lem}

In particular:

(1) If they satisfy the same initial conditions:
$$p_i^0=x_i^0, \ \ q_0^j=y^j_0 \qquad (i,j \geq 0)$$
 then the uniform estimates hold:
 $$  |p_i^{j+1} - x_i^{j+1}| , \ \    |q_{i+1}^j - y_{i+1}^j| \leq P_{i+j}(c) T.$$

(2) If $c <1$ holds, then $  |p_i^{j+1} - x_i^{j+1}| $ and $   |q_{i+1}^j - y_{i+1}^j|$
are both uniformly bounded.

\vspace{3mm}

{\em Proof of lemma $2.10$:}
We verify the conclusion by induction on $i+j \geq 1$.

For $i+j = 1$, the estimates:
\begin{align*}
 |p_0^1 - x_0^1|   , \  |q^0_1-y^0_1| &  \leq c\max(  | q_0^0-y_0^0|,  \ |p_0^0-x_0^0|) + T \\
& \leq T+ c |\{p_i^0,x_i^0\}_i ; \{q^j_0, y^j_0\}_j| .
\end{align*}
hold. So the conclusion follows for $i+j = 1$.

Suppose the conclusion holds for all $(i,j)$ with $i+j \leq N$, 
and take any $(i,j)$ with $i+j=N\geq 1$.
 
Firstly assume both $i,j \ne 0$. Then the estimates hold:
\begin{align*}
  |p_i^{j+1} - x_i^{j+1}|  , \   |q_{i+1}^j  - y_{i+1}^j| 
  &  \leq c \max(  | q_i^j-y_i^j|,  \ |p_i^j-x_i^j|) + T \\
& \leq c( P_{i+j-1}(c)T    +   \tilde{c}^{i+j-1} |\{p_i^0,x_i^0\}_i ; \{q^j_0, y^j_0\}_j|) +T \\
&  =P_{i+j}(c)T+  \tilde{c}^{i+j}|\{p_i^0,x_i^0\}_i ; \{q^j_0, y^j_0\}_j|.
\end{align*}

Next suppose $j=0$ and $i=N \geq1$. Then:
\begin{align*}
 & |p_i^1 - x_i^1|  , \   |q_{i+1}^0  - y_{i+1}^0| 
    \leq c \max(  | q_i^0-y_i^0|,  \ |p_i^0-x_i^0|) + T \\
& \leq c\max( P_{i-1}(c)T    +   \tilde{c}^{i-1} |\{p_i^0,x_i^0\}_i ; \{q^j_0, y^j_0\}_j|, 
|\{p_i^0,x_i^0\}_i ; \{q^j_0, y^j_0\}_j|) +T \\
&  \leq P_{i+j}(c)T+ \tilde{c}^i|\{p_i^0,x_i^0\}_i ; \{q^j_0, y^j_0\}_j|
\end{align*}
since the estimates hold:
$$  |p_i^0-x_i^0|\leq 
|\{p_i^0,x_i^0\}_i ; \{q^j_0, y^j_0\}_j|\leq \tilde{c}^i|\{p_i^0,x_i^0\}_i ; \{q^j_0, y^j_0\}_j|$$
We can treat the case $i=0$ by the same way.
Thus we have verifyed the claim for $i+j\leq N+1$.
This finishes the induction step.

This completes the proof.

\vspace{3mm}

Let us take  four  sequences
$\{p_i\}_i, \{q^j\}_{j}, \{x_i\}_{i}, \{y^j\}_{j}$ as before.
Let $\alpha, \beta \geq 0$ be constants and put:
 $$\gamma=\max(\alpha, \beta).$$
Now we verify the general version which can be applied to 
estimate the orbits of  the state dynamics.
The argument is more complicated.
\begin{prop}
Suppose these satisfy the following estimates:
\begin{align*}
& |p_i^{j+1} - x_i^{j+1}|  \leq c \max( |q_i^j-y_i^j| , |p_i^j-x_i^j|, \dots,| p_{i+\alpha}^j-x_{i+\alpha}^j|)  + T \\
&  |q_{i+1}^j - y_{i+1}^j|  \leq c \max( |q_i^j-y_i^j| , |p_i^j-x_i^j|, \dots,| p_{i+\beta}^j-x_{i+\beta}^j|)  + T \\
\end{align*}
Then they satisfy the estimates:
\begin{align*}
  |p_i^{j+1} - x_i^{j+1}| , & \   |q_{i+1}^j - y_{i+1}^j| \\
& \leq P_{i+j(\gamma+1)}(c) T +  \tilde{c}^{i+1+j(\gamma+1)}|\{p_i^0,x_i^0\}_i ; \{q^j_0, y^j_0\}_j|.
\end{align*}
\end{prop}
{\em Proof:} We split the proof into several steps.
\\
{\bf Step 1:} Firstly we claim that the estimates below hold by induction on $i$:
$$|q_i^0-y_i^0|\leq P_{i-1}(c) T +\tilde{c}^i |\{p_i^0,x_i^0\}_i ; \{q^j_0, y^j_0\}_j|.$$ 
For $i=0$, the estimates $|q_0^0-y_0^0| \leq   |\{p_i^0,x_i^0\}_i ; \{q^j_0, y^j_0\}_j|$ hold
by definition.

Suppose the claim hold up to $i \geq 0$. Then:
\begin{align*}
 |q_{i+1}^0-y_{i+1}^0| & \leq c\max( |q_i^0-y_i^0|, 
  |p_i^0-x_i^0|, \dots,| p_{i+\beta}^0-y_{i+\beta}^0|)  + T \\
& \leq c( P_{i-1}(c) T +\tilde{c}^i |\{p_i^0,x_i^0\}_i ; \{q^j_0, y^j_0\}_j|)+T \\
& \leq P_i(c)T +\tilde{c}^{i+1}  |\{p_i^0,x_i^0\}_i ; \{q^j_0, y^j_0\}_j|.
\end{align*}
Thus they hold up to $i+1$. This verifies the claim.

Then we have the estimates:
\begin{align*}
|p_i^1-x_i^1| & \leq c \max( |q_i^0-y_i^0| , 
 |p_i^0-x_i^0|, \dots,| p_{i+\alpha}^0-y_{i+\alpha}^0|)  + T \\
&  \leq c( P_{i-1}(c) T +\tilde{c}^i |\{p_i^0,x_i^0\}_i ; \{q^j_0, y^j_0\}_j|)+T \\
& \leq P_i(c) T +\tilde{c}^{i+1} |\{p_i^0,x_i^0\}_i ; \{q^j_0, y^j_0\}_j|.
\end{align*}
\\
{\bf Step 2:}
Next we claim that   the estimates hold:
$$|q_i^1-y_i^1| \leq P_{i+\beta}(c)T+ \tilde{c}^{i+\beta+1} |\{p_i^0,x_i^0\}_i ; \{q^j_0, y^j_0\}_j|.$$
We proceed by induction on $i$.
For $i=0$, the estimates:
 $$|q_0^1-y_0^1| \leq  |\{p_i^0,x_i^0\}_i ; \{q^j_0, y^j_0\}_j|
\leq P_{\beta}(c)T+ \tilde{c}^{\beta+1} |\{p_i^0,x_i^0\}_i ; \{q^j_0, y^j_0\}_j|$$
 hold by definition, since $\tilde{c} \geq 1$ holds.
 
 Suppose the above estimates hold 
up to $i \geq0$. Then:
\begin{align*}
|q_{i+1}^1-y_{i+1}^1| 
&  \leq c\max(|q_i^1-y_i^1|, |p_i^1- x_i^1|,  
 \dots,  |p_{i+\beta}^1-x_{i+\beta}^1|)+T \\
& \leq c \max(|q_i^1-y_i^1|,  P_{i+\beta}(c)T  +\tilde{c}^{i+\beta+1}  |\{p_i^0,x_i^0\}_i ; \{q^j_0, y^j_0\}_j|  ) + T \\
& \leq c (  P_{i+\beta}(c)T  +\tilde{c}^{i+\beta+1}  |\{p_i^0,x_i^0\}_i ; \{q^j_0, y^j_0\}_j|  ) + T \\
 & =P_{i+\beta+1}(c) T +  \tilde{c}^{i+\beta+2}  |\{p_i^0,x_i^0\}_i ; \{q^j_0, y^j_0\}_j|   
\end{align*}
where we used step $1$ at the second  inequalities.
So the above estimates also hold for $i+1$, and we have verified the claim.

Then we have   the estimates:
\begin{align*}
|p_i^2-x_i^2|
 & \leq c \max( |q_i^1-y_i^1| , |p_i^1-x_i^1|,..,| p_{i+\alpha}^1-x_{i+\alpha}^1|) + T \\
& \leq c\max( P_{i+\beta}(c)T+\tilde{c}^{i+\beta+1} |\{p_i^0,x_i^0\}_i ; \{q^j_0, y^j_0\}_j| , \\
& \qquad \qquad  \qquad P_{i+\alpha}(c)T+ 
\tilde{c}^{i+\alpha+1}  |\{p_i^0,x_i^0\}_i ; \{q^j_0, y^j_0\}_j| ) +T \\
& = P_{i+\gamma+1}(c)T +\tilde{c}^{i+\gamma +2}  |\{p_i^0,x_i^0\}_i ; \{q^j_0, y^j_0\}_j| .
\end{align*}
\\
{\bf Step 3:}
Let us  verify the estimates for the general case by induction on $j$.
 In step $1, 2$,
we have verified the conclusions for $j\leq 1$ and all $i \geq 0$. 

Suppose the conclusions hold up to $j-1 \geq 1$ and all $i \geq 0$.
Firstly 
let us  consider the pair $(y_i^j,q_i^j)$. We claim that the estimates hold for all $i \geq 0$:
 $$ |q_i^j- y_i^j|
\leq P_{i-1+j(\gamma+1)}(c) T + 
 \tilde{c}^{i+j(\gamma+1)}|\{p_i^0,x_i^0\}_i ; \{q^j_0, y^j_0\}_j|.$$
The estimates  $|q_0^j-y_0^j| \leq  |\{p_i^0,x_i^0\}_i ; \{q^j_0, y^j_0\}_j| $ hold
by definition.
Let us proceed by induction on $i$. Suppose the above estimates hold 
up to $i \geq 0$. Then:
\begin{align*}
& |q_{i+1}^j  -y_{i+1}^j|
 \leq c \max(|q_i^j-y_i^j|, |p_i^j- x_i^j|,  
 \dots, |p_{i+\beta}^j-x_{i+\beta}^j|)+T \\
& \leq c \max( P_{i-1+j(\gamma+1)}(c)T 
+ \tilde{c}^{i+(j+1)(\gamma+1)} |\{p_i^0,x_i^0\}_i ; \{q^j_0, y^j_0\}_j| , \\
& \qquad \qquad P_{i+\beta+( j-1)(\gamma+1)}(c)T+ \tilde{c}^{i+\beta+1+(j-1)(\gamma+1)} 
 |\{p_i^0,x_i^0\}_i ; \{q^j_0, y^j_0\}_j| ) + T \\
 &   =P_{i+j(\gamma+1)}(c)T 
+ \tilde{c}^{i+1+j(\gamma+1)}  |\{p_i^0,x_i^0\}_i ; \{q^j_0, y^j_0\}_j| .
\end{align*}
So the  estimates also hold for $i+1$, and we have verified the claim.

Then 
we claim that  the estimates:
$$|p_i^{j+1}-x_i^{j+1}|\leq P_{i+j(\gamma+1)}(c)T 
+\tilde{c}^{i+1+j(\gamma+1)} |\{p_i^0,x_i^0\}_i ; \{q^j_0, y^j_0\}_j| $$
hold. 
This follows from the following:
\begin{align*}
& |p_i^{j+1} -  x_i^{j+1}|
  \leq c \max( |q_i^j-y_i^j| , |p_i^j-x_i^j|,  \dots,| p_{i+\alpha}^j-x_{i+\alpha}^j|) +T \\
& \leq c\max( P_{i-1+j(\gamma+1)}(c)T+ \tilde{c}^{i+j(\gamma+1)} |\{p_i^0,x_i^0\}_i ; \{q^j_0, y^j_0\}_j| , \\
& \qquad \qquad P_{i+\alpha+(j-1)(\gamma+1)}(c)T   
 +\tilde{c}^{i+\alpha+1+(j-1)(\gamma+1)} |\{p_i^0,x_i^0\}_i ; \{q^j_0, y^j_0\}_j| )    +T \\
& = P_{i+j(\gamma+1)}(c)T +\tilde{c}^{i+1+j(\gamma+1)} |\{p_i^0,x_i^0\}_i ; \{q^j_0, y^j_0\}_j| .
\end{align*}
Thus we have verified the claim under the induction hypothesis up to  $j-1$.

This completes the proof.
\vspace{3mm} \\
{\bf 2.C.2 Uniform estimates between the orbits  for Mealy automata:}
Let $\psi$ and $\phi$ be a pair of relative $(\max, +)$-functions of two variables.
Let  $c$ and $M$ be  the maximums of the Lipschitz constants 
and the numbers of their components respectively.

Let $\psi_t$ and $\phi_t$  be the tropical 
 correspondences to $\psi$ and $\phi$ respectively, and
  consider  the systems of the equations:
\begin{align*}
& x(i,j+1)= \psi(y(i,j), x(i,j)), \\
& y(i+1,j)= \phi(y(i,j), x(i,j)), \\
& \\
 & x'(i,j+1)= \psi_t(y'(i,j), x'(i,j)), \\
& y'(i+1,j)= \phi_t(y'(i,j), x'(i,j)).
 \end{align*}
 with the same initial values 
 $x'(i,0)=x(i,0)=x_i$ and $y'(0,j)=y(0,j)=y^j$.
 
  \begin{lem}
 The uniform estimates hold:
$$|x(i,j+1)-x'(i,j+1)|, \  |y(i+1,j)-y'(i+1,j)| \
 \leq P_{i+j}(c) \log_t M.$$
\end{lem}
{\em Proof:}
By lemma $2.1$, both the estimates hold:
\begin{align*}
& |x(i,j+1) - x'(i,j+1)| 
 =|\psi(y(i,j),  x(i,j))-  \psi_t(y'(i,j) ,  x'(i,j))| \\
&\leq |\psi(y(i,j),  x(i,j))-  \psi(y'(i,j) , x'(i,j))| \\
& \qquad \qquad +
|\psi(y'(i,j),  x'(i,j))-  \psi_t(y'(i,j) , x'(i,j))| \\
& \leq c \max( | y(i,j)-y'(i,j)|,  |x(i,j)-x'(i,j)|) + \log_t M ,\\
& \\
& |y(i+1,j) - y'(i+1,j)| 
 =|\phi(y(i,j),  x(i,j))-  \phi_t(y'(i,j) ,  x'(i,j))| \\
&\leq |\phi(y(i,j),  x(i,j))-  \phi(y'(i,j) , x'(i,j))| \\
& \qquad \qquad +
|\phi(y'(i,j),  x'(i,j))-  \phi_t(y'(i,j) , x'(i,j))| \\
& \leq c \max( | y(i,j)-y'(i,j)|,  |x(i,j)-x'(i,j)|) + \log_t M .\\
\end{align*}
By applying lemma $2.10$ for
$p_i^j=x'(i,j)$, $x(i,j)=x_i^j$, $q_i^j=y'(i,j)$ and $y_i^j=y(i,j)$ with $T=\log_t M$,
one obtains the desired result.

This completes the proof.
\vspace{3mm} \\
{\bf 2.C.3 Uniform estimates for the state dynamics:}
Now let us consider the general case.
Let $\psi$ and $\phi$ be a pair of relative $(\max, +)$-functions as in $2.B$.
 
 Let us choose   the initial data  $\{x_i\}_{j \geq 0}$ and $\{y^j\}_{j \geq 0}$
respectively, and consider two state dynamics  given by the systems of the equations:
\begin{align*}
& x(i,j+1)= \psi(y(i,j), x(i,j), \dots, x(i+\alpha,j)) \\
& y(i+1,j)= \phi(y(i,j), x(i,j), \dots,x(i+\beta,j)),
& \\
& \\
& x'(i,j+1)= \psi_t(y'(i,j), x'(i,j),\dots,  x'(i+\alpha,j)) \\
& y'(i+1,j)= \phi_t(y'(i,j), x'(i,j),\dots  x'(i+\beta,j)),
\end{align*}
with the same  initial values:
 $$x(i,0)= x'(i,0)=x_i, \qquad y(0,j)= y'(0,j)=y^j.$$

\vspace{3mm}

We verify the following:
 \begin{prop}
 The uniform estimates hold:
$$ |x(i,j+1)-x'(i,j+1)|, \ 
  |y(i+1,j)-y'(i+1,j)| \leq P_{i+j(\gamma+1)}(c)\log_t M.$$
 \end{prop}
{\em Proof:}
The proof is parallel to lemma $2.12$.
By lemma $2.1$, one has the estimates:
\begin{align*}
& |x(i,j+1)-x'(i,j+1)|\\
& = |\psi(y(i,j), x(i,j),.., x(i+\alpha,j)) 
- \psi_t(y'(i,j),x'(i,j), .., x'(i+\alpha,j))| \\
& \leq  |\psi(y(i,j), x(i,j),.., x(i+\alpha,j)) 
- \psi(y'(i,j),x'(i,j), .., x'(i+\alpha,j))| \\
& \qquad  +  |\psi(y'(i,j), x'(i,j),.., x'(i+\alpha,j)) 
- \psi_t(y'(i,j),x'(i,j), .., x'(i+\alpha,j))| \\
& \leq c \max( |y(i,j)-y'(i,j)| , |x(i,j)-x'(i,j)|,  \\
&  \qquad \qquad \qquad \qquad \qquad \qquad
..,| x(i+\alpha,j)-x'(i+\alpha,j)|) 
+ \log_t M, \\
& \\
& |y(i+1,j)-y'(i+1,j)| \\
& =
|\phi(y(i,j), x(i,j),.., x(i+\beta,j)) -  \phi_t(y'(i,j),x'(i,j), .., x'(i+\beta,j))| \\
& \leq  |\phi(y(i,j), x(i,j),.., x(i+\beta,j))
 - \phi(y'(i,j),x'(i,j), .., x'(i+\beta,j))| \\
& \qquad +  |\phi(y'(i,j), x'(i,j),.., x'(i+\beta,j)) 
- \phi_t(y'(i,j),x'(i,j), .., x'(i+\beta,j))| \\
& \leq c \max(|y(i,j)-y'(i,j)|, x(i,j)- ,x'(i,j)|,  \\
&  \qquad \qquad \qquad \qquad \qquad \qquad
 .., |x(i+\beta,j)-x'(i+\beta,j)|)+\log_t M.
 \end{align*}
By applying proposition $2.11$ for 
$p_i^j=x'(i,j)$, $x(i,j)=x_i^j$, $q_i^j=y'(i,j)$ and $y_i^j=y(i,j)$ with $T=\log_t M$,
one obtains the result.

This completes the proof.
\vspace{3mm}\\
{\bf 2.C.4 Initial value dependence:}
Let $\phi$, $\psi$ and $M,c$ be as in $2.B$, and
 consider the state  dynamics:
\begin{align*}
& x(i,j+1)= \psi(y(i,j), x(i,j), \dots, x(i+\alpha,j)), \\
& y(i+1,j)= \phi(y(i,j), x(i,j), \dots,x(i+\beta,j)).
\end{align*}
Here we consider how the initial values influence on the long time
behavioir of the dynamics.
For $l=1,2$, let  $\{x_i(l)\}_{j \geq 0}$ and $\{y^j(l)\}_{j \geq 0}$ be two initial data, 
and denote the  corresponding solutions by $\{(x_l(i,j),y_l(i,j))\}_{l=1,2}$
with $x_l(i,0)=x_i(l)$ and $y_l(0,j)=y^j(l)$.
Recall $|\{x_i(1),x_i(2)\}_i ; \{y^j(1), y^j(2)\}_j| $ in $2.C$.
\begin{lem} The estimates hold:
\begin{align*}
  |x_1(i, j+1)- x_2(i,j+1)| ,  \  & |y_1(i+1 ,j) - y_2(i+1,j)| \\
& \leq   \tilde{c}^{i+1+j(\gamma+1)}|\{x_i(1),x_i(2)\}_i ; \{y^j(1), y^j(2)\}_j|.
\end{align*}
\end{lem}
{\em Proof:}
Let us consider the estimates:
\begin{align*}
& |x_1(i,j+1)-x_2(i,j+1)|=  \\
& |\psi(y_1(i,j),x_1(i,j), \dots , x_1(i+\alpha,j)) 
- \psi(y_2(i,j),x_2(i,j),  \dots , x_2(i+\alpha,j))| \\
& \leq c \max( |y_1(i,j)-y_2(i,j)| ,\dots ,| x_1(i+\alpha,j)-x_2(i+\alpha,j)|)  \\
& \\
& |y_1(i+1,j)-y_2(i+1,j)|= \\
& |\psi(y_1(i,j), x_1(i,j),.., x_1(i+\beta,j)) -  \psi(y_2(i,j),x_2(i,j), .., x_2(i+\beta,j))| \\
& \leq c \max(|y_1(i,j)-y_2(i,j)|,\dots , |x_1(i+\beta,j)-x_2(i+\beta,j)|).
 \end{align*}
Then by applying proposition $2.11$ for 
$p_i^j=x_1(i,j)$, $x_i^j=x_2(i,j)$
and $q_i^j=y_1(i,j)$, $y_i^j=y_2(i,j)$ with $T=0$, one obtains the
desired estimates. 

This completes the proof.
\vspace{3mm}\\
{\bf 2.D Rational dynamics:}
Let  $\varphi: {\mathbb R}^{\alpha+2} \to {\mathbb R}$
and $ \psi: {\mathbb R}^{\beta+2} \to {\mathbb R}$
be  relative  $(\max,+)$ functions as in $2.B$, with the constants 
 $M$ and $c$ as before.
 
Passing through the scale transform in tropical geometry, 
one obtains two parametrized rational functions
$f_t$ and $g_t$ with respect to $\varphi$ and $\psi$ respectively:
$$ f_t(\bar{z})=   \frac{ t^{\alpha_1} \bar{z}^{\bar{a}_1}+ \dots + t^{\alpha_{\alpha+2}} \bar{z}^{\bar{a}_{\alpha+2}}}
{ t^{\beta_1} \bar{z}^{\bar{b}_1}+\dots +  t^{\beta_{\alpha+2}} \bar{z}^{\bar{b}_{\alpha+2}}}, \ \
 g_t(\bar{w})=   
\frac{ t^{\gamma_1} \bar{w}^{c_1}+ \dots + t^{\gamma_{\beta+2}} \bar{w}^{\bar{c}_{\beta+2}}}
{ t^{\delta_1} \bar{w}^{\bar{d}_1}+ \dots + t^{\delta_{\beta+2}} \bar{w}^{d_{\beta+2}}}.$$

Let us take initial values:
$$0< z_i < \infty, \quad 0< w^j <\infty$$
for $i,j =0,1,2, \dots$

\begin{df} 
The dynamical system given by:
 \begin{align*}
& z(i,j+1)= f_t(w(i,j), z(i,j), \dots, z(i+\alpha,j)), \\
&   w(i+1,j)= g_t(w( i,j),   z(i,j), \dots, z(i+\beta, j)) 
\end{align*}
   is called the  state system of the rational dynamics, 
   where  the initial values  are given by $z(i , 0) = z_i$ and $w(0, j ) = w^j$.
\end{df}

Let $\psi_t$ and $\phi_t$  be the tropical 
 correspondences to $\psi$ and $\phi$ respectively.
 Let us  consider the state dynamics:
\begin{align*}
& x'(i,j+1)= \psi_t(y'(i,j), x'(i,j),\dots,  x(i+\alpha,j)) \\
& y'(i+1,j)= \phi_t(y'(i,j), x'(i,j),\dots  x(i+\beta,j))
\end{align*}
with the   initial values
 $x'(i,0)= \log_t z_i$ and $ y'(0,j) = \log_t w^j$.
 
 \begin{lem}
 The equalities hold:
 $$x'(i,j)= \log_t z(i,j), \quad y'(i,j)=\log_t w(i,j).$$
 \end{lem}
 {\em Proof:}
 This follows from Proposition $2.4$.
 This completes the proof.
 \vspace{3mm}\\
{\bf 2.D.2 Analysis on the equivalent dynamics  for simple case:}
Compare the contents here with  $2.B.2$.

Let  $\varphi^1$ and $\varphi^2$ be relative $(\max, +)$-functions,
which are mutually equivalent with thier  Lipschitz constant $c$.
Let $M$ be the bigger one of the  numbers of the components.

Let $g^1_t$ and $g^2_t$ be the corresponding   functions which are mutually
 tropically equivalent   (see definition $2.1$).

Let us take  $w_0 \in (0, \infty)$ and an initial sequence
 $\{z_0,z_1, \dots\} \in {\mathbb R}_{>0}^{\mathbb N}$ 
by positive numbers. 
 Then one considers the rational dynamics 
defined inductively by the iterations:
$$ w^l_{i+1}= g^l_t( w^l_i, z_i,  \dots, z_{i+ \beta}) \qquad (l=1,2)$$
with $w^l_0=w_0$.
Notice that these orbits take positive values.
\begin{lem}
The uniform estimates hold for all $i \geq 0$:
$$\frac{w_i^1}{w_i^2} , \ \ \frac{w_i^2}{w_i^1} \  \leq \ M^{2P_{i-1}(c)}$$
\end{lem}
{\em Proof:}
Let us put
$x_i=\log_t z_i$ and $y_0=\log_t w_0$, 
 and consider the discrete dynamics 
defined inductively by the iterations:
$$ y^l_{i+1}= \varphi^l( y_i, x_i,  \dots,  x_{i+\beta}), \quad
 (y^l)_{i+1}'= \varphi^l_t(y'_i, x_i,  \dots,  x_{i+\beta})$$
with $y_0^l=(y_0^l)'=y_0$.
The equalities:
 $$y_i^1=y_i^2$$
 hold for all $i \geq 0$,
since $ \varphi^1$ and $ \varphi^2$ are mutually equivalent.

By proposition $2.4$, the equalities hold:
$$(y^l)'_i  = \log_t w_i^l    \qquad (l=1,2)$$
The following estimates follow from lemma $2.9$:
\begin{align*}
\log_t  \{ \max( \frac{w_i^1}{w_i^2} , \ \ \frac{w_i^2}{w_i^1} ) \} & = 
| \log_t w_i^1 - \log_t w_i^2| 
 = |(y^1)_i' - (y^2)'_i|  \\
& \leq  |(y^1)_i' - y^1_i| + |y^1_i - y^2_i| + |y^2_i- (y^2)'_i|  \\
& \leq 2 P_{i-1}(c)  \log_t M = \log_t M^{2P_{i-1}(c)}.
\end{align*}
Since $\log_t$ are increasing, these estimates imply the desired one.

This completes the proof.

\vspace{3mm}

Later on we will use the notations:
$$( \frac{w}{w'})^{\pm1 } \equiv 
 \max( \ \frac{w}{w'} , \ \ \frac{w'}{w} \ ).$$
{\bf 2.D.3 Analysis on the equivalent dynamics:} 
For $l=1,2$,  let
 $\{w^j(l)\}_{j \geq 0}$ and $ \{z_i(l)\}_{i \geq 0}$ be sequences
 by positive numbers.

\begin{df}
The initial rate is given by the positive number:
$$[ \{(z_i(l), w^j(l)) \}_{l=1}^2] 
= \sup_{i,j} \max( \frac{z_i(1)}{z_i(2)} , \frac{z_i(2)}{z_i(1)} , \frac{w^j(1)}{w^j(2)},  \frac{w^j(2)}{w^j(1)} )
 \geq 1.$$
\end{df}

Let us put $\log_t z_i(l)= x_i(l)$ and $\log_t w^j(l)= y^j(l)$.
Then the equality holds:
$$\log_t [ \{(z_i(l), w^j(l)) \}_{l=1}^2] = |\{x_i(1),x_i(2)\}_i ; \{y^j(1), y^j(2)\}_j|.$$
For the later notation, see $2.C$.

\vspace{3mm}

Let $(\psi^1, \phi^1)$ and $(\psi^2, \phi^2)$ be paris of $(\max,+)$-functions
so that $\psi^1 \sim \psi^2$ and  $\phi^1 \sim \phi^2$ are mutually equivalent.
We say that $(\psi^1, \phi^1)$ and $(\psi^2, \phi^2)$ are {\em pairwisely equivalent}.

Let  $c_{\varphi}$ and $M_{\varphi}$ be
 the Lipschitz constants and the numbers of the components  for $\varphi$ respectively,
 and  put:
$$ c= \max(c_{\psi^1},c_{\psi^2}, c_{\phi^1}, c_{\phi^2}), \ \
M= \max(M_{\psi^1},M_{\psi^2}, M_{\phi^1}, M_{\phi^2}).$$
Let $(f_t^1,g_t^1)$ and $(f_t^2,g_t^2)$ be the  relatively elementary functions
with respect to $(\psi^1, \phi^1)$ and $(\psi^2, \phi^2)$.

Let $\{w^j(l)\}_{j \geq 0}$ and $\{z_i(l)\}_{i \geq 0}$ be the initial sequences by positive numbers, 
and denote the  solutions by $(z_i^j(l),w_i^j(l))$
 to the state systems of the rational dynamics:
 \begin{align*}
& z_i^{j+1}(l)= f^l_t(w_i^j(l), z_i^j(l), \dots, z_{i+\alpha}^j(l)), \\
& w_{i+1}^j(l)= g^l_t(w_i^j(l),   z_i^j(l),\dots, z_{i+\beta}^j(l) ).
\end{align*}
with the  initial values 
$z_i^0(l)=z_i(l)$ and $w_0^j(l)=w^j(l)$  respectively.

\begin{prop}
The uniform estimates hold:
\begin{align*}
\max(  & \frac{z_i^{j+1}(1)}{z_i^{j+1}(2)} , \   \frac{z_i^{j+1}(2)}{z_i^{j+1}(1)} , \
 \frac{w_{i+1}^j(1)}{w_{i+1}^j(2)}, \ \frac{w_{i+1}^j(2)}{w_{i+1}^j(1)}) \\
 & \\
& \leq      M^{2P_{i+j(\gamma+1)}(c)} [ \{(z_i(l), w^j(l)) \}_{l=1}^2] ^{\tilde{c}^{i+1+j(\gamma+1)}} .
\end{align*}
\end{prop}
In particular if the initial values are the same: 
$$z_i(1)=z_i(2),  \ \ w^j(1)=w^j(2)$$
 then they satisfy the uniform estimates:
$$
\max(   \frac{z_i^{j+1}(1)}{z_i^{j+1}(2)} , \   \frac{z_i^{j+1}(2)}{z_i^{j+1}(1)} , \
 \frac{w_{i+1}^j(1)}{w_{i+1}^j(2)}, \ \frac{w_{i+1}^j(2)}{w_{i+1}^j(1)}) \\
\leq      M^{2P_{i+j(\gamma+1)}(c)} .$$

Notice that for the Mealy case ($\alpha=\beta=0$), ther rates are bounded by: 
$$M^{2P_{i+j}(c)}   [ \{(z_i(l), w^j(l)) \}_{l=1}^2] ^{\tilde{c}^{i+j+1}}.$$

{\em Proof of proposition $2.17$:}
Let us consider the solutions to the equations:
\begin{align*}
& x_l'(i,j+1)= \psi_t^l(y_l'(i,j), x_l'(i,j), \dots, x_l'(i+\alpha,j)) \\
& y_l'(i+1,j)= \phi_t^l(y_l'(i,j), x'_l(i,j), \dots, x_l'(i+\beta,j))
\end{align*}
with the initial values $x_l'(i,0)= \log_t z_i(l)$ and $y_l'(0,j)=\log_t w^j(l)$
for $l=1,2$.

With the same initial values, let us also consider  the solutions to another  equations:
\begin{align*}
& x_l(i,j+1)= \psi^l(y_l(i,j), x_l(i,j), \dots, x_l(i+\alpha,j)) \\
& y_l(i+1,j)= \phi^l(y_l(i,j), x_l(i,j), \dots, x_l(i+\beta,j))
\end{align*}
with  $x_l(i,0)= x_l'(i,0)$ and $y_l(0,j)= y_l'(0,j)$.

By proposition $2.13$,  the estimates hold for $l=1,2$:
$$ |x_l(i,j+1)-x_l'(i,j+1)|, \ 
 |y_l(i+1,j)-y_l'(i+1,j)| \leq P_{i+j(\gamma+1)}(c) \log_t M.$$
 On the other hand by lemma $2.14$, the estimates hold:
\begin{align*}
  |x_1(i, j+1)- x_2(i,j+1)| ,  \  & |y_1(i+1 ,j) - y_2(i+1,j)| \\
& \leq \tilde{c}^{i+1+j(\gamma+1)}|\{x_i(1),x_i(2)\}_i ; \{y^j(1), y^j(2)\}_j|.
\end{align*}

Thus combining with these,  the estimates hold:
\begin{align*}
 |x_1' & (i,j)-x_2'(i,j)|  \\
&  \leq
  |x_1'(i,j)-x_1(i,j)| + |x_1(i,j)-x_2(i,j)| +|x_2(i,j)-x_2'(i,j)| \\
&  \leq 2P_{i+(j-1)(\gamma+1)}(c) \log_t M + 
 \tilde{c}^{i+1+(j-1)(\gamma+1)}|\{x_i(1),x_i(2)\}_i ; \{y^j(1), y^j(2)\}_j|, \\
 |y_1'  & (i,j)-y_2'(i,j)| \\
&  \leq   |y_1'(i,j)-y_1(i,j)| + |y_1(i,j)-y_2(i,j)|+  |y_2(i,j)-y_2'(i,j)| \\
& \leq 2 P_{i-1+j(\gamma+1)}(c)\log_t M + \tilde{c}^{i+j(\gamma+1)}|\{x_i(1),x_i(2)\}_i ; \{y^j(1), y^j(2)\}_j| .
\end{align*}
Since
$\log_t F^l(i,j)= x_l'(i,j)$ and $ \log_t G^l(i,j)=y_l'(i,j)$ hold
by proposition $2.4$, 
these estimates
verify the conclusions.
This completes the proof.
\vspace{3mm} \\
{\bf 2.E Dynamical inequalities:}
Let  $\varphi^1$ and $\varphi^2$ be
relatively $(\max,+)$-functions with 
 $M = \max(M_{\varphi^1},  M_{\varphi^2}  )$.

Suppose they are mutually equivalent and 
so  the equality $c=c_{\varphi^1} = c_{\varphi^2}  $ holds.
Let us denote tropically equivalent functions by
 $g^1_t$ and $g^2_t$ correspondingly.

Notice that there are some cases where the inequalities 
$g^1_t \leq g^2_t$ hold. For example if $\varphi^2$ has the presentation
as $\varphi^2 = \max(\varphi^1, \varphi^1)$, 
then $g^1_t \leq g^2_t = 2g^1_t$ holds.

Let us analyze orbits which admit  dynamical inequalities.
Here we start from the simple case  as in $2.D.2$.
Let us take $w_0 \in (0, \infty)$ and 
 an initial sequence $\{z_0,z_1, \dots\} \subset (0,  \infty) $. Then
consider the solutions $\{w^l_i\}_i$ to the equations with $w^l_0=w_0$:
$$w^l_{i+1}= g^l_t( w^l_i,z_i, \dots, z_{i+\beta}).$$

\begin{lem} With the same $w_0$,
suppose  that another sequence $\{w_i\}_i$ satisfies the dynamical inequality:
$$ g^1_t(w_i,z_i, \dots, z_{i+\beta}) \leq  w_{i+1}  \leq  g^2_t( w_i,z_i, \dots, z_{i+\beta}) $$
for all $i \geq 0$. Then the uniform estimates hold:
$$(\frac{w_i}{w_i^1})^{\pm 1}, \ (\frac{w_i}{w_i^2})^{\pm 1} \  \leq \ M^{2P_{i-1}(c)}.$$
\end{lem}
{\em Proof:}  
Let us put $x_i= \log_t z_i$, $y_i=\log_t w_i$,
  $y^l_i = \log_t w^l_i$ and:
$$\bar{x}_i =(x_i, \dots, x_{i+\beta}).$$
Then the equations hold by proposition $2.4$:
 $$y^l_{i+1}=  \varphi_t^l ( y^l_i, \bar{x}_i) \equiv 
\varphi_t^l ( y^l_i, x_i, \dots, x_{i+\beta}).$$ 
Since $\log_t$ are increasing, the estimates hold:
$$\varphi_t^1 (y_i,\bar{x}_i) \ \leq  \ y_{i+1} 
 \ \leq \  \varphi_t^2 ( y_i,\bar{x}_i).$$
Notice that 
 the equivalent functions take the same values:
$$\varphi^1(y,\bar{x}) = \varphi^2( y,\bar{x}) .$$

Now we claim that the estimates:
$$|y_i - y_i^1| \leq 2P_{i-1}(c) \log_t M$$
hold for all $i \geq 0$. For $i=0$, 
$y_0=y_0^1$ holds.

Let us take any $i\geq 0$, and  divide into two cases.

Firstly 
suppose $y_{i+1} \geq y_{i+1}^1$ hold. Then the estimates hold by  lemma $2.1$:
\begin{align*}
0 & \leq  y_{i+1}-y_{i+1}^1 \leq \varphi_t^2 (y_i, \bar{x}_i) - \varphi_t^1( y^1_i,\bar{x}_i) 
 =  |\varphi_t^2 (y_i,\bar{x}_i) - \varphi_t^1(y^1_i, \bar{x}_i)| \\
& \leq  |\varphi_t^2 ( y_i,\bar{x}_i) - \varphi^2(y_i,\bar{x}_i) |\\
& \qquad \qquad  +
|\varphi^1(y_i,\bar{x}_i) - \varphi^1( y^1_i,\bar{x}_i) | + 
|\varphi^1 (y^1_i,\bar{x}_i) - \varphi_t^1(y^1_i,\bar{x}_i)| \\
& \leq 2 \log_t M + c|y_i- y_i^1|.
\end{align*}

Conversely suppose $y_{i+1} \leq y_{i+1}^1$ hold. Then the estimates hold:
\begin{align*}
0 & \leq  y^1_{i+1}-y_{i+1} \leq \varphi_t^1 (y^1_i, x_i) - \varphi_t^1(y_i,x_i) 
 =  |\varphi_t^1 (y^1_i,x_i) - \varphi_t^1( y_i,x_i)| \\
& \leq  |\varphi_t^1 ( y^1_i,x_i) - \varphi^1(y^1_i,x_i) |\\
& \qquad \qquad +
|\varphi^1(y^1_i,x_i) - \varphi^1(y_i,x_i) | +
|\varphi^1 (y_i,x_i) - \varphi_t^1( y_i,x_i)| \\
& \leq 2 \log_t M + c|y_i- y_i^1|.
\end{align*}

Thus the estimates 
$|y_{i+1}-y_{i+1}^1| \leq 2 \log_t M + c|y_i-y_i^1|$ hold  in any case.
By iteration, 
\begin{align*}
& |y_{i+1}- y_{i+1}^1|  \leq 2 \log_t M + c|y_i-y_i^1| \\
& \ \  \leq 2 \log_t M + c\{2  \log_t M + c|y_{i-1}-y_{i-1}^1|\} \\
& =2 (1+c) \log_t M+ c^2 |y_{i-1}-y_{i-1}^1| \\
& \leq \dots \leq 2 P_i(c) \log_t M+ c^N|y_0-y_0^1| =2 P_i(c)\log_t M
\end{align*}
 hold, since  $y_0=y_0^1 = \log_t w_0$. 
 This verifies the claim.
 
The left hand side is equal to $\log_t  ( \frac{w_{i+1}}{w^1_{i+1}})^{\pm 1}$ and 
the right hand side is equal to $\log_t M^{2P_i(c)}$.
Since $\log_t$ are distance increasing, we obtain the estimates:
$$(\frac{w_{i+1}}{w_{i+1}^1})^{\pm 1} \  \leq \ M^{2P_i(c)}.$$
The estimates  $(\frac{w_{i+1}}{w_{i+1}^2})^{\pm 1} \  \leq \ M^{2P_i(c)}$
are obtained by the same way, and we omit repetition.
This completes the proof.
\vspace{3mm} \\
{\bf 2.E.2 Dynamical inequalities for state dynamics:}
Let $(\psi^1, \phi^1)$ and $(\psi^2, \phi^2)$ be 
paris of $(\max,+)$-functions
so that $\psi^1 \sim \psi^2$ and  $\phi^1 \sim \phi^2$ are pairwisely equivalent.
Let $(f_t^1,g_t^1)$ and $(f_t^2,g_t^2)$ be the pairs of the corresponding relatively elementary functions.

Let us take the initial data $\{w^j\}_{j \geq 0}$ and $\{z_i\}_{i \geq 0}$ by positive numbers, 
and consider the solutions to the state systems of the rational dynamics:
 \begin{align*}
&  F^l(i,j+1)= f^l_t(G^l(i,j), F^l(i,j), \dots, F^l(i+\alpha,j)), \\
& G^l(i+1,j)= g^l_t(G^l( i,j),   F^l(i,j),\dots, F^l(i+\beta,j))
\end{align*}
with $F^l(i,0)=z_i$ and $G^l(0,j)=w^j$ for $l=1,2$.

Now we study the dynamical inequalities:
\begin{thm}
Let us take  another sequences $\{w_i^j\}_{i,j}$ and $\{z_i^j\}_{i,j}$
by positive numbers. Suppose they  satisfy the dynamical 
inequalities:
\begin{align*}
& f^1_t(w_i^j, z_i^j, \dots, z_{i+\alpha}^j)
 \leq z_i^{j+1} \leq  f^2_t(w_i^j, z_i^j, \dots, z_{i+\alpha}^j)  \\
& g^1_t(w_i^j,   z_i^j, \dots,z_{i+\beta}^j) \leq w_{i+1}^j \leq g^2_t(w_i^j,   z_i^j, \dots , z_{i+\beta}^j)
\end{align*}
with the same initial values $z_i^0=z_i$ and $w^j_0=w^j$.

Then the uniform estimates hold for $l=1,2$:
$$( \frac{F^l(i,j+1)}{z_i^{j+1}} )^{ \pm 1} ,   \quad
 (\frac{G^l(i+1,j)}{w_{i+1}^j})^{\pm 1}   \\
 \leq  \ \  M^{2P_{i+j(\gamma+1)}(c)}.$$
\end{thm}
{\em Proof:}
Let us put $p_i^j = \log_t z_i^j$ and $q_i^j= \log_t w_i^j$.
Then 
since $\log_t$ are monotone increasing, the estimates:
\begin{align*}
& \psi^1_t(q_i^j, p_i^j, \dots, p_{i+\alpha}^j) \leq p_i^{j+1}
 \leq  \psi^2_t(q_i^j, p_i^j, \dots,p_{i+\alpha}^j) \\
& \phi^1_t(q_i^j,   p_i^j, \dots, p_{i+\beta}^j) \leq q_{i+1}^j
 \leq \phi^2_t(q_i^j,   p_i^j,  , \dots, p_{i+\beta}^j)
\end{align*}
hold by proposition $2.4$, where $p_i^0=x_i=\log_t z_i$ and $q_0^j=y^j=\log_t w^j$.

Let us put:
  $$x_l'(i,j)= \log_t F^l(i,j),  \quad y_l'(i,j)=\log_t G^l(i,j).$$
 Then 
  $x_l'(i,0)=p_i^0$ and $y_l'(0,j)=q_0^j$ hold, 
  and they satisfy the state  systems of the equations for $l=1,2$:
 \begin{align*}
 &  x_l'(i,j+1)= \psi_t^l(y_l'(i,j), x_l'(i,j), \dots, x_l'(i+\alpha,j)) \\
&  y_l'(i+1,j)= \phi_t^l(y_l'(i,j), x_l'(i,j), \dots, x_l'(i+\beta,j))
\end{align*}

We claim that  the uniform estimates hold for $l=1,2$:
$$|x_l'(i,j+1)-p_i^{j+1}|, \  |y_l'(i+1,j)- q_{i+1}^j| \  \leq \ 2P_{i+j(\gamma+1)}(c) \log_t M.$$

Firstly suppose $p_i^{j+1} \geq x_1(i,j+1) $ hold. 
Then the estimates hold:
\begin{align*}
& 0  \leq  p_i^{j+1} - x_1'(i,j+1) =| p_i^{j+1} - x_1'(i,j+1)|
\\
& \leq |\psi_t^2(q_i^j, p_i^j, \dots, p_{i+\alpha}^j) 
- \psi_t^1(y_1'(i,j), x_1'(i,j), \dots, x_1'(i+\alpha,j))| \\
& \leq 2 \log_t M +
 |\psi^2(q_i^j, p_i^j, \dots, p_{i+\alpha}^j) 
- \psi^1(y_1'(i,j), \dots, x_1'(i+\alpha,j))| \\
& \leq 2 \log_t M  + c \max(|q_i^j-y_1'(i,j)|,
  | p_i^j -  x_1'(i,j)|, \dots,  | p_{i+\alpha}^j -  x_1'(i+\alpha,j)|) .
\end{align*}
Conversely suppose  $p_i^{j+1} \leq x_1'(i,j+1) $ hold. Then the estimates hold:
\begin{align*}
0 & \leq  x_1'(i,j+1)  -p_i^{j+1}=|x_1'(i,j+1)  -p_i^{j+1} |  \\
& \leq |\psi_t^1(q_i^j, p_i^j, \dots, p_{i+\alpha}^j) 
- \psi_t^1(y_1'(i,j), x_1'(i,j), \dots, x_1'(i+\alpha,j))| \\
& \leq 2 \log_t M +
 |\psi^1(q_i^j, p_i^j, \dots, p_{i+\alpha}^j) 
- \psi^1(y_1'(i,j), \dots, x_1'(i+\alpha,j))| \\
& \leq 2 \log_t M 
 + c \max(|q_i^j-y_1'(i,j)|,
  | p_i^j -  x_1'(i,j)|, \dots,  | p_{i+\alpha}^j -  x_1'(i+\alpha,j)|) .
\end{align*}
Thus in any cases, the estimates:
$$|  x_1'(i,j+1)  -p_i^{j+1}| \leq 
 2 \log_t M 
 + c \max(|q_i^j-y_1'(i,j)|,
  \dots,  | p_{i+\alpha}^j-  x_1'(i+\alpha,j)|) $$
hold. By use of the similar argument, one obtains the estimates:
$$|  x_2'(i,j+1)  -p_i^{j+1}| \leq 
 2 \log_t M 
 + c \max(|q_i^j-y_2'(i,j)|,
  \dots,  | p_{i+\alpha}^j-  x_2'(i+\alpha,j)|) .$$
By the same way one obtains the  estimates   for $l=1,2$:
\begin{align*}
&  |  y_l'(i+1,j)   -q_{i+1}^j| \\
 &  \leq  
 2 \log_t M 
 +  c \max(|q_i^j-y_l'(i,j)|, |p_i^j - x_l'(i,j)|, 
  \dots,  | p_{i+\beta}^j -  x_l'(i+\beta,j)|).
  \end{align*}
By applying proposition $2.11$ for
$x_i^j=x_l'(i,j)$ and $y_i^j=y_l'(i,j)$ with $T= \log_t M^2$, 
one obtains the estimates:
$$|x_l'(i,j+1)-p_i^{j+1}|, \  |y_l'(i+1,j)- q_{i+1}^j| \  \leq \ 2P_{i+j(\gamma+1)}(c) \log_t M.$$
Thus we have verified the claim.

From proposition $2.4$,  one has the inequalities:
$$| \log_t \frac{z_i^{j+1}}{F^l(i,j+1)}) | , \ | \log_t \frac{w_{i+1}^j}{G^l(i+1,j)}) |
 \leq \log_t M^{2P_{i+j(\gamma+1)}(c)}.$$
By removing $\log_t$ from the both sides,  the conclusion follows.

This completes the proof.

\vspace{3mm}

Let  us consider the case when the initial data take different values.
Recall $\tilde{c}=\max(1,c)$.

Let us consider the solutions to the state systems:
\begin{align*}
&  F^l(i,j+1)= f^l_t(G^l(i,j), F^l(i,j), \dots, F^l(i+\alpha,j)) \\
& G^l(i+1,j)= g^l_t(G^l( i,j),   F^l(i,j),\dots, F^l(i+\beta,j))
\end{align*}
and denote $F^l(i,0)=z_i(1)$ and $G^l(0,j)=w^j(1)$.
\begin{cor}
Suppose another sequences $\{w_i^j\}_{i,j}$ and $\{z_i^j\}_{i,j}$ satisfy the dynamical 
inequalities:
\begin{align*}
& f^1_t(w_i^j, z_i^j, \dots, z_{i+\alpha}^j) \leq z_i^{j+1}
 \leq  f^2_t(w_i^j, z_i^j, \dots, z_{i+\alpha}^j) , \\
& g^1_t(w_i^j,   z_i^j, \dots,z_{i+\beta}^j) \leq w_{i+1}^j \leq g^2_t(w_i^j,   z_i^j, \dots , z_{i+\beta}^j)
\end{align*}
with $z_i^0 \equiv z_i(2)$ and $w^j_0 \equiv w^j(2)$.

Then the uniform estimates hold  for $l=1,2$:
\begin{align*}
( \frac{F^l(i,j+1)}{z_i^{j+1}})^{ \pm 1},  & \quad
( \frac{G^l(i+1,j)}{w_{i+1}^j})^{\pm 1} \\
&  \leq      M^{4P_{i+j(\gamma+1)}(c)}   [ \{(z_i(k), w^j(k)) \}_{k=1}^2] ^{\tilde{c}^{i+1+j(\gamma+1)}}.
\end{align*}
\end{cor}
{\em Proof:}
Let  us consider another solutions to the systems of  the equations:
\begin{align*}
&  J^l(i,j+1)= f^l_t(K^l(i,j), J^l(i,j), \dots, J^l(i+\alpha,j)), \\
& K^l(i+1,j)= g^l_t(K^l( i,j),   J^l(i,j),\dots, J^l(i+\beta,j))
\end{align*}
with the initial values  $J^l(i,0)=z_i(2)$ and $K^l(0,j)=w^j(2)$.
Then by proposition $2.17$, the estimates hold for $l=1,2$:
$$(  \frac{F^l(i,j+1)}{J^l(i,j+1)} )^{ \pm 1} , \ \
( \frac{G^l(i+1,j)}{K^l(i+1,j)})^{ \pm 1}  \\
 \leq      M^{2P_{i+j(\gamma+1)}(c)} [ \{(z_i(k), w^j(k)) \}_{k=1}^2] ^{\tilde{c}^{i+1+j(\gamma+1)}} .$$
On the other hand by theorem $2.19$, the estimates hold:
$$( \frac{J^l(i,j+1)}{z_i^{j+1}} )^{ \pm 1} ,   \ \
( \frac{K^l(i+1,j)}{w_{i+1}^j})^{\pm 1} 
 \leq  \ \  M^{2P_{i+j(\gamma+1)}(c)}.$$
Then by multiplying these, one obtains the desired estimates:
\begin{align*}
& ( \frac{F^l(i,j+1)}{z_i^{j+1}})^{\pm 1} =
 ( \frac{J^l(i,j+1)}{z_i^{j+1}})^{\pm 1}  ( \frac{F^l(i,j+1)}{J^l(i,j+1)})^{\pm 1} \\
& \leq  M^{2P_{i+j(\gamma+1)}(c)}    
M^{2P_{i+j(\gamma+1)}(c)} [ \{(z_i(k), w^j(k)) \}_{k=1}^2] ^{\tilde{c}^{i+1+j(\gamma+1)}} \\
& =    M^{4P_{i+j(\gamma+1)}(c)} [ \{(z_i(k), w^j(k)) \}_{k=1}^2] ^{\tilde{c}^{i+1+j(\gamma+1)}} .
 \end{align*} 
The desired estimates for $ (\frac{G^l(i+1,j)}{w_{i+1}^j})^{\pm 1}$ are obtained 
by the same way.

This completes the proof.

\vspace{3mm}

Let us take two sets of functions
$\{f_t^k(l)\}_{k,l=1,2} $ and $\{g_t^k(l)\}_{k,l=1,2}$
such that they are all  tropically equivalent respectively:
$$f_t^1(1)\sim  f_t^2(1) \sim f_t^1(2) \sim f_t^2(2) ,
\quad g_t^1(1) \sim g_t^2(1) \sim  g_t^1(2) \sim g_t^2(2).$$
Let $c$ and $M$ be the biggest numbers of their Lipschitz constants 
and the numbers of their components respectively.

Suppose the  sets of sequences $(\{w_i^j(l)\}_{i,j} , \{z_i^j(l)\}_{i,j})$ satisfy the dynamical 
inequalities:
\begin{align*}
& f^1_t(l)(w_i^j(l), z_i^j(l), \dots, z_{i+\alpha}^j(l)) \leq z_i^{j+1}(l)
 \leq  f^2_t(l)(w_i^j(l), z_i^j(l), \dots, z_{i+\alpha}^j(l)) , \\
& g^1_t(l)(w_i^j(l),   z_i^j(l), \dots,z_{i+\beta}^j(l))
 \leq w_{i+1}^j (l)\leq g^2_t(l)(w_i^j(l),   z_i^j(l), \dots , z_{i+\beta}^j(l))
\end{align*}
with $z_i^0(l) \equiv z_i(l)$ and $w^j_0(l) \equiv w^j(l)$.

\begin{cor}
Supppose the above conditions. 
Then the uniform estimates hold:
\begin{align*}
( \frac{z_i^{j+1}(1)}{z_i^{j+1}(2)})^{ \pm 1},  & \quad
( \frac{w_{i+1}^j(1)}{w_{i+1}^j(2)})^{\pm 1} \\
&  \leq      M^{6P_{i+j(\gamma+1)}(c)}   [ \{(z_i(k), w^j(k)) \}_{k=1}^2] ^{\tilde{c}^{i+1+j(\gamma+1)}}.
\end{align*}
\end{cor}
{\em Proof:}
Let $(F_k^1(i,j), G_k^1(i,j)$ be as in theorem $2.19$.
Then the uniform estimates:
\begin{align*}
( \frac{F^1_k(i,j+1)}{z_i^{j+1}(k)})^{ \pm 1},  & \quad
( \frac{G^1_k(i+1,j)}{w_{i+1}^j(k)})^{\pm 1} \\
&  \leq      M^{2P_{i+j(\gamma+1)}(c)}  .
\end{align*}

On the other hand by proposition $2.17$, 
the uniform estimates hold:
\begin{align*}
( \frac{F^1_1(i,j+1)}{F^1_2(i, j+1)})^{ \pm 1},  & \quad
( \frac{G^1_1(i+1,j)}{G^1_2(i+1,j)})^{\pm 1} \\
&  \leq      M^{2P_{i+j(\gamma+1)}(c)}    
[ \{(z_i(k), w^j(k)) \}_{k=1}^2] ^{\tilde{c}^{i+1+j(\gamma+1)}}.
\end{align*}

Thus combing with these estimates, we obtain the following:
\begin{align*}
( \frac{z_i^{j+1}(1)}{z_i^{j+1}(2)})^{ \pm 1}
& = ( \frac{z_i^{j+1}(1)}{F^1_1(i,j+1)})^{ \pm 1}
( \frac{F^1_1(i,j+1)}{F^1_2(i, j+1)})^{ \pm 1}
( \frac{F^1_2(i,j+1)}{z_i^{j+1}(2)})^{ \pm 1},\\
&  \leq      M^{6P_{i+j(\gamma+1)}(c)}  [ \{(z_i(k), w^j(k)) \}_{k=1}^2] ^{\tilde{c}^{i+1+j(\gamma+1)}}.
\end{align*}

We can estimate $( \frac{w_{i+1}^j(1)}{w_{i+1}^j(2)})^{ \pm 1}$ by the same way.
This completes the proof.

\section{Asymptotic comparison for PDE}
{\bf 3.A Rough approximations by discrete dynamics:}
Let us describe our general procedure to approximate 
solutions to the systems of PDE by the state systems of rational dynamics.
 
 Let $f_t  =  \frac{a_t}{b_t}$ and $ g_t  =  \frac{c_t}{d_t}$ be 
  relatively elementary functions of $\alpha+2$
 and $\beta+2$ variables respectively,
 where $a_t,b_t,c_t,d_t$ are all elementary.
Let us
 consider the state systems of the rational dynamics:
\begin{align*}
& z_i^{j+1}= f_t(w_i^j, z_i^j, \dots, z_{i+\alpha}^j), \\
& w_{i+1}^j= g_t(w_i^j,   z_i^j, \dots, z_{i+\beta}^j).
\end{align*}
We will be interpreted the dynamical systems as   the   approximations of the systems 
of the partial differential equations as below.

Let us choose  constants $0< \epsilon \leq 1$. 
Let us 
consider a $C^{\mu+1}$ function
$u: [0, \infty) \times [0, \infty) \to (0, \infty)$, and  take the Taylor expansions: 
 \begin{align*}
  u & (x  + i \epsilon   , s+ j \epsilon)
   = u+i  \epsilon u_x +  j \epsilon u_s+
    \frac{(i\epsilon)^2}{2} u_{2x} +  \frac{(j \epsilon)^2}{2} u_{2s} \\
    & +i j \epsilon^2  u_{xs} +
    \dots  + \frac{(i\epsilon)^{\mu}}{\mu !} 
    u_{\mu x} + \frac{(j \epsilon)^{\mu}}{\mu !} u_{\mu s} \\
&  + \frac{(i\epsilon)^{(\mu +1)}}{(\mu +1)!} u_{(\mu +1)x}(\xi_{ij}) +
  \dots
 + \frac{(j\epsilon)^{(\mu +1)}}{(\mu +1)!} u_{(\mu +1)s}(\xi_{ij}) \\
 \end{align*}
 where 
  $|(x,s)-\xi_{ij}| \leq ( |i| +|j|) \epsilon$ hold.

Let us  introduce the change of variables by:
$$ i = \frac{x}{\epsilon}, \quad j= \frac{s}{\epsilon}, \quad u(x,s) = z(i,j),
\quad v(x,s) =w(i,j)$$
Then we take  the difference, and   insert the Taylor expansions:
\begin{align*}
& z_i^{j+1}- f_t(w_i^j,z_i^j, \dots, z_{i+\alpha}^j) \\
 & =  u(  x,s+\epsilon) 
-   f_t(  v(x,s),   u(x,s), \dots,  u(x+\alpha \epsilon,s)) \\
& = 
\frac{ P_1(\epsilon,t, u,v, u_s, u_x ,\dots, u_{\mu x}) + R_1(\epsilon, t, u,v, \dots, u_{(\mu+1)x}(\xi))}
{ b_t(  v(x,s),    u(x,s), \dots, u(x+\alpha \epsilon,s))}  \\ 
& \equiv  {\bf L}_1(\epsilon, t,u,v,u_s, \dots, u_{\mu x}) 
 + \epsilon^{\mu+1} {\bf E}_1(\epsilon, t,u,v, \dots,
\{ u_{\bar{a}}(\xi_{ij}) \}_{\bar{a},i,j} ) .
\end{align*}
where $P_1$ and $R_1$ are polynomials, and 
each monomial in $R_1$ contains derivatives of $u$ of order $\mu+1$.

Similarly we have the expansions:
\begin{align*}
& w_{i+1}^j- g_t(w_i^j,z_i^j, \dots, z_{i+\beta}^j) \\
 & =  v(  x+\epsilon,s) 
-   g_t( v(x,s),   u(x,s), \dots,  u(x+\beta \epsilon,s)) \\
& =
\frac{ P_2(\epsilon,t, u,v,u_x, v_x, \dots, u_{\mu x}) +  R_2(\epsilon, t, u,v, \dots, u_{(\mu+1)x}(\xi'))}
{ d_t( v(x,s),    u(x,s), \dots, u(x+ \beta \epsilon,s))}  \\ 
& \equiv  {\bf L}_2(\epsilon, t,u,v,u_x, v_x, \dots, u_{\mu x}) \\
& \qquad \qquad + \epsilon^{\mu+1} {\bf E}_2(\epsilon, t,u,v,  \dots,
\{ u_{\bar{a}}(\xi'_{ij}) \}_{\bar{a},i,j} , \{ v_{\bar{a}}(\xi'_{ij}) \}_{\bar{a},i,j} ) 
\end{align*}
where each monomial in $R_1$ contains derivatives of $u$ or  $v$ of order $\mu+1$.

We say that ${\bf L}_i$ and ${\bf E}_i$ are  the {\em leading} and {\em error} terms  respectively.

Once one has chosen a pair of relatively elementary functions $(f_t,g_t)$, 
then the above process determines
a system of PDEs $ P_1=P_2=0$,
while tropical geometry provides  an automaton given by  $(\max,+)$ functions $(\psi, \phi)$.
So the pair $(f_t,g_t)$ plays a role of a bridge to connect between systems of  PDEs and  automata.

Conversely 
a pair of  $(\max,+)$ functions $(\psi, \phi)$ and 
 the order $\mu$  determine the parametrized systems of PDEs:
 \begin{align*}
& P_1(\epsilon, t,  u,v, u_x,u_s,   u_{2x}, \dots, u_{\mu x})=0 \\
& P_2(\epsilon, t,  v,u,u_x,v_x,   u_{2x}, \dots, u_{\mu x})=0 \\
\end{align*}
which  are called 
   the {\em induced  systems of partial differential equations} of order $\mu$,
with respect to the pairs  $(\psi, \phi)$.
\vspace{3mm} \\
{\em Remark 3.1:}
In this paper we  treat  systems of PDEs 
with the same order.
In particular  the Mealy systems are of  first order.
However  a general construction produces
systems of PDEs with various orders, where
say $P_1$ has order $\mu_1$ and $P_2$ order $\mu_2$.
In such cases we will have to choose different scaling parameters
as $ i = \frac{x}{\epsilon^p}$ and $ j= \frac{s}{\epsilon^q}$.
Such cases are treated in [K3] in a general form.
 \vspace{3mm} \\
{\bf 3.A.2 Mealy automaton:}
Let us consider the Mealy dynamics given by:
$$ x(i,j+1)= \psi(y(i,j), x(i,j)), \quad 
 y(i+1,j)= \phi(y(i,j), x(i,j)).$$
  Let $f_t$ and $g_t$ be relatively elementary functions corresponding to $\psi$ and $\phi$ respectively,
  and consider the systems of the rational dynamics:
$$ z_i^{j+1}= f(w_i^j, z_i^j),\quad  w_{i+1}^j= g(w_i^j, z_i^j).$$
\begin{df}
The induced  first order systems of the equations:
\begin{align*}
&  \epsilon \ u_s = f_t(v,u) -u, \\
 &  \epsilon \ v_s  =g_t(v,u) -v
 \end{align*}
are called  the hyperbolic Mealy systems.
\end{df}
Notice that the error term are the followings:
$${\bf E}_1 =  \frac{1}{2} \frac{\partial^2 u}{\partial s^2}, \quad
{\bf E}_2 = \frac{1}{2} \frac{\partial^2 v}{\partial x^2}$$
respectively.
\vspace{3mm} \\
{\bf 3.B Higher distorsions:}
Our main interest here is to study analysis  of  asymptotic growth of  solutions to different PDEs
with respect to higher derivatives  and initial conditions.
In this paper we study globally rough-analytic behaviour of positive solutions to  systems of PDEs,
in terms of  two data introduced below.

 Let $f_t  $ and $ g_t $ be 
  relatively elementary functions of $\alpha+2$
 and $\beta+2$ variables respectively.

\vspace{3mm}

({\bf A}) Let  
$u,u',v,v': [0, \infty) \times [0, \infty) \to (0, \infty)$ be four  functions.
 The {\em initial rate} of the pairs with respect to $\epsilon >0$ is defined by:
\begin{align*}
  [(u,v): & (u',v')]_{\epsilon} \\
 & = \max [   \sup_{(x,s) \in [0, \infty) \times [0, \epsilon] }
\max\{ \frac{u}{u'}(x,s), \frac{u'}{u}(x,s)\}, \\
 & \qquad \qquad \qquad \qquad \qquad  \sup_{(x,s) \in [0,  \epsilon] \times [0, \infty)}
\max\{ \frac{v}{v'}(x,s), \frac{v'}{v}(x,s)\} ].
\end{align*}

\vspace{3mm}

({\bf B})   Let $(u,v): [0, \infty)  \times[ 0, \infty) \to (0, \infty)^2$
be a pair of  functions of   $C^{\mu + 1}$ class.
We  introduce the pointwise  norms  by:
$$  ||(u,v)||_{\mu, \alpha}^1   (x,s) =  
 \begin{cases}
    \max [ 
\sup_{(x,s)  \in [x, x+\alpha \epsilon ]\times \{s\} } 
   |\frac{\partial^{ \mu+1} u}{\partial_x^{\mu+1 }}|(x,s) , & \\
\qquad \qquad \qquad  \sup_{(x,s)  \in\{x\} \times [s, s+ \epsilon ] }
\quad    |\frac{\partial^{ \mu+1} u}{\partial_s^{\mu+1 }}|(x,s) ], &  \alpha  \geq 1 ,  \\
& \\
\sup_{(x,s)  \in\{x\} \times [s, s+ \epsilon ] }
\quad    |\frac{\partial^{ \mu+1} u}{\partial_s^{\mu+1 }}|(x,s) ], &  \alpha  =0 ,
\end{cases}$$
$$
 ||(u,v)||_{\mu,  \beta}^2   (x,s) =  
\begin{cases}
  \max [ 
\sup_{(x,s)  \in [x, x+\beta \epsilon ]\times \{s\} } 
\quad    |\frac{\partial^{ \mu+1} u}{\partial_x^{\mu+1 }}|(x,s) , & \\
\qquad \qquad \qquad \sup_{(x,s)  \in [x, x+ \epsilon ]\times \{s\} } 
    |\frac{\partial^{ \mu+1} v}{\partial_x^{\mu+1 }}|(x,s) ], & \beta \geq 1, \\
& \\
\sup_{(x,s)  \in [x, x+ \epsilon ]\times \{s\} } 
    |\frac{\partial^{ \mu+1} v}{\partial_x^{\mu+1 }}|(x,s) ], & \beta =0 .
\end{cases}$$

   The {\em higher distorsion}  is given by:
$$ 
K(u,v) \equiv  \sup_{(x,s) \in [0, \infty)^2 }   
\max[\frac{||(u,v)||_{ \mu, \alpha}^1}{  u(x,s+ \epsilon)}, \ \
\frac{||(u,v)||_{ \mu, \beta}^2}{  v(x+ \epsilon,s)}].$$

The  {\em error constants} $C$ is the minimum number so that
  the pointwise estimates hold for any functions $u,v$  independently of $\epsilon$ and $t$:
\begin{align*}
& |{\bf E}_1  (\epsilon, t,u,v,  \dots,
\{ u_{\bar{a}}(\xi_{ij}) \}_{\bar{a},i,j})| (x,s) \
    \leq C    ||(u,v)||_{ \mu,  \alpha}^1(x,s)  , \\
& \\
& |{\bf E}_2  (\epsilon, t,u,v, \dots,
\{ u_{\bar{a}}(\xi_{ij}) \}_{\bar{a},i,j} , \{ v_{\bar{a}}(\xi'_{ij}) \}_{\bar{a},i,j})| (x,s) \
    \leq C   ||(u,v)||_{ \mu,  \beta}^2(x,s)  
\end{align*}
\quad
\vspace{3mm} \\
{\em Example 3.1:}
For the Mealy case,  the error constant 
is  $C=  \frac{1}{2}$, and 
we have the followings:
\begin{align*}
&  ||(u,v)||_{\mu, \alpha}^1    =  
\sup_{(x,s)  \in\{x\} \times [s, s+ \epsilon ] }
   |\frac{\partial^2 u}{\partial^2 s}|(x,s) ,  \\
&  ||(u,v)||_{\mu,  \beta}^2    =  
\sup_{(x,s)  \in [x, x+ \epsilon ]\times \{s\} } 
    |\frac{\partial^2 v}{\partial^2 x}|(x,s) .
\end{align*}
 \quad
\vspace{3mm} \\
{\bf 3.C Asymptotic comparisons:}
Let us take four relatively elementary functions $f^1,f^2$ and $g^1,g^2$ so that
$f^1 \sim f^2$ and $g^1 \sim g^2$ are tropically equivalent mutually.
For $l=1,2$, let:
\begin{align*}
& P_1^l(\epsilon, t,  u, v, u_s, \dots, u_{\mu x})=0 , \\
&  P_2^l(\epsilon, t,  u, v, v_x, \dots, u_{\mu x})=0 
\end{align*}
be the induced  systems of PDEs of order $\mu$. 

let $M$ be the largest numbers of their components.

\begin{thm}  For $l=1,2$, let $C$ be the bigger one of their  error constants. 

Let  $(u^l,v^l): [0, \infty) \times [0, \infty) \to (0, \infty)$
be the solutions to  the above systems 
respectively, so that the estimates:
 $$0  \leq   CK(u^l,v^l) \leq (1-  \delta) \epsilon^{-1}$$
are satisfied for some positive $ \delta >0$.

 Then
 they satisfy the asymptotic estimates
 for all $(x,s) \in [0,  \infty) \times [0, \infty)$:
 \begin{align*}
 (\frac{u^1}{u^2})^{\pm 1}(x,s), & \  \ (\frac{v^1}{v^2})^{\pm 1}(x,s) \\
& \leq (N_0 M)^{6P_{\epsilon^{-1}( x+s (\gamma+1))}(c)}  \ \
 ([ (u^1,v^1):(u^2,v^2)]_{\epsilon})^{\tilde{c}^{\epsilon^{-1}(x+s(\gamma+1))+1}}  .
 \end{align*}
where $N_0$  is any integer  with $N_0 \geq  \max(\delta^{-1}, 2-  \delta)$.
 \end{thm}
{\em Proof:}
Let us  choose $0 \leq a,b \leq 1$, and
 put the domain lattices   by:
$$ L_{\epsilon}(a,b)=\{( \ (l_1+a)\epsilon, (l_2+b)\epsilon \ ) \in  [0, \infty) \times [0, \infty)
: l_1,l_2 \in {\mathbb N} \}.$$

Let us put: 
\begin{align*}
& z_i^j(l) = u^l((i+a)\epsilon, (j+b)\epsilon), \\
& w_i^j(l) =v^l((i+a)\epsilon, (j+b)\epsilon)
\end{align*}
and  consider the Taylor expansions:
\begin{align*}
& z_i^{j+1}(l)- f^l_t(w_i^j(l),z_i^j(l), \dots, z_{i+\alpha}^j(l)) \\
& =
  {\bf L}_1(\epsilon, t,u^l,v^l,u^l_s, \dots, u^l_{\mu x}) 
 + \epsilon^{\mu+1} {\bf E}_1(\epsilon, t,u^l,v^l, \dots,
\{ u^l_{\bar{a}}(\xi_{ij}) \}_{\bar{a},i,j} ) , \\
& \\
& w_{i+1}^j(l)- g_t(w_i^j(l),z_i^j(l), \dots, z_{i+\beta}^j(l) ) \\
 & =   {\bf L}_2(\epsilon, t,u^l,v^l,u^l_x, v^l_x, \dots, u^l_{\mu x}) \\
& \qquad \qquad + \epsilon^{\mu+1} {\bf E}_2(\epsilon, t,u^l,v^l,  \dots,
\{ u^l_{\bar{a}}(\xi'_{ij}) \}_{\bar{a},i,j} , \{ v^l_{\bar{a}}(\xi'_{ij}) \}_{\bar{a},i,j} ) 
\end{align*}

By the assumption, both the leading terms satisfy   the equalities:
$${\bf L}_1^l(\epsilon, t,u^l,v^l, u^l_s,  \dots) =0, \quad
{\bf L}_2^l(\epsilon, t,u^l,v^l, u^l_x,  \dots) =0.$$ 
Moreover the error terms satisfy the estimates:
\begin{align*}
&   |{\bf E}_1^l(\epsilon, t,u^l,v^l,  \dots) |(x,s)
\leq C K(u^l,v^l) \ u^l(x,s+\epsilon), \\
&   |{\bf E}_2^l(\epsilon, t,u^l,v^l, \dots) |(x,s)
\leq C K(u^l,v^l) \ v^l(x+\epsilon,s).
\end{align*}

Then combining with these, one obtains the estimates:
\begin{align*}
& |z_i^{j+1}(l)- f^l_t(w_i^j(l),z_i^j(l), \dots, z_{i+\alpha}^j(l)) |
\leq 
 \epsilon^{\mu+1} | {\bf E}_1(\epsilon, t,u^l,v^l, \dots, ) |  \\
& \leq \epsilon CK(u^l,v^l) z_i^{j+1}(l) 
\leq (1- \delta) z_i^{j+1}(l), \\
& \\
& |w_{i+1}^j(l)- g_t(w_i^j(l),z_i^j(l), \dots, z_{i+\beta}^j(l) ) | \leq
  \epsilon^{\mu+1} | {\bf E}_2(\epsilon, t,u^l,v^l,  \dots ) | \\
& \leq \epsilon CK(u^l,v^l) w_{i+1}^j(l)
\leq (1- \delta) w_{i+1}^j(l) .
\end{align*}
In particular there is some integer  $N_0 \geq  \max(\delta^{-1}, 2-  \delta)$
 so that the inequalities hold:
\begin{align*}
&  \frac{1}{N_0}  f^l_t(w_i^j(l),z_i^j(l), \dots, z_{i+\alpha}^j(l))
  \leq  z_i^{j+1}(l) \leq N_0
 f^l_t(w_i^j(l),z_i^j(l), \dots, z_{i+\alpha}^j(l))\\
& \\
&  \frac{1}{N_0}  g_t(w_i^j(l),z_i^j(l), \dots, z_{i+\beta}^j(l) )
 \leq    w_{i+1}^j(l) \leq N_0
 g_t(w_i^j(l),z_i^j(l), \dots, z_{i+\beta}^j(l) ) .
  \end{align*}

Now  $\frac{1}{N_0}f$ and $N_0f$ are both tropically equivalent to $f$, 
with the bounds $\max(M_{\frac{1}{N_0}f}, M_{N_0f}) \leq N_0 M_f$.
Then 
it follows from corollary $2.21$ that the uniform estimates:
\begin{align*}
  (\frac{z_i^{j+1}(1)}{z_i^{j+1}(2)} )^{\pm 1} ,   \ \ &
(  \frac{w_{i+1}^j(1)}{w_{i+1}^j(2)})^{ \pm 1}   
  \leq  (N_0M)^{6P_{i+j(\gamma+1)}(c)}  [ (u^1,v^1):(u^2,v^2)]_{\epsilon}^{\tilde{c}^{i+1+j(\gamma+1)}}   \\
& \leq   (N_0M)^{6P_{\epsilon^{-1}(x+s(\gamma+1))}(c)}
  [ (u^1,v^1):(u^2,v^2)]_{\epsilon}^{\tilde{c}^{\epsilon^{-1}(x+s(\gamma+1) )+1}}  
\end{align*}
where $(x,s)=((i+a)\epsilon, (j+b)\epsilon)$.

Since the right hand side does not depend on $a,b$,
one obains the uniform estimates:
\begin{align*}
(\frac{u^1}{u^2})^{\pm 1}(x,s), & \ (\frac{v^1}{v^2})^{\pm 1}(x,s) \\
& \leq  (N_0M)^{6P_{\epsilon^{-1}( x+s (\gamma+1))}(c)}  
[ (u,v):(u',v')]_{\epsilon}^{c^{\epsilon^{-1}(x+s(\gamma+1))+1}}  
\end{align*}
This completes the proof.
\vspace{3mm}

\section{State dynamics and their scale transforms}
{\bf 4.A Automaton:}
Let us recall the state dynamics. 
Let $S$ and $Q$ be  finite sets 
called the {\em alphabet set} and the {\em state set} respectively.
Let $\alpha , \beta \geq 0$ and consider a pair of functions:
$$\psi: Q \times S^{\alpha+1} \to S, \quad
\phi: Q \times S^{\beta+1} \to Q$$
where we call $\psi$ as the {\em output map} and $\phi$ as the {\em transition map}.
In this paper we call such a pair of functions as an {\em automaton}.

 Let: 
$$X_S =\{ \bar{k} = (k_0,k_1, \dots) : k_i \in S\}$$
be all the set of one-sided strings of infinite length.
Then each state $q \in Q$ induces a continuous map:
$${\bf A}_q: X_S \to X_S$$
given by:
${\bf A}_q(k_0,k_1, \dots) =(k_0',k_1', \dots)$, where $k_i'$ are inductively defined 
as below with $q_0=q$:
$$ k_i' = \psi(q_i,k_i, \dots, k_{i+\alpha}), \quad
 q_{i+1}=\phi(q_i,k_i, \dots, k_{i+\beta}).$$
 We call it the {\em state dynamics} with respect to ${\bf A}$.
 \vspace{3mm} \\
 {\bf 4.A.2 Automata groups:}
 A  {\em Mealy automaton}  ${\bf A}$ over $S$ is given by two functions  
of the form ([GNS]):
$$ \psi: Q \times S \to S, \quad
 \phi: Q \times S \to Q$$
which is  a special case of the state dynamics.

 In this case the state dynamics are given by
${\bf A}_q(k_0,k_1, \dots) =(k_0',k_1', \dots)$, where $k_i'$ are inductively defined 
with $q_0=q$:
$$ k_i' = \psi(q_i,k_i), \quad q_{i+1}=\phi(q_i,k_i).$$
Notice that the Mealy dynamics induce  the level-set actions as
${\bf A}_q: X_A^{N+1} \to X_A^{N+1}$,  where
 $$X_A^{N+1} = \{ \bar{k}^* =(k_0,k_1, \dots, k_N): k_i \in S \}$$
are the set of words of length $N+1$.

Let $m = \sharp |A|$ be the cardinality of $S$, and 
$T_m$ be the rooted regular $m$-tree.
The set of all vertices of $T_m$ can be identified with $X_A^{\infty} \equiv \cup_N X_A^N$.
Thus ${\bf A}_q$ give the actions:
$${\bf A}_q: T_m \to T_m.$$

Let us say that ${\bf A}$ is {\em invertible}, if $\psi(q, \quad) : S \cong S$ are one to one and onto 
for all $q \in Q$. An invertibe automaton ${\bf A}$ gives automorphisms
${\bf A}_q: T_m \cong T_m$.
 \begin{df} Let ${\bf A}$ be invertible.
 The group generated by the set of states:
$$G({\bf A}) = \text{ gen } \{{\bf A}_q : T_m \cong T_m : q \in Q\}$$
is called the automata group.
 \end{df}
\quad
\vspace{3mm} \\
{\em Remark 4.1:}
 (1) Automata groups are subgroups of   Aut$(T_m)$.

(2) General state dynamics give the  maps as ${\bf A}_q: \partial T_m \to \partial T_m$.

\vspace{3mm}

Later on we will always assume that Mealy automata are invertible, and
both $S$ and $Q$ are subsets of the real number:
$$S, \ Q \ \subset \ {\mathbb R}.$$

For $ X_Q^L= \{\bar{q}^L= (q^0, \dots, q^L) : q^i \in Q\}$,
let us put:
$$ {\bf A}_{\bar{q}^L} = {\bf A}_{q^L} \circ \dots {\bf A}_{q^0}: T_m \cong T_m$$

If we denote
 ${\bf A}_{\bar{q}^L}(k_0,k_1, \dots) =(k_0^L, k_1^L, \dots)$, 
 then the state dynamics given by the Mealy automata 
 are exactly given by  ${\bf A}_{\bar{q}^L}$, so that 
  $\{k_i^j\}_{i,j}$ are the orbits 
given by the state dynamics
with the initial data 
$(q^0,q^1, \dots)$ and $(k_0,k_1, \dots)$.
\vspace{3mm} \\
{\bf 4.A.3  Generalization of  Mealy automata:}
One of the most basic properties commonly shared among all Mealy automata
is that the groups are residually finite, since their restrictions induce the level-set actions.
Here we consider a canonical generalization whose groups are not necessarily the case,
which arise from the state dynamics.

Let us consider an automaton ${\bf A}$ given by the two functions:
$$ \psi: Q \times S^{\alpha+1} \to S, \quad
 \phi: Q \times S^{ \beta+1} \to Q$$
 and consider the corresponding state dynamics ${\bf A}_q: X_S \to X_S$
 for each $q \in Q$.
${\bf A}_q: \partial T_m \to \partial T_m$ does not induce   level-set actions in general.

Let us say that ${\bf A}$ is {\em invertible}, if 
${\bf A}_q: X_S \cong X_S$ are isomorphisms
 for all $q \in Q$.
 An invertibe automaton ${\bf A}$ gives automorphisms
${\bf A}_q: \partial T_m \cong \partial T_m$, and the group generated by the set of states is also
denoted by $G({\bf A})$:
$$G({\bf A}) = \text{ gen } \{{\bf A}_q : \partial T_m \cong \partial T_m : q \in Q\}$$
which is a subgroup of the boundary automorphism group  Aut$(\partial T_m)$.
\vspace{3mm}\\
{\em Example 4.1:}
The following automaton is not  Mealy, since
dynamics are not determined by  every level set of the trees.

Let us choose $\alpha=1$. Let $S=\{ s_0,s_1\}$ and $Q=\{q^0,q^1\}$, and 
$\epsilon_0$ be the permutation between two elements.
Let us consider the functions:
\begin{align*}
&  \psi(q^1, \quad, s)=\epsilon_0, \quad 
\psi(q^0, \quad , s) =
\begin{cases}
\epsilon_0 & s=s_0 \\
id & s=s_1
\end{cases}
\\
& \phi(\quad ,s,s')=\epsilon_0
\end{align*}
where $s,s' \in S$.

\begin{lem}
The state dynamics by this automaton is invertible.
\end{lem}
{\em Proof:}
Let us  consider $ \psi(q^0, \quad, \quad)$, and
suppose $(q^0,s,t) \ne (q^0, s',t')$.
The equalities $ \psi(q^0, s, t) = \psi(q^0, s', t')$ hold only  when
$(s,s')=( s_1,s_0) \text{ or } (s_0,s_1)$ and 
$t \ne t'$.

On the other hand at the next level, $\phi(q^0,s,t)=\phi(q^0,s',t')=q^1$ hold, and so: 
$$\psi(q^1, t, \quad) \not= \psi(q^1, t', \quad).$$
This implies ${\bf A}_{q^0}: \partial T_m \to \partial T_m$ is injective.

Since $ \psi(q^1, \quad, \quad)$ does not depend on the second variable,
injectivity of 
 ${\bf A}_{q^1}: \partial T_m \to \partial T_m$ follows from the one  of ${\bf A}_{q^0}$.

Let us consider surjectivity, and take any $(x'_0,x'_1,\dots) \in X_S$.
We seek for some elements $(x_0,x_1 \dots) \in X_S$ with 
${\bf A}_q(x_0,x_1,\dots) =(x_0',x_1', \dots)$.
Notice that the states change periodically as $(q^0,q^1,q^0,q^1, \dots)$
or $(q^1,q^0,q^1, \dots)$.

Both cases can be considered similarly, and we treat  the first case only.
Firstly we choose $x_{2i+1}= x'_{2i+1}+1$ mod $2$ for all $i \geq 0$.
Then we can choose $x_{2i}$ uniquely so that the equalities 
$ \psi(q^0, x_{2i}, x_{2i+1}) =x_{2i}'$ hold for all $i$.
Thus ${\bf A}_q$ are surjective, and so
they are  isomorphisms.
 This completes the proof.
\vspace{3mm}\\
{\bf 4.B Extensions of  automata:}
Let $Q, S \subset {\mathbb R}$ be finite sets, and choose 
an automaton $\psi: Q \times S^{\alpha+1} \to S$ and
$\phi: Q \times S^{\beta+1} \to Q$.

Let: 
$$\tilde{\psi}: {\mathbb R} \times {\mathbb R}^{\alpha+1} \to {\mathbb R}, \quad
\tilde{\phi}:  {\mathbb R} \times {\mathbb R}^{\beta+1} \to {\mathbb R}$$
be two maps.

Let us say that
the pair of functions $(\tilde{\psi}, \tilde{\phi})$  extends the automaton, 
if their restrictions coincide with each other:
$$\tilde{\psi}|Q\times S^{\alpha+1} = \psi, \quad  \tilde{\phi} |Q \times S^{\beta+1} = \phi.$$

Conversely if 
 there are finite subsets $Q,S \subset {\mathbb R}$ so that their restrictions
 of the pair of  functions $(\tilde{\psi}, \tilde{\phi})$
induce the functions:
$$\psi \equiv \tilde{\psi}: Q \times S^{\alpha+1} \to S, \quad
\phi \equiv \tilde{\phi}: Q \times S^{\beta+1} \to Q$$
then we say that  $(\tilde{\psi}, \tilde{\phi})$
restricts to an automaton $(\psi,\phi)$.

Notice  that if a  pair of bounded functions  induces maps as:
$$\tilde{\psi}: {\mathbb Z} \times {\mathbb Z}^{\alpha+1} \to {\mathbb Z}, \quad
\tilde{\phi}:  {\mathbb Z} \times {\mathbb Z}^{\beta+1} \to {\mathbb Z}$$
then they restricts to some automaton 
for some $Q,S \subset {\mathbb Z}$.
\vspace{3mm} \\
{\bf 4.B.2 Lamplighter group:}
Let:
\begin{align*}
& \psi: \{q^0,q^1\} \times \{s_0,s_1\} \to \{s_0,s_1\},
 \quad \psi(q^0, \quad) = id, \  \psi(q^1, \quad) =  \epsilon, \\
& \phi:  \{q^0,q^1\} \times \{s_0,s_1\} \to \{q^0,q^1\}, \quad \phi(\quad , s_i) = q^i \ \ (i=0,1)
\end{align*}
be the  pair which gives a Mealy automaton ${\bf A}$.
It is known ([GZ]) that the associated automata group is isomorphic to 
 the Lamplighter group which acts on the rooted binary tree.

Let us extend the automaton  to  a a pair of $(\max,+)$-functions as below.
Let us put $q^0=s_0=0$ and $q^1=s_1=3 \in {\mathbb Z}$.

Now consider the piecewise-linear functions:
\begin{align*}
& P(k)= - \max(-3, -\max(0, -3(k-2))), \\
& Q(k)= -\max(-3, -\max(0, 3(k-1))).
\end{align*}

Then an extension of the pair $(\psi, \phi)$ is given by the function:
$$ \psi(q,s) =  \max(P(s), Q(s-q)), \quad \phi(q,s)= Q(s)$$
as maps from ${\mathbb Z}^2$ to ${\mathbb Z}$.

The extension $ (\psi,\phi)$
 holds the following properties:

 \vspace{2mm}

(1)  They are constants on the  intervals $[-1,1]$ and $[2,4]$.
Namely the extension satisfies stability defined in $4.C$.

(2)  They are bounded.

 \vspace{3mm}

Now let us put the functions:
$$\alpha_a(z)= (t^{-3}+ \frac{1}{1+ z^3 t^a})^{-1}.$$
The tropical pair of the rational functions $(f_t,g_t)$ is given by:
$$g_t(w,z)= \alpha_{-3}(z), \quad
f_t(w,z) = \alpha_6(z^{-1}) +\alpha_{-3}(z w^{-1})).$$

The induced system of the partial differential equations is given by:
$$  \epsilon u_s = f_t( v,u) -u, \quad
   \epsilon v_x= g_t(u) -v$$
where $f_t$ and $g_t$ are as above.
\vspace{3mm} \\
{\bf 4.B.3 Extensions by relatively $(\max,+)$-functions:}
Let $\varphi: {\mathbb Z}^n \to {\mathbb Z}$ be a map.
Let us note that there are two canonical ways  of the extensions:

(1) By connecting these integer values by segments, 
one can extend $\varphi$ straightforwardly to a piecewise linear map:
$$\varphi: {\mathbb R}^n \to {\mathbb R}.$$

(2) If $\varphi: {\mathbb Z}^n \to {\mathbb Z}$ is equipped with its presentation as a $(\max,+)$-function, 
  then it is canonically extended as  a map $\varphi: {\mathbb R}^n \to {\mathbb R}$.
   \vspace{3mm}  \\
For our purposes in this paper we will choose the method  (2) above.

\begin{lem}
Let us  consider a Mealy automaton
$\psi: Q \times S \to S,  \ \phi: Q \times S \to Q$ with $Q,S \subset {\mathbb R}$.
Then there is a bounded  extension: 
$$\tilde{\psi}: {\mathbb R} \times {\mathbb R} \to {\mathbb R}, \quad
\tilde{\phi}:  {\mathbb R} \times {\mathbb R}  \to {\mathbb R}$$
so that both functions $\tilde{\psi}$ and $ \tilde{\phi}$ 
are represented by relatively $(\max, +)$-functions.
\end{lem}
{\em Proof:}
We construct  local cone extensions inductively.

\vspace{3mm}

{\bf Step 1:}
Let $\{m_1,m_2,m_3\} \subset Q \times S$ be three points which are mutually different,
and consider the triangle $\Delta \subset {\mathbb R}^3$ 
 whose vertices are given   by:
$$O_i =(m_i, \psi(m_i))  \ \in \ {\mathbb R}^3 \qquad (i=1,2,3).$$ 
Let $l_i \subset \Delta$ be segments 
whose end points are $O_i$ and $O_{i+1}$ mod $3$.

Let us choose planes $L_i$ which contain $l_i$ so that they are represented by 
graphs of affine linear
functions   $\varphi_i: {\mathbb R}^2 \to {\mathbb R}$ as
 $$L_i=\{(x_1,x_2, \varphi_i(x_1,x_2)) : (x_1,x_2) \in {\mathbb R}^2\}.$$

Let $ g \in \Delta$ be the barycenter, and $l:[0, \infty)$ be the half lines starting from $g$
which are vertical to ${\mathbb R}^2$ in ${\mathbb R}^3= {\mathbb R}^2 \times {\mathbb R}$.
 Exactly there are two choices $l_b$ and $l_a$
which go below  or above to  $\Delta$.

For $i=1,2,3$, let $L_i^b$ and $L_i^a$ be another planes which contain $l_i$ and 
intersect with $l_b$ and $l_a$ at $l_b(t)$ or $l_a(t)$ respectively.
So three places $\{ L_1^b, L_2^b,L_3^b\}$ intersect at one point $l_b(t)$.
Similar for $a$.
We will choose these planes with sufficiently large $t>>1$.
Let $\varphi_i^b$
and $\varphi_i^a$ be their representations
by relative $(\max, +)$-functions.

Now we have two types of the cones:
$$C_b=\text{ graph } \max( \varphi_1^b, \varphi_2^b, \varphi_3^b)$$
and $C_a$ are similar.
$C_b$ are concave and $C_a$ are convex.
Notice that we may choose arbitrarily sharp slopes of the cones.

Let us choose a large $C$ with
$C \geq \max \{ |\psi(m) | : m  \in Q  \times S \}$.
Let $\varphi$ be another affine linear function whose graph contains $\Delta$.

Now we put the bounded functions:
 \begin{align*}
&  \psi_{\Delta}^b = - \max(-C, -\max(\varphi,  \varphi_1^b, \varphi_2^b, \varphi_3^b)), \\
&  \psi_{\Delta}^a =  \max(-C, - \max(- \varphi,  - \varphi_1^a, - \varphi_2^b, - \varphi_3^b)).
 \end{align*}
Notice that both graphs contain  $\Delta$.

\vspace{3mm}

{\bf Step 2:}
Let us order all the points $\{m_{-1},m_0, m_1,m_2,\dots, m_l \}= Q\times S$,
so that there are families of  two dimensional polytopes   $M_i$
such  that: 

\vspace{2mm}

(1) $M_i$ are given by unions of $M_{i-1}$ with a single triangle $\Delta_i$
which contain the vertices  $m_i$ and

(2) $M_i$ do not contain the sets $\{m_{i+1}, \dots, m_l\}$.
\vspace{2mm}  \\
Now we inductively construct relaive $(\max,+)$-functions
$\tilde{\psi}_i: {\mathbb R}^2 \to {\mathbb R}$ so that  the equalities:
$$\tilde{\psi}_i(m_k)=\psi(m_k)$$
hold  for all $1 \leq k \leq i$.

We have constructed $\psi_1^a=\psi_{\Delta_1}^a $ by Step $1$,
with $\Delta_1= \Delta$.

Suppose we have obtained $\tilde{\psi}_i$, and let us construct $\tilde{\psi}_{i+1}$.
Let us divide into two cases:

Suppose $\tilde{\psi}_i(m_{i+1}) \geq \psi(m_{i+1})$.
Let $\psi_{\Delta_{i+1}}^b$ be as in Step $1$.
By choosing sufficiently sharp slopes, we may assume that the estimates:
$$\psi_{\Delta_{i+1}}^b(m_k) \geq  \tilde{\psi}_i(m_k)$$
hold for all $1 \leq k \leq i$.
Then we put the bounded function by:
$$\tilde{\psi}_{i+1}= - \max( - \psi_{\Delta_{i+1}}^b,  - \tilde{ \psi}_i ).$$

Next suppose $\tilde{\psi}_i(m_{i+1}) \leq \psi(m_{i+1})$.
Let $\psi_{\Delta_{i+1}}^a$ be as in Step $1$.
Again by choosing sufficiently sharp slopes, we may assume that the estimates:
$$\psi_{\Delta_{i+1}}^a(m_k) \leq  \tilde{\psi}_i(m_k)$$
hold for all $1 \leq k \leq i$.
Then we put the bounded function by:
$$\tilde{\psi}_{i+1}=  \max(  \psi_{\Delta_{i+1}}^b,  \tilde{ \psi}_i ).$$

In any cases $\tilde{\psi}_{i+1}$ above satisfy the desired properties.
 This finishes the induction step.
This completes the proof.
\vspace{3mm} \\
{\bf 4.C Stability and  extensions:}
  Let us itroduce  stability under iterations of maps and study their long-time bahaviours.
  
  Let $Q,S \subset {\mathbb R}$ be finite sets,
and consider  a map   $\varphi: Q \times S^a \to S$
with its extension $\varphi: {\mathbb R} \times {\mathbb R}^a \to {\mathbb R}$.

\begin{df} 
$\varphi: {\mathbb R} \times {\mathbb R}^a \to {\mathbb R}$
is  stable with respect to $(Q,S)$,  if there are $0< \delta < 1$
and $0 \leq  \mu < 1$ so that
 the estimates hold:
$$ |\varphi(y, \bar{x})-\varphi(q,\bar{s})| \leq \mu \ d((y, \bar{x}) , (q,\bar{s}))$$
for any $(y, \bar{x}) \in  {\mathbb R} \times {\mathbb R}^a$ and $(q,\bar{s}) \in Q \times S^a$
with $d((y, \bar{x}) , (q,\bar{s})) < \delta$.
  \end{df}
Notice that $\varphi$ is $(\delta', \mu)$-stable for all $0<\delta' \leq \delta$,  whenever
it is  $(\delta, \mu)$-stable.

  \vspace{3mm}
    Let ${\bf A}$ be  an automaton given by a pair of two functions
$ \psi: Q \times S^{\alpha+1} \to S$ and $ \phi: Q \times S^{\beta+1} \to Q$.
A {\em stable extension} of ${\bf A}$ is given by two stable extensions of the functions:
$$\psi: {\mathbb R} \times {\mathbb R}^{\alpha+1} \to {\mathbb R}, \quad
\varphi:  {\mathbb R} \times {\mathbb R}^{\beta+1} \to {\mathbb R}.$$
  \vspace{3mm} \\
{\em Examples 4.1:}
(1) 
 $\psi: {\mathbb R} \times {\mathbb R}^{\alpha+1} \to {\mathbb R}$ 
 is a stable extension of $ \psi: Q \times S^{\alpha+1} \to S$ with $\mu=0$,
if it is locally constant on small neighbourhoods of each  $(q,\bar{s}) \in Q \times S^{\alpha+1}$.

(2)  $\psi: {\mathbb R} \times {\mathbb R} \to {\mathbb R}$ 
 is a stable extension of $ \psi: Q \times S \to S$, 
if it is of the form with $|\alpha_{q,s}| + |\beta_{q,s}| <1$
on small neighbourhoods of each  $(q,s) \in Q \times S$:
$$\psi(y, x) = \alpha_{q,s}(x-s)+\beta_{q,s}(y-q) + \psi(q,s) .$$

\begin{prop}
Let $Q,S \subset {\mathbb R}$ be finite sets and $(\psi, \phi)$ be an  automaton:
$$\psi: Q \times S \to S, \quad
\phi: Q \times S \to S.$$
Then there is a bounded  and stable extension: 
$$\tilde{\psi}: {\mathbb R} \times {\mathbb R} \to {\mathbb R}, \quad
\tilde{\phi}:  {\mathbb R} \times {\mathbb R}  \to {\mathbb R}$$
so that both functions $\tilde{\psi}$ and $ \tilde{\phi}$ 
are represented by relative $(\max, +)$-functions.
\end{prop}
{\em Proof:}
The proof consists of  a minor modification of  lemma $4.2$.

Let us choose disjoint union of squares  in ${\mathbb R}^2$ 
so that each square contains unique point $(q,s) \in Q \times S$ in its interior.
Let us denote such square by $D(q,s)\subset {\mathbb R}^2$.

Firstly let us put locally constant functions by:
$$\tilde{\psi}(y,x)= \psi(q,s), \quad
\tilde{\phi}(y,x)= \phi(q,s)$$
for $(x,y) \in D(q,s)$.
This determines functions 
$\tilde{\psi}$ and $\tilde{\phi}$ on $\cup_{(q,s) \in Q \times S} D(q,s)$.

For the rest, one can follow the same argument as the proof of lemma $4.2$, 
and extend the domains of these functions inductively.

We omit the repetition.
This completes the proof.
\vspace{3mm} \\
{\bf 4.D Stablity and state dynamics:}
Let $Q,S \subset {\mathbb R}$ be finite sets, 
and consider an automaton
${\bf A}$ defined by 
$\psi: Q \times S^{\alpha+1} \to S$ and $  \phi: Q \times S^{\beta+1} \to Q$.
For each $q \in Q$, let
${\bf A}_q: X_S \to X_S$ be  the continuous dynamics
given by:
\begin{align*}
& {\bf A}_q(k_0,k_1, \dots) =(k_0',k_1', \dots), \\
& k_i' = \psi(q_i,k_i, \dots, k_{\alpha+1}), \quad q_{i+1}=\phi(q_i,k_i, \dots, k_{\beta+1})
\end{align*}
with $q_0=q$. 
For a sequence $\bar{q}^l= (q^0,q^1, \dots, q^l)$, let us denote the 
compositions of the corresponding dynamics
${\bf A}_{\bar{q}^l}\equiv {\bf A}_{q^l} \circ \dots \circ {\bf A}_{q^0}: X_S \to X_S$.

\vspace{3mm}

 Let us introduce  the space of the real sequences:
 $$X_{\mathbb R}= \{ \bar{x}= (x_0,x_1, \dots):  \ x_i \in {\mathbb R}\}$$
and equip with  the uniform distance between 
two elements $\bar{x} =(x_0,x_1, \dots)$ and $  \bar{x}' =(x_0',x_1', \dots) \in X_{\mathbb R}$ by:
$$d(\bar{x}, \bar{x}')= \sup_{0 \leq i <\infty}   |x_i -x_i'| \ \in \  [0, \infty]$$

Let:
 $$\psi: {\mathbb R} \times {\mathbb R}^{\alpha+1} \to {\mathbb R}, \quad
\phi:  {\mathbb R} \times {\mathbb R}^{\beta+1} \to {\mathbb R}$$ 
be a stable extension  of ${\bf A}$ with the constants $(\delta, \mu)$. 
For each $y \in {\mathbb R}$, let us extend 
the above dynamics as
${\bf A}_y : X_{\mathbb R} \to X_{\mathbb R}$
given by the same rule:
\begin{align*}
& {\bf A}_y(x_0,x_1, \dots) =(x_0',x_1', \dots), \\
& x_i' = \psi(y_i,x_i, \dots, x_{\alpha+1}), \quad y_{i+1}=\phi(y_i,x_i, \dots, x_{\beta+1})
\end{align*}
with $y_0=y$. 
For $\bar{y}^l= (y^0,y^1, \dots, y^l) \in {\mathbb R}^{l+1}$, we put:
$${\bf A}_{\bar{y}^l}\equiv {\bf A}_{y^l} \circ \dots \circ {\bf A}_{y^0}: X_{\mathbb R} \to X_{\mathbb R}.$$
This is exactly the state dynamics in $2.B$.

For finite sequences $\bar{q}^l \in X_Q^{l+1}$ and $\bar{y}^l \in X_{\mathbb R}^{l+1}$,
we also equip  the same uniform norms by
$d(\bar{q}^l, \bar{y}^l)= \sup_{0 \leq i \leq l} |q^j - y^j|\in [0, \infty)$.

\begin{lem}
Let us choose pairs of the sequences,
  $\bar{y}^l \in {\mathbb R}^{l+1}$ with   $\bar{q}^l \in X_Q^{l+1}$, and
  $\bar{x} \in X_{\mathbb R}$ with $\bar{k} \in X_S$, 
so that they have bounded distances by $0< \delta <1$  from each other:
$$d(\bar{x}, \bar{k}) , \  \ 
d(\bar{y}, \bar{q}) \  < \ \ \delta .$$

Then the estimates hold:
$$d({\bf A}_{\bar{q}}(\bar{k}) , {\bf A}_{\bar{y}}(\bar{x})) \ < \ \delta.$$
  \end{lem}
  {\em Proof:} 
We split the proof into two steps.

  {\bf Step 1:}
    Firstly let us consider the case $l=0$, and put $y_0=y^0 (=y^0_0)$.

    $d(q^0,y^0) <   \delta$ hold by the assumption.
    So let us verify  the estimates:
  $$ d(q_{i+1},y_{i+1}) , \quad d(k_i',x_i')   \  <  \delta$$
    by induction on  $i=0, 1,2, \dots$

Let us start from the former estimates,
and   suppose they hold up to $i$.
Then we have  the estimates:
\begin{align*}
d(q_{i+1}, y_{i+1}) &  = d( \phi(q_i,k_i, \dots, k_{\beta+i}),  \phi(y_i,x_i, \dots, k_{\beta+i})) \\
  &  \leq \ \mu  \  \max_{ i \leq j \leq \beta+i} \{ d(q_i,y_i), d(k_j , x_j)\}  < \mu \delta < \delta.
  \end{align*}
So it holds also at $i+1$. Thus 
$ d(q_i,y_i)     <  \delta$ hold for all $i$ by induction.

  Next let us consider the latter estimates.
  By use of the former ones, we have the desired estimates:
    \begin{align*}
 d(k_i' , x_i') & = d( \psi(q_i,k_i, \dots, k_{\alpha+i}),  \psi(y_i,x_i, \dots, k_{\alpha+i})) \\
&     \leq \ \mu  \    \max_{ i \leq j \leq \alpha+i}(d(q_i,y_i) , d(k_j , x_j) )  <  \delta. \\
     \end{align*}
  Thus we obtain  the bounds:
    $$ d({\bf A}_{q^0}(\bar{k}) , {\bf A}_{y^0}(\bar{x})) \ \ < \delta.$$

  {\bf Step 2:} 
    Let us replace the pairs
 $(\bar{k} , \bar{x})$ by 
  $({\bf A}_{q^0}(\bar{k}) , {\bf A}_{y^0}(\bar{x}))$ and
 $(q^0,y^0)$ by $(q^1,y^1)$
   in  Step $1$. We  apply the same process and obtain the estimates:
    $$d({\bf A}_{q^1}({\bf A}_{q^0}(\bar{k})) , {\bf A}_{y^1}({\bf A}_{y^0}(\bar{x}))) \ < \
  \delta.$$
By iterating this process $l$ times, one obtains the conclusion.

This completes the proof.
\vspace{3mm} \\
{\bf 4.D.2 Change of automata structure:}
For $l=1,2$, let:
$${\bf A}_l: \psi^l: Q \times S^{\alpha+1} \to S,\quad  \phi^l: Q \times S^{\beta+1} \to Q$$
be two automata.
For each $q \in Q$, let
$({\bf A}_l)_q: X_S \to X_S$ be  the continuous dynamics
and denote
$({\bf A}_l)_q(k_0,k_1, \dots) =(k_0'(l),k_1'(l), \dots)$, where $k_i'(l)$ are inductively determined by: 
$$ k_i'(l) = \psi^l(q_i(l),k_i(l), \dots, k_{\alpha+1}(l)), 
\quad q_{i+1}(l)=\phi^l(q_i(l),k_i(l), \dots, k_{\beta+1}(l))$$
with $q_0(l)=q$. 
It happens quite often that  two different automata  ${\bf A}_1$  and ${\bf A}_2$
give the same dynamics. In such situation,  the equalities always hold:
$$(k_0'(1),k_1'(1), \dots) =(k_0'(2),k_1'(2), \dots)$$
while the state sequences may differ from each other:
$$(q^0(1),q^1(1), \dots ) \not= (q^0(2),q^1(2),\dots).$$

\begin{df}
Let us say that  ${\bf A}_1$ and ${\bf A}_2$  are equivalent, if
they give the same dynamics in the sense as above.
  \end{df}
  Let $R \subset Q$ be a subset. If the equalities 
  $({\bf A}_1)_q= ({\bf A}_2)_q$ hold for all $q \in R$, then 
  we say that  ${\bf A}_1$ and ${\bf A}_2$ are equivalent over $R$.

\vspace{3mm}

Let us give stable extensions of ${\bf A}_1$ and ${\bf A}_2$ 
by  relative  $(\max,+)$ functions
$(\phi^1, \psi^1)$ and $ (\phi^2, \psi^2)$ respectively.

 Let  us consider  the corresponding state dynamics:
\begin{align*}
& x_l(i,j+1)= \psi^l(y_l(i,j), x_l(i,j), \dots, x_l(i+\alpha,j)), \\
& y_l(i+1,j)= \phi^l(y_l(i,j), x_l(i,j), \dots,x_l(i+\beta,j))
\end{align*}
with the initial values $\bar{x}(l) = \{x_i(l)\}_i$ and $\bar{y}(l) = \{y^j(l)\}_j$ respectively.
For $l=1,2$,
let us denote:
$$\bar{x}^j(l)=(x_l(0,j),x_l(1,j), \dots)  \in X_{\mathbb R}.$$

For an infinite sequence  $\bar{q}= (q^0,q^1, \dots)\in X_Q$,
we denote its restrictions as  $\bar{q}^l= (q^0,q^1, \dots, q^l)\in X_Q^{l+1}$.

\begin{cor} Suppose  the following conditions:
 
 (1)  There exist  subsets $R \subset Q$ so that
 ${\bf A}_1$ and ${\bf A}_2$ are equivalent over $R$.

(2) The initial data satisfy  the uniform estimates:
$$d(\bar{x}(l), \bar{k}) ,  \quad
d(\bar{y}(l), \bar{q}) \  < \ \ \delta  $$
for some  $\bar{k} \in X_S$ and $ \bar{q}  \in X_R$ for $l=1,2$.

Then the uniform estimates hold for all $0  \leq j  < \infty$:
$$d(\ \bar{x}^j(1) ,\  \bar{x}^j(2)\ ) \ < \ 2  \delta.$$
  \end{cor}
  {\em Proof:} 
  The condition (1) implies that the equalities:
   $$({\bf A}_1)_{\bar{q}^l} = ({\bf A}_2)_{\bar{q}^l} :X_S \to X_S \qquad (l=0,1,2, \dots)$$
 hold for all  $\bar{q}^l \in X_R^{l+1}$.
    Let us denote:
  $$ {\bf A}_{q^j} \circ \dots \circ {\bf A}_{q^0} (k_0,k_1, \dots) =(k_0^j, k_1^j, \dots)\equiv \bar{k}^j.$$

By lemma $4.4$, the estimates:
$$d(\bar{x}^j(l), \bar{k}^j)  \ < \ \  \delta$$
hold for $l=1,2$.
So the estimates:
$$ d(\bar{x}^j(1), \bar{x}^j(2))  \ \leq  
  d(\bar{x}^j(1), \bar{k}^j) + d(\bar{x}^j(2), \bar{k}^j)  \ < \ 2  \delta$$
  hold. This completes the proof.  
\vspace{3mm} \\
{\em Remark 4.2:}
In particular for two different stable extensions of the same automaton, 
we can still apply this and obtain uniform estimates between their orbits.
    \vspace{3mm} \\
{\bf 4.D.3 Uniform estimates for stable dynamics:}
 Let: 
$${\bf A}: \psi: Q \times S^{\alpha+1} \to S, \quad  \phi: Q \times S^{\beta+1} \to Q$$
be an automaton, and choose  stable extensions with the constants $(\delta, \mu)$:
$$\psi: {\mathbb R} \times {\mathbb R}^{\alpha+1} \to {\mathbb R}, \quad
\phi:  {\mathbb R} \times {\mathbb R}^{\beta+1} \to {\mathbb R}$$  
represented by $(\max, +)$-functions.
Let   $(\psi_t, \phi_t)$ be the tropical  correspondences  to $(\psi, \phi)$,
and $M$ be the bigger one of the numbers of their components.
 
Let us consider the corresponding systems of the state dynamics:
\begin{align*}
& x(i,j+1)= \psi(y(i,j), x(i,j), \dots, x(i+\alpha,j)), \\
& y(i+1,j)= \phi(y(i,j), x(i,j), \dots,x(i+\beta,j)), \\
& \\
& x'(i,j+1)= \psi_t(y'(i,j), x'(i,j), \dots, x'(i+\alpha,j)), \\
& y'(i+1,j)= \phi_t(y'(i,j), x'(i,j), \dots,x'(i+\beta,j))
\end{align*}
with the same  initial values $\bar{x} = \{x_i\}_i$ and $\bar{y} = \{y^j\}_j$ respectively.

  \begin{lem}
Suppose  $t_0>>1$ satisfies the estimates:
$$ \mu \delta +2 \log_{t_0} M < \delta.$$

Then for all $t \geq t_0$ and
   any initial data with the uniform bounds:
$$d(\bar{x}, \bar{k}) ,  \quad
d(\bar{y}, \bar{q}) \  < \ \  \frac{\delta}{2} $$
for some  $\bar{k} \in X_S$ and  $\bar{q} \in X_Q$,
their orbits satisfy  the  uniform estimates:
$$ |x(i,j)-x'(i,j)|, \ 
  |y(i,j)-y'(i,j)|  < \frac{\delta}{2}.$$
 \end{lem}
{\em Proof:}
{\bf Step 1:}
The first condition  can be realized for sufficiently large $t_0>>1$,
 since $ \mu <1$ holds by stability assumption.

Let us denote the orbits by $(\{k_i^j\}, \{q^j_i\})$ determined by:
$$ k_i^{j+1}= \psi(q^j_i, k_i^j, \dots, k_{i+\alpha}^j), \ \
 q_{i+1}^j= \phi(q_i^j, k_i^j, \dots,k_{i+\beta}^j)$$
with the initial values $\bar{k}$ and $\bar{q}$ as above.
By lemma $4.4$, the uniform estimate:
$$|x(i,j)-k_i^j|, \quad |y(i,j)-q_i^j| < \frac{\delta}{2}$$
hold for all $i,j$.

\vspace{3mm}

    Let us consider the case $y(i,0)$ and $x(i,1)$.

    $y(0,0)=y'(0,0)$ hold by the assumption.
    So let us verify  the estimates:
  $$ d(y'(i+1,0),y(i+1,0)) , \quad d(x(i,1),x'(i,1))   \  < \frac{ \delta}{2}$$
    by induction on  $i=0, 1,2, \dots$
So  suppose they hold up to $i$.
Notice the estimates:
  $$d(y'(i,0), q_i^0) \leq d(y(i,0), q_i^0) + d(y'(i,0), y(i,0)) < \delta.$$
Then we have the estimates:
   \begin{align*}
&  d(y(i+1,0) , y'(i+1,0))  = 
  d(\phi(y(i,0), x_i, \dots, x_{i+\alpha}),\phi_t(y'(i,0), x_i, \dots, x_{i+\alpha})) \\
&     \leq \ 
d(\phi(y(i,0), x_i, \dots, x_{i+\alpha}),\phi(y'(i,0), x_i, \dots, x_{i+\alpha}))\\
& \qquad + d(\phi(y'(i,0), x_i, \dots, x_{i+\alpha}),\phi_t(y'(i,0), x_i, \dots, x_{i+\alpha}))\\
& \leq \mu d(y(i,0),y'(i,0)) + \log_t M
       < \mu \frac{\delta}{2} + \log_t M < \frac{\delta}{2}.
           \end{align*}
So we have verified the estimates $d(y'(i,0),y(i,0))  <  \frac{\delta}{2}$ for all $i$ by induction.

Then  we have the estimates:
   \begin{align*}
&  d(x(i,1) , x'(i,1))  = 
  d(\psi(y(i,0), x_i, \dots, x_{i+\alpha}),\psi_t(y'(i,0), x_i, \dots, x_{i+\alpha})) \\
&     \leq 
  d(\psi(y(i,0), x_i, \dots, x_{i+\alpha}),\psi(y'(i,0), x_i, \dots, x_{i+\alpha}))  \\
&\qquad +   d(\psi(y'(i,0), x_i, \dots, x_{i+\alpha}),\psi_t(y'(i,0), x_i, \dots, x_{i+\alpha}))\\
& \leq \mu d(y(i,0),y'(i,0))+ \log_t M        < \mu   \frac{\delta}{2} + \log_t M < \frac{\delta}{2}.
           \end{align*}
So we have verified the estimates $d(x'(i,1),x(i,1))  <  \frac{\delta}{2}$.

\vspace{3mm}

{\bf Step 2:} 
Let us put the sequences $\bar{x}^j \equiv (x(0,j),x(1,j), \dots)$ and similar for others.
 Let us verify  the estimates:
  $$d((\bar{y}')^j, \bar{y}^j)  ,   \ \ d((\bar{x}')^{j+1}, \bar{x}^{j+1}) \quad   \  < \frac{ \delta}{2}$$
    by induction on  $j=0, 1,2, \dots$

We are done for $j=0$ at Step $1$.
Suppose the above estimates hold up to $j-1$. 

Let us start from the former estimates.
$y(0,j)=y'(0,j)$ hold by the assumption.
Suppose  the estimates $d(y'(i,j),y(i,j))  <  \frac{\delta}{2}$
 hold up to $i$.
Then we have  the estimates:
   \begin{align*}
&  d(y(i+1,j) , y'(i+1,j))  = \\
&  d(\phi(y(i,i), x(i,j), \dots, x(i,i+\beta)),\phi_t(y'(i,j), x'(i,j), \dots, x'(i+\beta)) \\
&     \leq \ 
 d(\phi(y(i,j), x(i,j), \dots, x(i,i+\beta)),\phi(y'(i,j), x'(i,j), \dots, x'(i+\beta)) \\
 &+  d(\phi(y'(i,j), x'(i,j), \dots, x'(i,i+\beta)),\phi_t(y'(i,j), x'(i,j), \dots, x'(i+\alpha)) \\
& \leq \mu  \  \max_{ i \leq l \leq \beta+i} \{ d(y(i,j),y'(i,j)), d(x(l,j),x'(l,j))\} 
      + \log_t M   \\
& < \mu \frac{\delta}{2} + \log_t M < \frac{\delta}{2}.
           \end{align*}
So we have verified the estimates 
$d((\bar{y}')^j, \bar{y}^j)  <  \frac{\delta}{2}$ for all $i$ by induction.

Then  we have the estimates:
 \begin{align*}
&  d(x(i,j+1) , x'(i,j+1))  = \\
&  d(\psi(y(i,j), x(i,j), \dots, x(i+\alpha,j)),\psi_t(y'(i,j), x'(i,j), \dots, x'(i+\alpha,j))) \\
&     \leq \ 
  d(\psi(y(i,j), x(i,j), \dots, x(i+\alpha,j)),\psi(y'(i,j), x'(i,j), \dots, x'(i+\alpha,j))) \\
& +  d(\psi(y'(i,j), x'(i,j), \dots, x'(i+\alpha,j)),\psi_t(y'(i,j), x'(i,j), \dots, x'(i+\alpha,j))) \\
      & \leq \mu   \  \max_{ i \leq l \leq \alpha+i} \{ d(y(i,j),y'(i,j)), d(x(l,j),x'(l,j))\} 
      + \log_t M \\
      &  < \mu \frac{\delta}{2} + \log_t M < \frac{ \delta}{2}.
           \end{align*}
So we have verified the estimates $ d((\bar{x}')^{j+1}, \bar{x}^{j+1})  <  \frac{\delta}{2}$.

So we have finished the induction step on $j$.
  This completes the proof.
\vspace{3mm} \\
{\bf 4.D.4 Rational dynamics and  change of automata:}
  Let $S \subset {\mathbb R}$ be a finite set, and take 
  $\bar{a}=(a_0,a_1, \dots) \in X_S$.
  The {\em exponential sequence}  
 is given by the sequence of positive numbers parametrized by $t >1$:
 $$t^{\bar{a}}=(t^{a_0},t^{a_1}, \dots).$$
 Its $C\geq 1 $ neighbourhood is given by the set:
 $$N_C( t^{\bar{a}})= \{\bar{z}  \in X_{\mathbb R}: 
C^{-1}  t^{a_i} < z_i  < C t^{a_i} \} \ \ \supset \ \ t^{\bar{a}}.$$

\vspace{3mm}

For $l=1,2$,  let:
$${\bf A}_l: \psi^l: Q \times S^{\alpha+1} \to S, \quad  \phi^l: Q \times S^{\beta+1} \to Q$$
be two automata, and choose their stable extensions: 
$$\psi^l: {\mathbb R} \times {\mathbb R}^{\alpha+1} \to {\mathbb R}, \quad
\phi^l:  {\mathbb R} \times {\mathbb R}^{\beta+1} \to {\mathbb R} \qquad (l=1,2)$$
with the constants $(\delta, \mu)$ by relative $(\max,+)$-functions.

\vspace{3mm}

Let  $(f_t^l,g_t^l)$ and $(\psi^l_t, \phi^l_t)$ be the tropical 
 correspondences  to $(\psi^l, \phi^l)$ respecively,
and  $M$ be the bigger one of their  numbers of the components.

Let us   consider  the state systems of the rational dynamics:
 \begin{align*} 
& z_i^{j+1}(l)= f_t^l(w_i^j(l), z_i^j(l), \dots,   z_{i+ \alpha}^j(l)) , \\
&  w_{i+1}^j(l)= g_t^l(w_i^j(l),   z_i^j(l),  \dots, z_{i+ \beta}^j(l))
 \end{align*}
with the initial values $\bar{z}(l)= \{z_i(l)\}_i $ and $\bar{w}(l)= \{w^j(l)\}_j$.

\vspace{3mm}

\begin{thm} Suppose ${\bf A}_1$ and ${\bf A}_2$ are equivalent over some $R \subset Q$,
and choose stable extensions with the constants $( \delta, \mu)$.

Then for any large $C >> 1$, there exists $t_0 >1$ so that 
for all  $t \geq t_0$ and any  initial values which are  are contained in:
  $$\bar{z}(l) \in N_C( t^{\bar{k}}), \quad  \bar{w}(l) \in  
  N_C( t^{\bar{q}}) \qquad (l=1,2)$$ 
for some $\bar{k} \in X_S$ and $\bar{q} \in X_R$,
then the uniform estimates hold:
$$\max \{  \frac{z_i^j(1)}{z_i^j(2)} , \  \frac{z_i^j(2)}{z_i^j(1)} \}
 <  C^4.$$
 \end{thm}
{\em Proof:} We split the proof into two steps.

{\bf Step 1:}
Notice that the pairs
$(\psi^l,  \phi^l)$ are  $(\delta', \mu)$-stable for all $0<\delta' \leq \delta$.
Let us choose large $C$ so that the estimates hold:
$$M < C^{1-   \mu}.$$
Then  choose  $t_0 >>1$ so that  the estimates hold: 
$$ \log_{t_0} C  \leq \frac{ \delta}{2}.$$
Now let us choose and fix any $t  \geq t_0$, and put
$ \delta' =  \log_t C$.
Then the estimates:
$$ \mu \delta'  +  \log_{t_0} M \leq  \mu \log_t C +\log_t M   < \log_t C = \delta'$$
hold by the above inequality. So the condition in lemma $4.6$ is satisfied.

Let us regard that the  pairs $(\psi^l,\phi^l)$ are $(2\delta', \mu)$-stable.

\vspace{3mm}

{\bf Step 2:}
Let us put:
$$x_i(l)= \log_t z_i(l), \quad y^j(l)=\log_t w^j(l) \qquad (l=1,2)$$
and consider the corresponding systems of the state dynamics:
\begin{align*}
& x_l(i,j+1)= \psi^l(y_l(i,j), x_l(i,j), \dots, x_l(i+\alpha,j)), \\
& y_l(i+1,j)= \phi^l(y_l(i,j), x_l(i,j), \dots,x_l(i+\beta,j)), \\
& \\
& x'_l(i,j+1)= \psi^l_t(y'_l(i,j), x'_l(i,j), \dots, x'_l(i+\alpha,j)), \\
& y'_l(i+1,j)= \phi^l_t(y'_l(i,j), x'_l(i,j), \dots,x'_l(i+\beta,j))
\end{align*}
with the initial values $\bar{x}(l) = \{x_i(l)\}_i$ and $\bar{y}(l) = \{y^j(l)\}_j$ respectively.

Then the estimates:
$$d(\bar{x}(l) , \bar{k}) = \sup_i |x_i(l)-k_i| 
=\sup_i  \log_t ( \frac{z_i(l)}{t^{k_i}})^{\pm 1} <  \log_t C =\delta'$$
hold. Similarly the estimates $d(\bar{y}(l) , \bar{q})  < \delta'$ hold.

By corollary $4.5$, the estimates:
$$|x_1(i,j) - x_2(i,j)| < 2 \delta'$$ hold for all $i,j$. 

On the other hand by lemma $4.6$, 
$$|x_l(i,j) - x_l'(i,j)| <  \delta'$$
 hold for $l=1,2$ and all $i,j$. 
 
 Then combining with these, we have the estimates:
\begin{align*}
 | x_1'(i,j)- x_2'(i,j)|  & \leq  | x_1'(i,j)- x_1(i,j)| +  \\
 &  \qquad | x_1(i,j)- x_2(i,j)|  + | x_2(i,j)- x_2'(i,j)| \\
& < 4 \delta' = \log_t C^4.
\end{align*}
Since the equalities:
$$ | x_1'(i,j)- x_2'(i,j)| = \log_t (\frac{z_i^j(1)}{z_i^j(2)})^{\pm 1}$$
 hold, this verifies the desired estimates.
 This completes the proof.
\vspace{3mm} \\
{\bf 4.D.5 Dynamical inequalities under change of automata:}
For $l=1,2$, let:
$${\bf A}_l: \psi^l: Q \times S^{\alpha+1} \to S, \quad  \phi^l: Q \times S^{\beta+1} \to Q$$
be automata which are  equivalent over $R \subset Q$.

For each $l$, 
let us  take two pairs of their 
  stable extensions with the constants $(\delta, \mu)$
by $(\max,+)$-functions:
$$\psi^{l,m}: {\mathbb R} \times {\mathbb R}^{\alpha+1} \to {\mathbb R}, \quad
\phi^{l,m}:  {\mathbb R} \times {\mathbb R}^{\beta+1} \to {\mathbb R} \quad (m=1,2)$$
Suppose that for each $l=1,2$:
$$(\psi^{l,1}, \phi^{l,1}) \sim (\psi^{l,2}, \phi^{l,2})$$
 are pairwisely  tropically equivalent. 
Let  $(f_t^{l,m},g_t^{l,m})$ and $(\psi^{l,m}_t, \phi^{l,m}_t)$ be the tropical 
 correspondences  to $(\psi^{l,m}, \phi^{l,m})$ respectively.

\begin{cor} Assume the above conditions.
Then for any large $C \geq 1$, there exists $t_0 >1$ and $D\geq 1$ so that 
the following holds for all  $t \geq t_0$:

Suppose 
 two sequences 
$\{w_i^j(l)\}_{i,j}$ and $\{z_i^j(l)\}_{i,j}$ satisfy the dynamical 
inequalities:
\begin{align*}
& f^{l,1}_t(w_i^j(l), z_i^j(l), \dots, z_{i+\alpha}^j(l)) \leq z_i^{j+1}(l) \leq  
f^{l,2}_t(w_i^j(l), z_i^j(l), \dots, z_{i+\alpha}^j(l)) , \\
& g^{l,1}_t(w_i^j(l),   z_i^j(l), \dots,z_{i+\beta}^j(l)) \leq w_{i+1}^j(l) \leq
 g^{l,2}_t(w_i^j(l),   z_i^j(l), \dots , z_{i+\beta}^j(l))
\end{align*}
Moreover suppose the initial values are contained in the same $C$ neighbourhoods:
  $$\bar{z}(l)=(z_0(l), z_1(l), \dots)  \in N_C( t^{\bar{k}}), 
  \quad  \bar{w}(l)=(w^0(l),w^1(l), \dots)  \in  
  N_C( t^{\bar{q}}) $$ 
for some $\bar{k} =(k_0,k_1, \dots) \in X_S$ and $\bar{q} =(q^0,q^1, \dots) \in X_R$.

Then the uniform estimates hold:
$$\max \{  \frac{z_i^j(1))}{z_i^j(2)} , \  \frac{z_i^j(2)}{z_i^j(1)} \}
 <  D.$$
\end{cor}
{\em Proof:}
The proof is long and we split it into several steps,
but the idea is quite parallel to theorem $2.19$ and corollary $2.21$.

Let us choose and fix large $t>>1$.
By replacing   $\delta$ by a smaller one as in Step $1$ in the proof of theorem $4.7$,
we may assume the followings:

\vspace{2mm}

(1) $\delta = \log_t C$ and the estimates:
$$M^2 < C^{1- \mu}$$
holds. In particular  the estimates
$2\log_t M + \mu \delta < \delta$ holds.

\vspace{2mm}

(2) $(\psi^{l,m}, \phi^{l,m})$ are all $(3\delta,  \mu)$- stable.

\vspace{3mm}

Let us put $p_i^j(l) = \log_t z_i^j(l)$ and $o_i^j(l)= \log_t w_i^j(l)$.
Then we have  the estimates:
\begin{align*}
& \psi^{l,1}_t(o_i^j(l), p_i^j(l), \dots, p_{i+\alpha}^j(l)) \leq p_i^{j+1}(l)
 \leq  \psi^{l,2}_t(o_i^j(l), p_i^j(l), \dots,p_{i+\alpha}^j(l)) , \\
& \phi^{l,1}_t(o_i^j(l),   p_i^j(l), \dots, p_{i+\beta}^j(l)) \leq o_{i+1}^j(l)
 \leq \phi^{l,2}_t(o_i^j(l),   p_i^j(l),  , \dots, p_{i+\beta}^j(l))
\end{align*}
where $p_i^0(l)=x_i(l) \equiv \log_t z_i(l)$ and $o_0^j(l)=y^j(l) \equiv \log_t w^j(l)$.

Notice that  both the equalities:
$$\psi^{l,1}=\psi^{l,2}, \quad \phi^{l,1}=\phi^{l,2}$$ hold as functions
 by the assumption.
We will use the notation $\psi^l$ and $\phi^l$ when no cofusion occurrs.

Let us consider  another  solutions to the state systems:
\begin{align*}
& k_i^{j+1}(l)
= \psi^l(q_i^j(l), k_i^j(l), \dots, k_{i+\alpha}^j(l))  \\
& q_{i+1}^j(l) = \phi^l(q_i^j(l),   k_i^j(l),  , \dots, k_{i+\beta}^j(l))
\end{align*}
with $k_i^0(l)=k_i $ and $q_0^j(l)=q^j  $.

We verify   the uniform estimates: 
$$|p_i^j(l)-k_i^j(l)|, \quad |o_i^j(l) - q_i^j(l)| < \delta$$
by several steps.

\vspace{3mm}

{\bf Step 1:}
Let us introduce another state systems:
 \begin{align*}
 &  x_{l,m}'(i,j+1)= \psi_t^{l,m}(y_{l,m}'(i,j), x_{l,m}'(i,j), \dots, x_{l,m}'(i+\alpha,j)),  \\
&  y_{l,m}'(i+1,j)= \phi_t^{l,m}(y_{l,m}'(i,j), x_{l,m}'(i,j), \dots, x_{l,m}'(i+\beta,j)), \\
& \\
&  x_l(i,j+1)= \psi^l(y_l(i,j), x_l(i,j), \dots, x_l(i+\alpha,j)),  \\
&  y_l(i+1,j)= \phi^l(y_l(i,j), x_l(i,j), \dots, x_l(i+\beta,j)), \\
\end{align*}
with the same initial values 
 $x_{l,m}'(i,0)=x_l(i,0)=x_i(l)$ and $y_{l,m}'(0,j)=y_l(0,j)= y^j(l)$.

By lemma $4.4$, the estimates hold:
$$|x_l(i,j) - k_i^j(l)|, \quad |y_l(i,j)-q_i^j(l)| <  \delta.$$

On the other hand by lemma $4.6$, another estimates hold:
$$|x_l(i,j)-x'_{l,m}(i,j)|, \quad |y_l(i,j)-y'_{l,m}(i,j)| < \delta.$$

So combining with these estimates, one obtains the estimates:
$$|x_{l,m}'(i,j) - k_i^j(l)|, \quad |y_{l,m}'(i,j)-q_i^j(l)| < 2 \delta.$$

\vspace{3mm}

{\bf Step 2:}
Next let us  verify  the uniform estimates for $l=1,2$:
$$|x_{l,m}'(i,1)-p_i^1(l)|, \  |y_{l,m}'(i,0)- o_i^0(l)| \  < \delta . $$
Notice the equalities   $y_{l,m}'(0,0)=o^0_0(l)$ and $x_{l,m}'(i,0)= p_i^0(l)$.

   Let us verify 
  $ |y_{l,m}'(i,0)- o_i^0(l)| \  < \delta  $
    by induction on  $i=0, 1,2, \dots$
    Assume that the estimates hold up to $i$.

Firstly suppose $o_{i+1}^0(l) \geq  y_{l,1}'(i+1,0)$ hold. 
Then the estimates hold:
\begin{align*}
& 0  \leq  o_{i+1}^0(l) -  y_{l,1}'(i+1,0)
=| o_{i+1}^0(l) -  y_{l,1}'(i+1,0)| \\
& \leq |\phi_t^{l,2}(o_i^0(l), p_i^0(l), \dots,p_{i+\beta}^0(l)) -
\phi_t^{l,1}(y_{l,1}'(i,0), x_{l,1}'(i,0), \dots, x_{l,1}'(i+\beta,0))| \\
& \leq 2 \log_t M + \\
&  |\phi^l(o_i^0(l), p_i^0(l), \dots,p_{i+\beta}^0(l)) -
\phi^l(y_{l,1}'(i,0), x_{l,1}'(i,0), \dots, x_{l,1}'(i+\beta,0))| \\
& \leq 2 \log_t M  + 
 \mu 
|o_i^0(l) - y_{l,1}'(i,0)| 
 \leq 2\log_t M + \mu \delta < \delta.
\end{align*}
Conversely suppose   $o_{i+1}^0(l) \leq  y_{l,1}'(i+1,0)$ hold. 
 Then the estimates hold:
\begin{align*}
& 0  \leq    y_{l,1}'(i+1,0) -  o_{i+1}^0(l) 
=|  y_{l,1}'(i+1,0)  - o_{i+1}^0(l)   | \\
& \leq |\phi_t^{l,1}(y_{l,1}'(i,0), x_{l,1}'(i,0), \dots, x_{l,1}'(i+\beta,0)) - 
\phi_t^{l,1}(o_i^0(l), p_i^j(l), \dots,p_{i+\beta}^0(l))| \\
& \leq 2 \log_t M + \\
&  |\phi^l(y_{l,1}'(i,0), x_{l,1}'(i,0), \dots, x_{l,1}'(i+\beta,0)) - 
\phi^l(o_i^0(l), p_i^0(l), \dots,p_{i+\beta}^0(l))| \\
& \leq 2 \log_t M  + \mu 
|o_i^0(l) - y_{l,1}'(i,0)| \leq 2\log_t M + \mu \delta < \delta.
\end{align*}
Thus in any cases, the estimates
    $ |y_{l,1}'(i+1,0)- o_{i+1}^0(l)| \  < \delta  $  hold.
      By the induction step, we have verified the claim.
    By the same way  the estimates
      $ |y_{l,2}'(i,0)- o_i^0(l)| \  < \delta  $ also hold for all $i$.

Then  by use of the above estimates, 
we follw a parallel argument as below.
Suppose $p_i^1(l) \geq x_{l,1}'(i,1)$ hold for some $i$. Then we have the estimates:
\begin{align*}
& 0  \leq  p_i^1(l) -  x_{l,1}'(i,1)
=| p_i^1(l) -  x_{l,1}'(i,1)| \\
& \leq |\psi_t^{l,2}(o_i^0(l), p_i^0(l), \dots,p_{i+\alpha}^0(l)) -
\psi_t^{l,1}(y_{l,1}'(i,0), x_{l,1}'(i,0), \dots, x_{l,1}'(i+\beta,0))| \\
& \leq 2 \log_t M + \\
&  |\psi^l(o_i^0(l), p_i^0(l), \dots,p_{i+\alpha}^0(l)) -
\psi^l(y_{l,1}'(i,0), x_{l,1}'(i,0), \dots, x_{l,1}'(i+\beta,0))| \\
& \leq 2 \log_t M  +
\mu  
|o_i^0(l) - y_{l,1}'(i,0)|  \leq 2\log_t M + \mu \delta < \delta.
\end{align*}
Conversely suppose    $p_i^j(l) \leq x_{l,1}'(i,1)$ hold. 
 Then we have  the estimates:
 \begin{align*}
& 0  \leq   x_{l,1}'(i,1)- p_i^1(l)
=| x_{l,1}'(i,1)- p_i^1(l) | \\
& \leq |\psi_t^{l,1}(y_{l,1}'(i,0), x_{l,1}'(i,0), \dots, x_{l,1}'(i+\beta,0))-
\psi_t^{l,1}(o_i^0(l), p_i^0(l), \dots,p_{i+\alpha}^0(l)) | \\
& \leq 2 \log_t M + \\
&  |
\psi^l(y_{l,1}'(i,0), x_{l,1}'(i,0), \dots, x_{l,1}'(i+\beta,0))
- \psi^l(o_i^0(l), p_i^0(l), \dots,p_{i+\alpha}^0(l)) | \\
& \leq 2 \log_t M  +  \mu 
|o_i^0(l) - y_{l,1}'(i,0)| 
 \leq 2\log_t M + \mu \delta < \delta.
\end{align*}
Thus in any cases  we have the estimates 
$|x_{l,1}'(i,1)-p_i^1(l)|   < \delta  $ for all $i=0,1,2, \dots$
The estimates $|x_{l,2}'(i,1)-p_i^1(l)|   < \delta  $ are obtained by the same way.

       \vspace{3mm}

       {\bf Step 3:} Let us verify the estimates for all $i,j  \geq 0$:
       $$|x_{l,m}'(i,j)-p_i^j(l)|, \  |y_{l,m}'(i,j)- o_i^j(l)| \  < \delta  $$
       
       Let us put the sequences:
       $$\bar{x}'_{l,m}(j) =(x'_{l,m}(0,j),x'_{l,m}(1,j),x'_{l,m}(2,j), \dots)$$
       and similar for    $\bar{y}'_{l,m}(j)$, $\bar{p}^j(l)$ and $\bar{o}^j(l)$.
       
       Let us verify the estimates: 
       $$d(\bar{x}_{l,m}'(j+1), \bar{p}^{j+1}(l)), \  d(\bar{y}_{l,m}'(j) , \bar{o}^j(l)) \  < \delta  $$
       by induction on $j=0,1,2,\dots$
We have verified the estimates for $j=0$ at step $2$.
So assume they hold up to $j-1$.

The initial conditions
$o_0^j(l)=  y_{l,1}'(0,j)$ hold. Let us verify the estimates
$|o_i^j(l)-   y_{l,1}'(i,j)| < \delta$ by induction on $i$.
    Suppose they hold up to $i$.
    
Firstly suppose $o_{i+1}^j(l) \geq  y_{l,1}'(i+1,j)$ hold. 
Then the estimates hold:
\begin{align*}
& 0  \leq  o_{i+1}^j(l) -  y_{l,1}'(i+1,j)
=| o_{i+1}^j(l) -  y_{l,1}'(i+1,j)| \\
& \leq |\phi_t^{l,2}(o_i^j(l), p_i^j(l), \dots,p_{i+\beta}^j(l)) -
\phi_t^{l,1}(y_{l,1}'(i,j), x_{l,1}'(i,j), \dots, x_{l,1}'(i+\beta,j))| \\
& \leq 2 \log_t M + \\
&  |\phi^l(o_i^j(l), p_i^j(l), \dots,p_{i+\beta}^j(l)) -
\phi^l(y_{l,1}'(i,j), x_{l,1}'(i,j), \dots, x_{l,1}'(i+\beta,j))| \\
& \leq 2 \log_t M  + \\
& \mu  \max(
|o_i^j(l) - y_{l,1}'(i,j)|, 
| p_i^j(l) - x_{l,1}'(i,j)|,  \dots,  | p_{i+\beta}^j(l)  - x_{l,1}'(i+\beta,j)|) \\
& \leq 2\log_t M + \mu \delta < \delta.
\end{align*}
The converse case can be estimated by the same way as step $2$, and we omit repetition.
So by the induction step   we have the estimates
    $ |y_{l,1}'(i,j)- o_i^j(l)| \  < \delta  $ for all $i$.
    By the same way we can verify the estimates
      $ |y_{l,2}'(i,j)- o_i^j(l)| \  < \delta  $. 
      
Then  we have the estimates 
$|x_{l,m}'(i,j+1)-p_i^{j+1}(l)|   < \delta  $ for all $i=0,1,2, \dots$
again by the same way as step $2$.

This completes the induction step on $j$, and so we have obtained 
the desired estimates:
       $$|x_{l,m}'(i,j)-p_i^j(l)|, \  |y_{l,m}'(i,j)- o_i^j(l)| \  < \delta  $$
for all $i,j=0,1,2, \dots$

       \vspace{3mm}

{\bf Step 4:}
Combining step $1\sim 3$,  we obtain the estimates:
$$|p_i^j(l) - k_i^j(l)| < 3 \delta.$$
In particular all the values  of these sequences lie
within $\mu$-Lipschitz constants of
$\psi^{l,m}$ and $\phi^{l,m}$.

Let us consider the state systems of the rational dynamics:
\begin{align*}
& F^l(i,j+1)= f^{l,1}_t(G^l(i,j), F^l(i,j), \dots, F^l(i+\alpha,j)) \\
& G^l(i,j+1)= g^{l,1}_t(G^l(i,j),   F^l(i,j), \dots,F^l(i+\beta,j))
\end{align*}
with the initial values $F^l(i,0)=z_i(l)$ and $G^l(0,j)= w^j(l)$ respectively.
Then we apply    theorem  $2.17$, and obtain   the uniform estimates:
$$
\max( \frac{F^l(i,j+1)}{z_i^{j+1}(l)} ,  \  \frac{z_i^{j+1}(l)}{F^l(i,j+1)}) 
 \leq  \ \  M^{2P_{i+j(\gamma+1)}(\mu)}\leq C'  \qquad (l=1,2)$$
for some $C'$, since
$P_i(\mu)$ are uniformly bounded for $0< \mu <1$.
 On the other hand by theorem $4.7$, the estimates hold:
$$\max(   \frac{F^1(i,j+1)}{F^2(i,j+1)} , \  \frac{F^2(i,j+1)}{F^1(i,j+1)}, \ \leq    C^4.$$

Thus combining with these, one obtains the desired estimates:
$$(\frac{z_i^j(1)}{z_i^j(2)} )^{\pm 1}  
= (\frac{z_i^j(1)}{F^1(i,j)})^{\pm 1} (\frac{F^2(i,j)}{z_i^j(2)} )^{\pm 1}
 \leq  \ \  C' \ C^4 \equiv D. $$
This completes the proof.
\vspace{3mm} \\
{\bf 4.E  Quasi recursivity:}
Let us consider recursivity for rational dynamics. 
We start from a simpler case which is  time-independent and one dimensional.

Let $\varphi$ be a $(\max,+)$ function of $n$ variable, and consider the 
one dimensional dynamics:
$$x_N =\varphi(x_{N-n}, \dots, x_{N-1}) $$
with the initial data $\bar{x}_0=(x_0,\dots, x_{n-1}) \in {\mathbb R}^n$. 
It is {\em recursive} if there is some $M\geq 0$
so that any orbits $\{x_i\}_i$ with any initial values are periodic of period $M$.
Namely equalities $x_{i+M}=x_i$ hold  for all $i \geq 0$.

Let $f_t$ be the tropical correspondence to $\varphi$, and consider the parallel
dynamics:
$$z_N =f_t(z_{N-n}, \dots, z_{N-1})$$
 with the initial data
$x_0 = \log_tz_0, \dots, x_{n-1}= \log_t z_{n-1}$.

It is known that 
$\varphi$ is always recursive whenever $f_t$ is the case for all $t >1$
([K2]).
However the converse is not true.

Let us say that   $f_t$ is a {\em quasi recursive} of period $M$,
if there are constants  $C \geq 1$ and   $M \geq 0$
independently of $t$ and initial values,
 so that   the uniform  estimates hold:
   $$\max (\frac{z_{N+M}}{z_N},  \frac{z_N}{z_{N+M}})  \leq C, \quad 
   N =0,1, \dots$$

   \begin{prop}[K2]
Suppose $f_t$ corresponds to a relative
$(\max, +)$-function $\varphi$. 
Then $f_t$ is  quasi recursive of minimum period $M$, if and only if 
$\varphi$ is recursive of the same minimum  period.
\end{prop}
     
For example:
\begin{align*} & 
\varphi^1(x_0,x_1)= \max(0,x_1)-x_0, \\
& \varphi^2(x_0,x_1)= \max(x_1,-x_1)-x_0
\end{align*}
 are 
recursive of periods $5$ and $9$ respecively.
The corresponding relatively elementary  functions are given by:
\begin{align*}
& f_t^1(w,z)= w^{-1}(1+z), \\
&  f_t^2(w,z)= w^{-1}(z^{-1}+z)
\end{align*}
 respectively.
It turns out that the former is also  recursive, but the second is not the case, and so it only satisfies 
 quasi recursivity. 
\vspace{3mm} \\
{\bf 4.E.2  Finite order group elements:}
Let:
$${\bf A}: \psi: Q \times S^{\alpha+1} \to S, \ \  \phi: Q \times S^{\beta+1} \to Q$$
be an automaton, and 
$\psi: {\mathbb R} \times {\mathbb R}^{\alpha+1} \to {\mathbb R}, \ \
\phi:  {\mathbb R} \times {\mathbb R}^{\beta+1} \to {\mathbb R}$   be   stable extensions
with the constants $(\delta, \mu)$.

Let
$(f_t,g_t)$ and $(\psi_t, \phi_t)$ be the 
 corresponding functions to $(\psi, \phi)$ respectively, 
and $M$ be the biggest one of the numbers of their components.

For $\bar{q}^m =(q^0, \dots, q^m) \in X_Q^{m+1}$, let us denote 
$l$ times iterations of $\bar{q}^m$ by
$l \bar{q}^m \equiv (q^0, \dots, q^m, q^0, \dots ,q^m, \dots,  q^0, \dots, q^m) \in X^{(m+1)l}_Q$.
We also denote the infinite  times iterations of $\bar{q}^m$ by:
$$\bar{q}^m_{per} \equiv (q^0, \dots, q^m, q^0, \dots ,q^m, \dots,  q^0, \dots, q^m, \dots) \in X_Q.$$

 Let us consider the state dynamics:
 $$ z_i^{j+1}= f_t(w_i^j, z_i^j, \dots, z_{i+\alpha}), \quad
 w^j_{i+1}= g_t(w_i^j   ,z_i^j, \dots, z_{i+\beta}^j).$$

\begin{prop}
Suppose 
${\bf A}_{\bar{q}^m}: X_S \to X_S$ is of finite order with period $p $:
$${\bf A}_{p \bar{q}^m} = ({\bf A}_{\bar{q}^m}) \circ  \dots  \circ ({\bf A}_{\bar{q}^m}) 
 \equiv ({\bf A}_{\bar{q}^m})^p = \text{ id }.$$
Then 
for any $C  \geq 1$, there exists $t_0>1$  so that
  for all $t \geq t_0$ and 
any   initial values:
  $$\{z_i\}_i \subset N_C(t^S), \quad \{w^j\}_j \subset N_C(t^{ \bar{q}^m_{per}})$$ 
 the uniform bounds hold for all $i,j,l =0,1,2, \dots$:
$$(   \frac{z_i^j}{z_i^{j+p(m+1)l}} )^{\pm 1} \ \leq    C^4$$
\end{prop}
{\em Proof:} 
Let us choose large $t_0 >1$ so that the estimates hold:
$$\log_{t_0} C, \ \mu  \frac{\delta}{2} + \log_{t_0} M < \frac{\delta}{2}.$$

Recall that the pair $(\psi, \phi)$ is $(\delta', \mu)$-stable for any $0< \delta' \leq \delta$.
Let us fix $t \geq t_0$.
Then by replacing $\delta$ by $\delta' = 2\log_t C$, 
one mya assume the equality $\delta= 2\log_t C$.

By the assumption, there is 
$ \bar{k}=(k_0,k_1,  \dots)  \in X_S$ so that 
the initial value is contained as
$\{z_i\}_i \in N_C(t^{ \bar{k} })$. 
Let us rewrite  $ \bar{q}^m_{per}=(q^0,q^1,\dots)$

Let us consider three state systems of dynamics:
\begin{align*}
& x(i,j+1)= \psi(y(i,j), x(i,j), \dots, x(i+\alpha,j)), \\
& y(i+1,j)= \phi(y(i,j), x(i,j), \dots,x(i+\beta,j)), \\
& x'(i,j+1)= \psi_t(y'(i,j), x'(i,j), \dots, x'(i+\alpha,j)), \\
& y'(i+1,j)= \phi_t(y'(i,j), x'(i,j), \dots,x'(i+\beta,j)), \\
& k_i^{j+1}= \psi(q_i^j, k_i^j, \dots, k_{i+\alpha}^j), \\
& q_{i+1}^j= \phi(q_i^j, k_i^j, \dots, k_{i+\beta}^j),
\end{align*}
with the initial values $x(i,0)=x'(i,0)= \log_t z_i$, $y(0,j)=y'(0,j)= \log_t w^j$,
and $k_i^0=k_i$, $q_0^j=q^j$.
By the condition, 

\vspace{2mm}

(1) The estimates
 $|x(i,0) - k_i|, |y(0,j)-q^j| < \frac{\delta}{2}$ hold.

\vspace{2mm}

(2) periodicity $k_i^j = k_i^{{j+p(m+1)l}}$ hold for all $i,j,l$.

\vspace{3mm}

By lemma $4.4$,  the estimates 
$|x(i,j)-k_i^j| <\frac{\delta}{2}$ hold.
By lemma $4.6$, the estimates
$|x(i,j)-x'(i,j)| < \frac{\delta}{2}$ hold.
Combining with these, we obtain the estimates:
$$|x'(i,j)-k_i^j| < \delta.$$

Then we have the estimates:
\begin{align*}
 |x' & (i,j) -  x'(i, j+p(m+1)l)| \\
&  \leq |x'(i,j)-k_i^j| + |x'(i, j+p(m+1)l)- k_i^{j+p(m+1)}| 
<2 \delta = \log_t C^4.
\end{align*}
Since the left hand side is equal to 
$\log_t (\frac{z_i^j}{z_i^{j+p(m+1)l}})^{\pm 1}$, 
the conclusion holds.

This completes the proof.
\vspace{3mm} \\
{\bf 4.E.3 Infinite quasi-recursive dynamics:} 
The {\em Burnside problem} asks  existence of finitely generated and infinite torsion groups.
The first example was given by Adjan-Novikov ([AN]).
The second one is given by an automata group:
\begin{lem}[Al]
There eixsts a Mealy automaton with $2$ alphabets and $8$ states
such that the group generated by some two states $u,v $ is infinite and torsion.
\end{lem}

See also [Z] for the exposition of this group.

\vspace{3mm}

Let $(f,g)$ be a pair of (unparametrized) rational  functions,
and consider the state system of the rational dynamics
 with the initial sets:
\begin{align*}
&  z_i^{j+1}= f(w_i^j, z_i^j, \dots, z^j_{i+\alpha}), \quad
 w^j_{i+1}= g(w_i^j   ,z_i^j, \dots, z_{i+\beta}^j), \\
& X \subset \{(z_0,z_1, \dots) : z_i \in (0, \infty)\}, \quad
 Y \subset \{(w^0,w^1, \dots) : w^i \in (0, \infty)\}.
\end{align*}

The  system  is said to be {\em recursive}
with respect to the initial data $(X,Y)$,
if for any initial values $\{w^j_0\}_j \in Y$,
 there exist some $p\in {\mathbb N}$
so that 
any solutions $(\{z_i^j\}_{i,j}, \{w_i^j\}_{i,j})$ with the initial values $\{z_i^0\}_i \in X$
 and   $\{w^j_0\}_j$
 satisfy $p$-periodicity  
 $z_i^j = z_i^{j+pl}$ for all $ i,j, l =0,1,2, \dots$

{\em The rational Burnside problem} asks existence of pairs of rational functions
whose state dynamics are recursive with infinitely many such $p$.

So far it seems not known whether such a pair of rational functions exist.

\vspace{3mm} 

Let us introduce 
a variant of the rational Burnside problem,
 where we need to use 
pairs of parametrized rational functions.

Let $(f_t,g_t)$ be a pair of relatively elementary functions,
and consider the state system of the rational dynamics with initial sets:
\begin{align*}
&   z_i^{j+1}= f_t(w_i^j, z_i^j, \dots, z^j_{i+\alpha}), \quad
 w^j_{i+1}= g_t(w_i^j   ,z_i^j, \dots, z_{i+\beta}^j), \\
&  X_t \subset \{(z_0,z_1, \dots) : z_i \in (0, \infty)\}, \quad
 Y_t \subset \{(w^0,w^1, \dots) : w^i \in (0, \infty)\}
\end{align*}
\begin{df}
The state system of the parametrized rational dynamics by $(f_t,g_t)$
is quasi-recursive with respect to $(X_t,Y_t)$,
if  for any  $C, C'  \geq 1$, there exists $t_0 >1$  
so that for all $t \geq t_0$
and  any $\{w^j_0\}_j \in Y_t$, there exist some $p\in {\mathbb N}$
such  that:

(1)
any solutions $(\{z_i^j\}_{i,j}, \{w_i^j\}_{i,j})$ with 
 $\{z_i^0\}_i \in X_t$ 
 satisfy the uniform bounds:
$$ (\frac{z_i^{j+pl}}{z_i^j} )^{\pm 1} \leq C$$
  for all  $i,j, l =0,1,2, \dots  $, and 

(2) for any $1 \leq p' \leq p-1 $,
there are some  $\{z_i^0\}_i \in X_t$ so that the solutions 
 $(\{z_i^j\}_{i,j}, \{w_i^j\}_{i,j})$ satisfy the uniform lower bounds:
$$ (\frac{z_i^{j+p'}}{z_i^j} )^{\pm 1} \geq C' $$
 for all  $ j=0,1,2,  \dots$
 and  some  $ i$.

It is  infinitely quasi-recursive, if  infinitely many such $p$ exist.
 \end{df}
 Notice that lower bounds imply minimality of quasi-periods.

\vspace{3mm}
 
 Let us introduce {\em quasi periodic exponential lattices}:
 $$N_C^{per}(t^Q)=\{(w^0,w^1, \dots) : C^{-1}t^{q^j} \leq w^j \leq Ct^{q^j} , 
 \bar{q} \in X_Q \text{ are periodic} \}.$$

Let $S= \{s_0,s_1\} \subset {\mathbb Z}$ and 
$R = \{q^0,q^1\} \subset {\mathbb Z}$  be any embeddings.
Now we have the existence of parametrized infinitely quasi-recursive dynamics:
\begin{thm}
There exists a  pair of relatively elementary functions
$(f_t,g_t)$ so that 
  the state dynamics is 
 infinitely quasi-recursive.
 \end{thm}
 {\em Proof:}
Let us choose any embeddings $S= \{s_0,s_1\} \subset {\mathbb Z}$
and $Q = \{q^0, \dots, q^7 \} \subset {\mathbb Z}$.
Then we put:
  $$(X_t,Y_t) = (N_C(t^S), N_C^{per}(t^R)).$$
Notice that if two points $z,z' \in N_C(t^s)$, then their ratios satisfy the estimates
$(\frac{z}{z'})^{\pm 1} \leq C^2$.

 If  $z \in N_C(t^{s_0})$ and $z' \in N_C(t^{s_1})$,
then the estimates hold:
$$(\frac{z}{z'})^{\pm 1} \geq C^{-2}t.$$
In particular if $t >>1$ is sufficiently large,
then $C^{-2}t \geq C'$ holds.

Let us  consider 
 the automaton ${\bf A}$  by Aleshin in lemma $4.11$.
By proposition $4.3$,  there exists a stable extension of  ${\bf A}$.
Let $(f_t,g_t)$ be the pair of the corresponding relatively elementary functions.

Because the group generated by $\{q^0,q^1\}$ is infinite torsion, 
the conclusion follows from 
 proposition $4.10$.
This completes the proof.
\vspace{3mm} \\
{\bf 4.E.4 Estimates for dynamical inequalities:}
Let 
${\bf A}: \psi: Q \times S^{\alpha+1} \to S$, $  \phi: Q \times S^{\beta+1} \to Q$
be an automaton, and  choose
embeddings $Q,S \subset {\mathbb Z}$.

Let us take a
stable extension with the constants $(\delta, \mu)$:
\begin{align*}
& \psi: {\mathbb R} \times {\mathbb R}^{\alpha+1} \to {\mathbb R}, \\
& \phi:  {\mathbb R} \times {\mathbb R}^{\beta+1} \to {\mathbb R} .
\end{align*}
Let us represent the pair of  functions by two  relatively $(\max,+)$-functions
$(\psi^l,\phi^l)$  
for $l=1,2$,  and let 
 $(f_t^l,g_t^l)$  be the tropical  correspondences  respectively.
Thus
$(f_t^1, g_t^1) \sim (f_t^2, g_t^2)$
 are pairwisely  tropically equivalent.

\begin{cor} 
Suppose 
${\bf A}_{\bar{q}^m}: X_S \to X_S$ is of finite order with period $p $
for some $\bar{q}^m \in X_Q^{m+1}$.

Then 
for any $C  \geq 1$, there exists $t_0>1$ and $D\geq 1$   so that
  for all $t \geq t_0$ and 
any   pair of sequences $\{(z_i^j, w_i^j)\}_{i,j}$
which  satisfy the dynamical  inequalities:
\begin{align*}
& f^1_t(w_i^j, z_i^j, \dots, z_{i+\alpha}^j) \leq z_i^{j+1} \leq  
f^2_t(w_i^j, z_i^j, \dots, z_{i+\alpha}^j) , \\
& g^1_t(w_i^j,   z_i^j, \dots,z_{i+\beta}^j) \leq w_{i+1}^j \leq
 g^2_t(w_i^j,   z_i^j, \dots , z_{i+\beta}^j)
\end{align*}
with
the  initial values
  $\{z_i\}_i \subset N_C(t^S), \quad \{w^j\}_j \subset N_C(t^{ \bar{q}^m_{per}})$,
then they satisfy  the uniform bounds  for all $i,j,l =0,1,2, \dots$:
$$(   \frac{z_i^j}{z_i^{j+p(m+1)l}} )^{\pm 1} \ \leq    D.$$
\end{cor}
{\em Proof:}
Let us consider the state dynamics:
\begin{align*}
& F_i^{j+1}=  f^1_t(G_i^j, F_i^j, \dots, F_{i+\alpha}^j), \\
& G_{i+1}^j=  g^1_t(G_i^j, F_i^j, \dots, F_{i+\beta}^j)
\end{align*}
 with the initial values $F_i^0=z_i$ and $G^j_0=w^j$.

 By corollary $4.8$, the uniform estimates hold:
 $$(\frac{z_i^j}{F_i^j})^{\pm 1}  \leq D.$$

 By proposition $4.10$, the uniform bounds:
  $$(\frac{F_i^j}{F_i^{j+p(m+1)l}} )^{\pm 1} \leq C^4$$
 hold for all $i,j,l \geq 0$.
 
 Combining with these,  we obtain the desired estimates:
 $$(\frac{z_i^j}{z_i^{j+p(m+1)}} )^{\pm 1}=
  (\frac{z_i^j}{F_i^j })^{\pm 1}
  (\frac{F_i^j}{F_i^{j+p(m+1)}})^{\pm 1} 
( \frac{F_i^{j+p(m+1)}}{z_i^{j+p(m+1)}})^{\pm 1} \leq C^4D^2.$$
This completes the proof.
 \vspace{3mm}

\section{Analysis of hyperbolic systems of PDEs}
In section $5$, we develop basic analysis of the hyperbolic Mealy systems of PDE
in our sense, and then  apply the previous results to the large scale analysis of them.
\vspace{3mm} \\
 {\bf 5.A PDE systems and equivalent automata:}
 For $l=1,2$, let
$${\bf A}_l: \psi^l: Q \times S^{\alpha+1} \to S, \quad  \phi^l: Q \times S^{\beta+1} \to Q$$
be equivalent  automata over $R \subset Q$, 
and  choose embeddings $Q,S \subset {\mathbb R}$.

Let us take their stable extensions with the constants $(\delta, \mu)$:
$$\psi^l: {\mathbb R} \times {\mathbb R}^{\alpha+1} \to {\mathbb R}, \quad
\phi^l:  {\mathbb R} \times {\mathbb R}^{\beta+1} \to {\mathbb R} \qquad (l=1,2).$$  
and
$(f_t^l,g_t^l)$ and $(\psi^l_t, \phi^l_t)$ be the 
 corresponding functions  respecively.
 
 Let us
 consider the state systems of the rational dynamics:
\begin{align*}
&  z_i^{j+1}(l)= f^l_t(w_i^j(l), z_i^j(l), \dots, z_{i+\alpha}^j(l)), \\
& w_{i+1}^j(l)= g^l_t(w_i^j(l),  z_i^j(l), \dots, z_{i+\beta}^j(l)).
\end{align*}

Let us follow the process in $3.B$,
 and induce the systems of  partial differential equations of order $\mu$:
\begin{align*}
& P_1^l(\epsilon, t,  u^l,v^l, u_x^l,u_s^l,   \dots, u_{\mu s}^l)=0 \\
& P_2^l(\epsilon, t,  u^l,v^l, u_x^l,v_x^l,   \dots, u_{\mu x}^l)=0 
\end{align*}
  with  the scaling parameters:
   $$ i = \frac{x}{\epsilon}, \quad j= \frac{s}{\epsilon}, \quad
 u^l(x,s)= z_i^j(l) , \quad v^l(x,s)=w_i^j(l).$$

\vspace{3mm} 

Let us fix any positive number $0< \tau <0.5$, and  put the domains:
\begin{align*}
& D_1(\epsilon, \tau)=
\{ \  i\epsilon + a \in  [0, \infty)
: \  i  \in {\mathbb N} , \  \ |a| \  \leq  \epsilon \tau  \ \}, \\
& D(\epsilon, \tau)= D_1(\epsilon, \tau ) \times D_1( \epsilon, \tau ) \subset  [0, \infty)\times  [0, \infty).
\end{align*}
$D(\epsilon)$  are the disjoint unions of the squares.

We also put the initial domains:
$$I_x(\epsilon ,\tau) = [0, \infty) \times [0, \tau \epsilon) , 
 \quad   I_s(\epsilon ,\tau) = [0, \tau \epsilon) \times[0, \infty)  .$$
\vspace{3mm}

\begin{prop}
Let $C_0$ be the bigger one of their error constants.
Let  $u^l,v^l:  (0, \infty) \times [0, \infty)\to (0, \infty)$ be solutions
to the above systems respectively.
Then for any $C \geq 1$, there exists $t_0 >1$ and $D\geq 1$ so that 
the followings hold for all  $t \geq t_0$:

Suppose two conditions:
(1) The
estimates:
$$0 \leq CK(u^l,v^l) \leq (1- \mu) \epsilon^{-1}$$
are satisfied for some positive $\mu >0$.

(2) The inclusions hold:
$$u^l| I_x(\epsilon ,\tau) \subset N_C(t^S),  \quad
 v^l |   I_s(\epsilon ,\tau) \subset N_C(t^Q).$$
Then they 
satisfy  the uniform bounds:
 $$ (\frac{u^1}{u^2})^{\pm 1}(x,s)  \leq  D$$
 for all $(x,y) \in D(\epsilon, \tau)$.
\end{prop}
{\em Proof:}
We follow a similar argument as the proof of  theorem $3.2$.
Let $N_0 \geq \max(\mu^{-1}, 2-\mu)$ be an integer.

Let us  fix $0 \leq a, b \leq  \tau $, and
 put the domain lattices   by:
$$ L_{\epsilon, \tau}(a,b)=\{( \ (l_1+a)\epsilon, (l_2+b)\epsilon \ ) \in  [0, \infty) \times [0, \infty)
: l_1,l_2 \in {\mathbb N} \}.$$

Let us put:
$$z_i^j(l) =  u^l((i+a)\epsilon, (j+b)\epsilon), \quad
w_i^j(l)=  v^l((i+a)\epsilon, (j+b)\epsilon)$$

Then by the same way as the proof of theorem $3.2$, we obtain the 
inequalities:
\begin{align*}
& \frac{1}{N_0}  f^l_t(w_i^j(l), z_i^j(l), \dots, z_{i+\alpha}^j(l)) 
 \leq z_i^{j+1}(l) \leq N_0 
 f^l_t(w_i^j(l), z_i^j(l), \dots, z_{i+\alpha}^j(l)), \\
& \frac{1}{N_0}  g^l_t(w_i^j(l),  z_i^j(l), \dots, z_{i+\beta}^j(l)) \leq 
 w_{i+1}^j(l) \leq N_0 g^l_t(w_i^j(l),  z_i^j(l), \dots, z_{i+\beta}^j(l)).  
  \end{align*}
Notice that  $\frac{1}{N_0}f$ and $N_0f$ are both tropically equivalent to $f$.
 
By the assumption,  the initial values satisfy:
$$z_i^0(l) =u^l((i+a)\epsilon, b \epsilon) \in N_C(t^S), \ 
w_0^j(l) =v^l(a\epsilon, (j+b) \epsilon)  \in N_C(t^R).$$
It follows from  corollary $4.8$ that the uniform estimates:
$$ \max( \frac{z_i^j(1)}{z_i^j(2)},  \  \frac{z_i^j(2)}{z_i^j(1)}   \leq D.$$

Since the right hand side does not depend on $a,b$,
one obains the uniform bounds:
$$(\frac{u^1}{u^2})^{\pm 1}(x,s) \leq D$$
for all $(x,y) \in D(\epsilon, \tau)$. 
This completes the proof.
\vspace{3mm} \\
{\bf 5.B. Hyperbolic Mealy systems:}
In this section we study the basic analysis of 
$1$st order hyperbolic systems of PDE with $2$ variables.

Let $(f_t,g_t)$ be a pair of relatively elementary functions, 
corresponding to the pair of relative $(\max,+)$-functions $(\psi, \phi)$.

Let us  consider the corresponding hyperbolic Mealy systems:
\begin{align*}
&  \epsilon \ u_s = f_t(v,u) -u, \\
 &  \epsilon \ v_x  =g_t(v,u) -v.
 \end{align*}

  Recall the  higher distorsion for the Mealy systems:
$$ 
K(u,v) \equiv  \sup_{(x,s) \in [0, \infty)^2 }   
\max[\frac{||(u,v)||_{ \mu, 0}^1}{  u(x,s+ \epsilon)}, \ \
\frac{||(u,v)||_{ \mu, 0}^2}{  v(x+ \epsilon,s)}].$$

\begin{df}
The pair $(f_t, g_t)$ is  admissible, if 
there are $0< \delta < L <L'$, $ 0< \mu$ and constants $A$
so that there are  solutions 
with the initial values:
\begin{align*}
& u,v: [0,  \infty) \times [0,  \infty) \to [t^L, t^{L'}] ), \\
& \\
& ||u_x||C^0([0, \infty)  \times \{0 \}), \quad ||v_s||C^0( \{0\}  \times [0, \infty))  \leq A
\end{align*}
which satisfy the following estimates hold for all $0 \leq \alpha \leq 1$: 
 \begin{align*}
& [ |  ( f_t(v,u) -u) ( (f_t)_u(v,u)-1)|  + | v_s(f_t)_v(u,v)|]  (x,s)  
 < (2-  \mu) u(x,s+\alpha),  \\
& \\
& [ |  ( g_t(v,u) -v) ( (g_t)_u(v,u)-1)| +  |u_x(g_t)_v(u,v)|](x,s )   < (2-  \mu) v(x+\alpha,s).
 \end{align*}
\end{df}

\begin{cor}  For $l=1,2$, let $(f_t^l,g_t^l)$ be  admissible pairs, and 
  $(u^l,v^l): [0, \infty) \times [0, \infty) \to (0, \infty)$
be the solutions to  the above systems 
respectively.

 Then
 they satisfy the asymptotic estimates:
 \begin{align*}
 (\frac{u^1}{u^2})^{\pm 1}(x,s), & \  \ (\frac{v^1}{v^2})^{\pm 1}(x,s) \\
& \leq (M_0)^{6P_{\epsilon^{-1}( x+s (\gamma+1))}(c)}  \ \
 ([ (u^1,v^1):(u^2,v^2)]_{\epsilon})^{\tilde{c}^{\epsilon^{-1}(x+s(\gamma+1))+1}}  
 \end{align*}
for some  $M_0$ and  all $(x,s) \in [0,  \infty) \times [0, \infty)$.
 \end{cor}
{ \em Proof:}
For the Mealy systems, the error constants are always $ \frac{1}{2}$.
Moreover we have the equalities:
\begin{align*}
& \frac{\partial^2 u}{\partial_s^2} =  ( f_t(v,u) -u) ( (f_t)_u(v,u)-1)+  (f_t)_v(u,v) v_s, 
  \\
& \frac{ \partial^2 v}{\partial_x^2} =
  ( g_t(v,u) -v) ( (g_t)_v(v,u)-1)+ (g_t)_u(u,v) u_x.
\end{align*}

So in order to apply theorem $3.1$, the required consitions on the higher distorsion are: 
 \begin{align*}
& \frac{|  ( f_t(v,u) -u) ( (f_t)_u(v,u)-1)|  + | v_s(f_t)_v(u,v)| (x,s)  }{2u(x,s+\alpha)}
 < (1-  \delta ) \epsilon^{-1},  \\
& \frac{ |  ( g_t(v,u) -v) ( (g_t)_v(v,u)-1)| +  |u_x(g_t)_u(u,v)|(x,s )}{2v(x+\alpha,s) }   
< (1-  \delta) \epsilon^{-1}
\end{align*}
for some $0 < \delta <1$ and $0< \epsilon \leq 1$, which 
 follow from admissibilty.

This completes the proof.

\vspace{3mm}

In  $5.D$ we have concrete examples of admissible pairs
with the Lipschitz constants $1$. In particular we obtain the exponential 
estimates for their solutions.
\vspace{3mm} \\
{\bf 5.B.3 Refinement:}
Let: $${\bf A}: \psi: Q \times S^{\alpha+1} \to S, \quad  \phi: Q \times S^{\beta+1} \to Q$$
be an automaton over $R \subset Q$.
For $\bar{q}=(q^0, q^1, \dots) \in X_Q$
and $\bar{k}=(k_0,k_1, \dots) \in X_S$,
let us denote the orbits by $ \{k_i^j \}$ and $ \{q_i^j \}$.

Let
$\tilde{\psi}: {\mathbb R} \times {\mathbb R}^{\alpha+1} \to {\mathbb R}$ and 
$\tilde{\phi}:  {\mathbb R} \times {\mathbb R}^{\beta+1} \to {\mathbb R}$
be two maps.
Suppose the restrictions induce maps as below for some $a,b\in {\mathbb R}\cup \{ \pm \infty\}$:
 $$(\tilde{\psi}, \tilde{\phi}): [a,b]^2 \to  [a,b]^2 .$$ 
For $\bar{y}=(y^0, y^1, \dots) ,
\bar{x}=(x_0,x_1, \dots) \in X_{[a,b]}$, 
let us denote the orbits by $ \{x_i^j \}$ and $ \{y_i^j \}$.
Let us choose embeddings  $S,Q \subset [a,b]$.

\begin{df} $(\tilde{\psi}, \tilde{\phi})$ is an
$\epsilon$-refinement of the pair $(\psi, \phi)$, if
 there is   a posiive number $N$ 
so that for any $\bar{q}\in X_Q$
and $\bar{k}\in X_S$,
there are paths $y: \{0, 1, \dots\} \to {\mathbb R}$
and $x: \{0,1, \dots \} \to {\mathbb R}$
with:
\begin{align*}
& y(jN)= q^j , \quad |y(j+1)-y(j)|\leq \epsilon, \\
& x(iN) =k_i, \quad |x(i+1) - x(i)| \leq \epsilon
\end{align*}
 for all $i , j \in \{ 0,1, \dots \}$, such that the equalities hold:
$$x_{iN}^{jN}= k_i^j, \quad   y_{iN}^{jN} = q_i^j.$$
$(\psi, \phi)$ is refinable, if there is an $\epsilon$-refinement for any small $\epsilon >0$.

An $\epsilon$ refinement is  almost diagonal, if moreover
 they
 satisfy the estimates for all $(x,y) \in [a,b]$:
$$|(x,y) - (\tilde{\psi}(x,y) ,  \tilde{\phi}(x,y))| \leq \epsilon.$$

\end{df}
By use of refinement, we verify existence of admissible solutions in $5.D$,
and their exponential asymptotic estimates, which we will state below.
 \vspace{3mm}  \\
{\bf 5.B.4 Asymptotic comparisons with group actions:}
Let ${\bf A}$ be a Mealy automaton with $2$ alphabets, 
equipped with a representative $(\psi, \phi)$  by relatively $(\max,+)$-functions. 

For any $(s_0,s_1,\dots) \in X_2$, and $(q^0,q^1, \dots) \in X_{m+1}$,
let:
$${\bf A}_{(q^0, \dots,q^{l-1})} (s_0,s_1, \dots) = ((s_0^l, s_1^l, \dots ) \in X_2$$
be the orbits of the automata group actions.

For its refinement $(\bar{\psi}, \bar{\phi})$ and 
their tropical correspondences $( \bar{f}_t, \bar{g}_t)$, let us consider the hyperbolic Mealy systems:
$$u_s=  \bar{f}_t(v,u)-u, \quad v_x=  \bar{g}_t(v,u)-v.$$

\begin{thm}   For any $C>0$ and any $t \geq t(C) >1$, 
there are refinements 
 $(\bar{\psi}, \bar{\phi})$ of $( \psi,  \phi)$ with the   pairs of
 tropical correspondences $(\bar{f}_t,\bar{g}_t)$
so that:

(1) $(\bar{f}_t,\bar{g}_t)$ admits admissible solutions, 

(2)
for any another pairs $(f_t,g_t)$ toropically equivalent to 
$(\bar{f}_t,\bar{g}_t)$, 
 any admissible solutions to the equations:
$$u_s=f_t(v,u)-u, \quad v_x= g_t(v,u)-v$$
whose  initial values satisfy the inclusions for all $k =0,1,2, \dots$:
$$d(u(N k, 0) , S) ,\quad d(v(0,Nk ), Q) \leq C$$
then  they satisfy the asymptotic estimates  for some $M, c \geq 1$:
 \begin{align*}
 (\frac{u(N i, N j) }{s_i^j})^{\pm 1} & \  \ \leq M^{ P_{N(i+j) +1 } (c)}
 \end{align*}
 \end{thm}
{ \em Proof:}
This follows from proposition $5.12$ below.
 \vspace{3mm} \\
{\em Remark:}
(1) We can choose $c=1$ and hence obtain the exponential estimates, if
we can represent the transition functions by $1$-Lipschitz functions.

(2)
It would be quite possible to obtain better estimates, by 
tracing the orbits in detail, but one will be required more analysis to achieve it.
\vspace{3mm} \\
{\bf  On $\epsilon$-controll of higher distorsions:}
In order to apply the construction in section $3$,
we need to controll the second derivatives of solutions.

Let $(u, v): [ 0, \infty) \times [0, \infty) \to (0, \infty)^2$ be  a solution to 
 the hyperbolic Mealy system:
\begin{align*}
&  \epsilon \ u_s = f_t(v,u) -u, \\
 &  \epsilon \ v_x  =g_t(v,u) -v
 \end{align*}
with finite  higher distorsion $K(u,v) < \infty$.

In order to apply the asymptotic comparisons for solutions to the PDE systems, 
their solutions are required to satisfy some bounds:
$$CK(u,v) \leq (1- \delta)\epsilon^{-1}$$
for some $0 < \delta <1$, where $C$ is the error constant.

For small $0< \tau $,  let us consider
the reparametrized pair of the functions:
 $$(\tilde{u}, \tilde{v})(x,s) \equiv (u,v)(\tau x,\tau s)$$
which satisfy the equations:
\begin{align*}
&  \epsilon \tau^{-1} \ \tilde{u}_s = f_t(\tilde{v},\tilde{u}) -\tilde{u}, \\
 &  \epsilon \tau^{-1}  \ \tilde{v}_x  =g_t(\tilde{v},\tilde{u}) -\tilde{v}.
 \end{align*}

 \begin{lem}
 The pairs $(\tilde{u}, \tilde{v})$ satisfy the condition above on the higher distorsion
 if we choose sufficiently small $\tau >0$.
 \end{lem}
 {\em Proof:}
Notice that  the higher distortions are related as:
$$\tau^2 K_2(u,v)= K_2(\tilde{u},\tilde{v}).$$
The required condition with respect to the  reparametrized pair is given by:
$$ CK_2(\tilde{u},\tilde{v})
\leq (1- \delta) \tau \epsilon^{-1}   $$
which can be rewritten as:
 $$ \tau CK_2(u,v) \leq (1- \delta)  \epsilon   .$$
 So if we choose sufficiently small $\tau >0$, then 
 the pairs $(\tilde{u}, \tilde{v})$ satisfy the required condition.
 This completes the proof.

 \vspace{3mm}

On the other hand the estimates $0< \epsilon  \tau^{-1} \leq 1$
are required, and so at best we can take $ \tau = \epsilon$.
So for our later purpose 
 this reparametrization is not so effective.
We will  take:
$$ \epsilon =1$$
for the rest of this paper.
\vspace{3mm} \\
{\bf 5.C Basic PDE-analysis,  existence and uniqueness:}
 Let us start from some analytic properties of the hyperbolic Mealy systems.

 \begin{lem}
 (1) Suppose $\psi$ and $\phi$ both take bounded values  from both below and above.
 Then  $f_t$ and $g_t$ also satisfy the same properties for each $t >1$.

 \vspace{2mm}

Let us fix $t >1$, and
$(u,v) : [0, \infty) \times [0, \infty) \to {\mathbb R}^2$ be a solution
to the hyperbolic Mealy system,
 which take positive  initial values  on  $x , s \in [0, \infty)$:
$$u(x,0), v(0,s) \in (0, \infty).$$

(2) 
 Suppose  that 
there are constants $r < R$ so that 
 $f_t$ and $g_t$ both satisfy the uniform bounds:
$$r  \leq f_t(a,b), \ g_t(b,a)  \leq R.$$
 Then they take positive an bounded values on all the domain:
 $$u,v : [0, \infty) \times [0, \infty) \to (0, \infty).$$

(3)
If $f_t$ and $g_t$ both satisfy the bounds:
$$r  \leq f_t(a,b), \ g_t(b,a)  \leq R$$
for any $r\leq b\leq R$ and all $a>0$,
and  if we choose  the initial values with their ranges as:
$$r  \leq u(x,0), v(0,s)  \leq R$$
then their values are also contained in the range:
$$r  \leq u(x,s), v(x,s)    \leq R$$
for all $(x,s)  \in [0, \infty) \times [0, \infty) $.
\end{lem}
{\em Proof:}
(1) follows from lemma $2.1$ and proposition $2.4$.

Let us fix $x \in [0, \infty)$ and $v(x,s)$. Then we regard 
$u(x, \quad): [0,  \infty) \to { \mathbb R}$ satisfy 
 the ODE on $s$ variable:
 $$ \epsilon u_s (x,  \quad)=f_t(v(x,s), u(x, \quad))-u(x,  \quad).$$

For (2) and (3),
 because the range of $f_t$ is away from $0$, if $u(x, \quad)$ 
take small values, then $u_s(x,  \quad)$ become positive and so their values must increase.
Conversely  if $u(x, \quad)$ 
take large values, then $u_s(x,  \quad)$ become negative and  their values must decrease.

For $v(\quad , s)$ case, we can use the same argument.
From these observations, the conclusions follow immediately.
This completes the proof.

\vspace{3mm}

Let us say that the pair $(f_t,g_t)$ {\em restricts to a self-dynamics} over $[r,R]$ (for a fixed $t>1$), 
if there is  some $0< q <r$ so that they  satisfy the bounds
for  all $a >0$:
$$  f_t(a,b) -b , \  \ g_t(b,a) -b \ \
\begin{cases} 
 \geq  q & \ \ b \leq r+q  \\
\leq -q   &  \ \ b \geq R-q
\end{cases} $$
$$| f_t(a,b) -b |, \  \ | g_t(b,a) -b|  \leq R-r  \qquad r  \leq b  \leq R.$$
If $(f_t, g_t)$ satisfies the uniform bounds as $r \leq f_t , g_t \leq R$,
then the pair restricts to a self-dynamics over $[r-2q ,R+2q]$.

By lemma $5.5$,  if we take  the initial values as
$r+q  \leq u(x,0), v(0,s)  \leq R-q$,
then their values are also contained in the range:
$$r +q \leq u(x,s), v(x,s)    \leq R-q$$
for all $(x,s)  \in [0, \infty) \times [0, \infty) $.
\vspace{3mm}  \\
{\bf 5.C.2 Existence and uniqueness:}
Let us  study the existence of solutions to the hyperbolic Mealy systems.
We use the automatic  version of the
 Picard iteration method of successive approximation.

Let $(f_t,g_t)$ be a pair of  relatively elementary functions of $2$ variables,
which correspond to the relative $(\max, +)$-functions $(\psi,\phi) $.
Let us assume that the pair restricts to a self-dynamics over $[r,R]$ with $q<r$.

Let us introduce the following:
$$  D=  D_{r,R}(f_t,g_t) = \sup_{(v,u) \in  [r, R]^2 } \{|f_t(v,u)-u|, \  \ |g_t(v,u)-v|\}.$$

\vspace{3mm}

Let us put 
$\bar{f}_t (v,u) = f_t(v,u)-u$ and $\bar{g}_t (v,u) = g_t(v,u)-v$ so that 
the hyperbolic Mealy equations can be written as:
$$u_s =\bar{f}_t (v,u) , \quad v_x = \bar{g}_t (v,u).$$

\vspace{3mm}

Now let us fix $t>1$ and give the initial values:
$$u \ | \ [0, \infty) \times \{0\}, \qquad  v \ | \ \{0\} \times [0, \infty)$$
which satisfy the followings:

\vspace{3mm} 

(1) Their ranges are uniformly bounded both from above and below by:
$$r+q \leq u(x,0), v(0,s)  \leq R-q.$$

\vspace{3mm}

(2) Both have  uniformly bounded $C^2$ norms.
\vspace{3mm} \\
With $u$ amd $v$ as the initial values,
  let us solve the equations of the hyperbolic Mealy systems
by the following constructions (1),(2),(3), (4):
 \vspace{3mm} \\
{\bf (1):} Let us choose small $\tau >0$ so that:

\vspace{2mm}

(a)  the Lipschitz constants of both $\bar{f}_t$ and $\bar{g}_t$ 
are smaller than $(2\tau)^{-1}$,

\vspace{2mm}

(b) the estimates holds:
 $$\tau \leq  D^{-1} q.$$ 
{\bf (2):}
For  $m =0,1,2, \dots$, let us  put the small domains: 
$$D_m = [m \tau , (m+1)\tau] \times [0, \tau] \subset [0, \infty) \times [0,\infty)$$
and construct solutions inductively on $m=0,1, \dots$.

Suppose given solutions over $D_{m-1}$. With the original initial values,
  we have determined the values of solutions as:
$$u| D_{m-1} \cup [0, \infty) \times \{0\}  ,
 \quad v| | D_{m-1} \cup \{0\}\times  [0, \infty) .$$
{\bf (3):}
For $x_0 = m \tau$, we  have the initial values as:
 $$ u|  [x_0, x_0+\tau] \times \{0\} , \quad
v| \{x_0\} \times [0,\tau]$$
 and let us extend the solutions  over $D_m$.

For $(x,s) \in D_m$, let us put:
$$(u_0,v_0) = (u(x,0), v(x_0,s))$$
and define the sequences inductively by:
$$(u_n,v_n) = (u_0,v_0)+ ( \int_0^s \bar{f}_t(v_{n-1},u_{n-1})dt , \int_{x_0}^x \bar{g}_t(v_{n-1},u_{n-1})dy).$$

\begin{lem}
Suppose that  the pair $(f_t, g_t)$ restricts to a self-dynamics over $[r,R]$.
Then the  sequences $\{(u_n,v_n)\}_n$ converge uniformly on $D_m$, and:
$$(u,v)= \lim_{n \to \infty} (u_n,v_n)$$
give the solutions  which coincide with  the given initial values:
$$ u| [x_0, x_0+\tau] \times \{0\} ,\quad v| \{x_0\} \times [0,\tau].$$
Moreover they satisfy  the estimates:
$$|(u,v)-  (u_n,v_n)|  \leq \frac{1}{2^n}\max(   |u_1-u_0|, |v_1 - v_0| ).$$
\end{lem}
{\em Proof:}
{\bf Step 1:} We claim that the ranges of $u_n$ and $v_n$ are in the region $[r, R]$.
We verify it only for $u_n$. $v_n$ can be considered similarly.

By the assumption, the uniform bounds $r +q\leq u_0  \leq R-q$ hold.
Firstly we have the estimates:
$$|u_1 -u_0| \leq \int_0^{\tau} | f_t(v_0,u_0)-u_0| \leq \tau D \leq
 D^{-1} qD = q.$$
So $u_1$ admits the bounds $r \leq u_1 \leq R$.

Let us consider:
$$u_n = u_0 +  \int_0^s [f_t(v_{n-1},u_{n-1}) -u_{n-1}] dt$$
and assume the uniform  bounds
$r \leq u_{n-1} \leq R$. Then 
we have the estimates:
$$|u_n -u_0| \leq \int_0^{\tau} | f_t(v_{n-1},u_{n-1})-u_{n-1}| \leq \tau D \leq
 D^{-1} qD= q.$$
So we have verified the uniform bounds
$r \leq u_n \leq R$ for all $n$  by the induction step.
Since the estimates are independent of chioce of $x \in   [x_0,x_0+  \tau] $, this 
verifies the claim.

\vspace{3mm}

{ \bf Step 2:}
By step $1$, both $u_n$ and $v_n$ take their values in $[r,R]$.
So the Lipschitz constants 
$L_{\bar{f}}$ of $\bar{f}_t$ and $L_{\bar{g}}$ are both uniformly bounded at $(u_n,v_n)$.

Let us put $V_n =(u_n,v_n)$. Then:
\begin{align*}
& |V_{n+1}-V_n| \\
& = |( \int_0^s( \bar{f}_t(v_n,u_n) -  \bar{f}_t(v_{n-1},u_{n-1}))dt , 
\int_{x_0}^x(\bar{g}_t(v_n,u_n) -  \bar{g}_t(v_{n-1},u_{n-1})dy)| \\
& \leq \max( s L_{\bar{f}} |V_n - V_{n-1}|, (x-x_0)L_{\bar{g}} |V_n - V_{n-1}|) \\
& \leq \tau \max(L_{\bar{f}}, L_{\bar{g}})  |V_n - V_{n-1}|  <  \frac{1}{2} |V_n - V_{n-1}| .
\end{align*}
Thus 
they are contracting:
$$\max \{|u_n-u_{n-1}|,  \ |v_n-v_{n-1}|\}
\leq \frac{1}{2^n} \max\{|u_1- u_0|, \ |v_1 - v_0| \}$$
 and  $(u,v) =\lim_{n \to \infty} (u_n,v_n)$
exist uniformly.
One also obtains the estimates:
$$\max \{|u_n-u|,  \ |v_n-v| \}
\leq \frac{1}{2^n} \max\{|u_1- u_0|, \ |v_1 - v_0| \}.$$

Then they satisfy the integral equations:
$$(u,v)= (u_0,v_0)+ (\int_0^s \bar{f}_t(v,u)dt, \int_{x_0}^x \bar{g}_t(v,u)dy)$$
which are equivalent to the hyperbolic Mealy equations.

Let us check that these satisfy the boundary conditions.
On  $ [x_0, x_0+\tau] \times \{0\} \cup \{x_0\} \times [0, \tau]$,
\begin{align*}
& (u,v)(x,0)= (u(x,0) \ ,\ v(x_0,0)+  \int_{x_0}^x \bar{g}_t(v,u)dy) =(u(x,0),v(x,0)) \\
&  (u,v)(x_0,s)= (u(x_0,0)  +  \int_0^s \bar{f}_t(v,u)dy \ , \ v(x_0,s))=(u(x_0,s),v(x_0,s)).
\end{align*}
So the pair $(u,v)$ certainly satisfies the boundary condition.

This completes the proof.
\vspace{3mm} \\
{\bf (4):}
Let us continue the construction of solutions on $[0, \infty) \times [0, \infty)$.
By lemma $5.4$, let us extend solutions on $D_m$ inductively, and
one obtains solutions on $[0, \infty) \times [0, \tau]$.

Since the initial values  $u| [0 , \infty) \times \{0\}$ and 
$v|\{0\} \times [0, \infty)$  take  their ranges between $[r- q, R+q]$,
the range of both $u,v$ on $[0, \infty) \times [0, \tau]$ also the same, 
since the pair $(f_t,g_t)$ restricts to a self-dynamics over $[r,R]$.

Next let us regard that:
 $$[0, \infty) \times \{ \tau \} \cup \{0\} \times [\tau, 2\tau]$$
is the boundary equipped  with the boundary condition. Let us 
 iterate the same construction
over $[0, \infty) \times [ \tau, 2\tau]$.
By the above observation, the range of the initial conditions are also contained in 
$[r+q, R-q]$. 
So one can repeat the above process by the same way.

By repeating this process, one finally obtains the solutions with the given boundary condition:
$$u,v: [0, \infty) \times [0, \infty)  \to  { \mathbb R}   .$$
This completes the construction of  solutions to the hyperbolic Mealy systems.

\vspace{3mm}

Let us put the initial domain:
$$I_0 = [0, \infty) \times \{0\} \  \cup   \ \{0\} \times [0, \infty) .$$

\begin{thm}
Suppose that  the pair $(f_t, g_t)$ restricts to a self-dynamics over $[r,R]$, 
and give the positive initial values:
$$u , \ v :   \ I_0  \to [r+q, R- q].$$
Then:

(1)  there exists a positive solution 
 $$u,v : [0, \infty) \times [0, \infty) \to (0, \infty)$$
 with the uniform bounds:
 $$r +q \leq u(x,s), \ v(x,s) \leq R-q$$
 hold.

 (2) The solution is unique.
 \end{thm}
{\em Proof:}
The first statements follow from the above construction with lemma $5.2$.

Let us verify uniqueness.
Suppose two solutions $(u,v)$ and $(u',v')$ exist with the same initial values.
Let us put:
 $$a= |u-u'|, \ b=|v-v'| : [0, \infty)^2 \to [0, \infty)$$
and  verify $a=b\equiv 0$.
Notice that 
 $a| [0,\infty) \times \{0\} $ and $b| \{0\}  \times [0,\infty)$
both vanish.
Firstly, at $x=0$, the ODE:
 $$u_s(0,s) = f_t(v(0,s), u(0,s))-u(0,s)$$ has the unique solution 
on $s \in [0, \infty)$, where  $v(0,s)$ is the given initial value.
So $u(x,0)=u'(x,0)$ hold for all $x \in [0, \infty)$.

Similarly $v(0,s) =v'(0,s)$ holds.
In particular:
$$u(x,s) =u'(x,s) , \quad v(x,s)=v'(x,s) $$
hold for all $ (x,s) \in [0,\infty) \times \{0\} \cup  \{0\}  \times [0,\infty)$.

The solutions satisfy the integral equations:
$$ u(x,s)= u(x,0)+   \int_0^s \bar{f}_t(v,u)ds,  \quad
v(x,s) = v(0,s) + \int_0^x \bar{g}_t(v,u)dx.$$
Since both $ \bar{f}_t$ and $  \bar{g}_t$ are Lipschitz, there exists $L$ so that the estimates hold:
$$a \leq L \int_0^s (a+b), \ \ b \leq L\int_0^x (a+b).$$

Let us denote the broken lines:
$$\gamma_{x,s} = \{ (\alpha,s )\cup (x, \beta): 0 \leq \alpha \leq x , \ 0 \leq \beta \leq s\}.$$ 
Then combination of the above inequalities gives the following inequalities:
$$(a+b)(x,s) \leq L \int_{\gamma_{x,s}} (a+b).$$
Let us choose $(x,s)$ so that
the estimate $L(x+s)<\delta < 1$ holds.

Now let $D \subset [0, \infty)^2$ be the rectangle
whose boundary is given  by:
  $$ \partial D=  \gamma_{x,s}  \cup
  [0,x] \times \{0\} \cup  \{0\}  \times [0,s].$$
Let us 
choose some point  $(x_0,s_0) \in D$ so that the equality
$(a+b)(x_0,s_0) = \sup_{(p,q) \in D} (a+b)(p,q)$ holds.
Then for $(x,s) \in D$, the estimates hold:
$$(a+b)(x,s) \leq (a+b)(x_0,s_0) \leq 
L \int_{\gamma_{x_0,s_0}} (a+b) \leq \delta  \ (a+b)(x_0,s_0).$$
Since $\delta<1$, this implies $(a+b)(x_0,s_0)=0$ and hence
$(a+b)(x,s) \equiv 0$ on $D$.

By changing the domains by  parallel transport as in the above construction,
we can follow the  same argument as above. 
By iterating the same process, 
we conclude that $a=b  \equiv 0$ hold on $ [0, \infty)^2 $.

This completes the proof.
\vspace{3mm} \\
{\bf 5.C.3 Energy estimates:}
Let us study $C^1$ estimates of solutions.
It is well known as the energy estimates for hyperbolic equations.
Here we   also give concrete estimates of   the constants which appear in 
 the asymptotic growth of $C^1$ derivatives of solutions.

For functions $u: [0,  \infty)^2 \to [0, \infty)$,
let us denote the norms of the first derivative by:
$$||u||\bar{C}^1 = \sup_{(x,s) \in [0, \infty)^2} 
\max \{ |\frac{\partial u}{\partial x}|(x,s),  |\frac{\partial u}{\partial s}|(x,s)  \}.$$

The following is elementary:
\begin{lem}
Let us give any  discrete initial values:
$$  u :  {\mathbb N} \times \{0\} \to  {\mathbb R} , \quad
  v :  \{0\}\times {\mathbb N}  \to  {\mathbb R} $$
such that the estimates holds for all $n =0,1,2,\dots$:
$$|u(n+1)-u(n)|, \   |v(n+1)-v(n)| \leq \mu.$$

Then there is a constant $C$ independent of choice of the initial values 
and $ \mu$  so that
 there are extensions of the initial values as:
 $$u : [0,  \infty)  \times  \{0 \} \to {\mathbb R}  ,  \quad  
v: \{0 \} \times [0,  \infty) \to {\mathbb R}$$
equipeed with uniform $ \bar{C}^1$  bounds:
 $$||u||\bar{C}^1([0,  \infty)  \times  \{0 \}) ,  \quad  ||v||\bar{C}^1(   \{0 \} \times [0,  \infty)) \leq C\mu  .$$
\end{lem}
{\em Proof:} 
Let us extend the initial values by connecting the discrete values by segments.
Then its $C^1$ approximations give the desired property.

This completes the proof.

\vspace{3mm}

Recall the number 
 $D=  D_{r,R}(f_t,g_t) $ in $5.C.2$ and 
 introduce another one:
$$ B = \max(  ||(f_t)_u -1||C^0 , ||(f_t)_v||C^0, \ ||(g_t)_v-1||C^0 , ||(g_t)_u||C^0) .$$

The next gives the exponential energy estimates:
\begin{prop}
Suppose that  the pair $(f_t, g_t)$ restricts to a self-dynamics over $[r,R]$, 
and give the  initial values:
$$u(\quad ,0) , \ v(0, \quad) : [0,  \infty)       \to [r+q, R- q]$$
with uniformly bounded $C^1$norms:
$$||u_x||C^0([0,  \infty)  \times  \{0 \}) ,  \quad  ||v_s||C^0(   \{0 \} \times [0,  \infty)) \leq A  <  \infty.$$

Then  there is a constant $C$ so that solutions 
 $u,v : [0, \infty) \times [0, \infty) \to (0, \infty)$
 have the  asymptotic $C^1$ bounds:
  \begin{align*}
  &||\frac{\partial u}{\partial x}||C^0( [0, \infty)\times \{m \}) , \  \
  ||\frac{\partial v}{\partial s}||C^1(  \{m\} \times [0, \infty) ) \  \leq \ 2^{\tau^{-1}m } (A  + 2D), \\
&   ||\frac{\partial u}{\partial s}||C^0 , \quad  ||\frac{\partial v}{\partial x}||C^0 \ \ \leq  \ \ D
\end{align*}
where:
$$ \tau \text{ Lip}_{\bar{f}_t,  \bar{g}_t} \leq \frac{1}{2}, 
\quad \tau \leq D^{-1}q , \quad \delta \equiv \tau B \leq \frac{1}{4}.$$
  \end{prop}
{\em Proof:}
Let us use the notations   in  $5.C.2$.
Firstly 
we split the domains into periodic stripes. Then
we verify some uniform estimates
on the first derivatives of the approximated solutions
$(u_n,v_n)$ over each stripe, which are  independent of $n$.
Secondly we verify that they converge to the solutions 
uniformly in $C^1$ over each stripe.

\vspace{3mm}

{\bf Step 1:}
By theorem $5.5$, the range of solutions $(u,v)$ are uniformly bounded
between $r$ to $R$.
Then by the defining equations, 
both  uniform bounds hold:
$$||u_s||C^0, ||v_x||C^0 \leq D.$$

Let us estimate $u_x$ as follows.
By theorem $5.7$, 
the solutions are unique with respect to the given initial values.

Let us recall the inductive construction of solutions over $[0, \infty) \times [0, \tau]$ in  $5.C.2$.
We split the domain into periodic squares $D_m = [m \tau , (m+1)\tau] \times [0, \tau]$,
and  construct solutions successively on each $D_m$ as $ \lim_{n  \to \infty} (u_n,v_n)$, where:
\begin{align*}
& u_n = u_0 + \int_0^s(  f_t(v_{n-1}, u_{n-1})-u_{n-1})dt, \\
& v_n = v_0+ \int_{x_0}^x (g_t(v_{n-1},u_{n-1}) - v_{n-1})dy
\end{align*}
for $0 \leq s  \leq \tau$ and $x_0 =m\tau \leq  x \leq (m+1)\tau$.

Let us denote $u_n'= \frac{\partial u_n}{\partial x}$. Then:
\begin{align*}
& u_n' = u_0' + \int_0^s( [ ( f_t)_u (v_{n-1}, u_{n-1})-1] u_{n-1}' + ( f_t)_v (v_{n-1}, u_{n-1})v_{n-1}' 
   )dt, \\
& v_n' = g_t(v_{n-1},u_{n-1}) - v_{n-1}.
\end{align*}

By lemma $5.6$, the estimates:
$$|v-v_n|, \ \ |u-u_n|  \leq \frac{1}{2^n}\max(  |u_1-u_0 |  , |v_1-v_0|) \leq \frac{q}{2^n} $$
hold. Since  the equality $ v_x = g_t(v,u) - v$  holds,
we obtain  the estimates:
 \begin{align*}
|v_n'|&   \leq |v'| + |v_n'-v'| \leq D + |(g_t(v,u)-v )- (g_t(v_{n-1},u_{n-1})-v_{n-1})|  \\
& \leq D+ \frac{q}{2^{n-1}}  Lip_{ \bar{g}_t}  \leq  D+ \frac{q}{2^n \tau}.
 \end{align*}

Let us put:
 $$ \tau B =  \delta \leq \frac{1}{4}.$$
Then  we have the estimates:
\begin{align*}
&  | u_n'|        \leq |u_0'| +  \tau B[  |u_{n-1}' |+  D+  \frac{q}{2^n \tau}] \\
&  = |u_0'| +(\delta D+ B  \frac{q}{2^n})  + \delta |u_{n-1}'|  
  \leq   \dots  \\
&  \leq(1+\delta + \delta^2 + \dots) ( |u_0'| + D)+
(\frac{1}{2^n} + \frac{1}{2^{n+1}}+ \dots) Bq  \\
& \leq  2( |u_0'|+  D+ \frac{Bq}{2^n} )
\leq 2(A+  D+ \frac{Bq}{2^n} ).
\end{align*}

Similarly we have the estimates
 $|\frac{\partial v_n}{\partial s}| \leq  2(A +D + \frac{Bq}{2^n})$
on $[0,   \tau] \times [0,  \infty)$.

By lemma $5.6$, we have the estimates:
$$|\frac{\partial u_n}{\partial s}|, \quad  |\frac{\partial v_n}{\partial x}| \ \leq D.$$

Notice that so far we have required the conditions on $\tau >0$ as:
$$\tau \text{ Lip}_{\bar{f}_t,  \bar{g}_t} \leq \frac{1}{2}, 
\quad \tau \leq D^{-1}q , \quad \delta = \tau B \leq \frac{1}{4}.$$ 

 \vspace{3mm}

{\bf Step 2:}
Let us verify $C^1$ uniform convergence of $\{u_n,v_n \}_n$.
For simplicity of the argument, we assume that 
 the Lipschitz constants of  $(\bar{f}_t)_u, (\bar{g}_t)_v$ and $(f_t)_v, (g_t)_u$
are all finite and bounded by $B'$.

Let us consider the estimates: 
$$ |v_n'|  \leq D, \quad
  | u_n'|        \leq 2A +2D + \frac{2Bq}{2^n} \equiv \alpha_n.$$
on the domains  $[0,   \infty) \times [0,  \tau]$.
Then we have the estimates:
\begin{align*} 
&
 | ( \bar{f}_t)_u (v_{n-1}, u_{n-1})u_{n-1}'-  (\bar{f}_t)_u (v_{n-2}, u_{n-2})u_{n-2}'|  \\
& \leq  | (\bar{f}_t)_u (v_{n-1}, u_{n-1})u_{n-1}'-  ( \bar{f}_t)_u (v_{n-1}, u_{n-1})u_{n-2}'| \\
& \qquad \qquad +  | ( \bar{f}_t)_u (v_{n-1}, u_{n-1})u_{n-2}'-  ( \bar{f}_t)_u (v_{n-2}, u_{n-2})u_{n-2}'|  \\
&\leq B|u_{n-1}' - u_{n-2}'| + B'(|u_{n-1}-u_{n-2}|+|v_{n-1}-v_{n-1}|) |u_{n-2}'| \\
&\leq B|u_{n-1}' - u_{n-2}'| + B' \alpha_{n-2} (|u_{n-1}-u_{n-2}|+|v_{n-1}-v_{n-2}|) ,
\end{align*}
\begin{align*}
& | ( f_t)_v (v_{n-1}, u_{n-1})v_{n-1}'-  ( f_t)_v (v_{n-2}, u_{n-2})v_{n-2}'|  \\
& \leq B|v_{n-1}' - v_{n-2}'| + B'(|u_{n-1}-u_{n-2}|+|v_{n-1}-v_{n-2}|) |v_{n-2}'| \\
& \leq B(  g_t (v_{n-2}, u_{n-2})-  g_t(v_{n-3}, u_{n-3}) |+ | v_{n-2} - v_{n-3}|)  \\
& \qquad  + B'D(|u_{n-1}-u_{n-2}|+|v_{n-1}-v_{n-2}|) \\
& \leq B\tau^{-1}(  |u_{n-2} -u_{n-3}|+  | v_{n-2} - v_{n-3}|)  \\
& \qquad  + B'D(|u_{n-1}-u_{n-2}|+|v_{n-1}-v_{n-2}|) .
\end{align*}

By lemma $5.6$, the estimates hold:
$$|u_n-u_{n-1}|, |v_n-v_{n-1}| \leq \frac{q}{2^n} .$$
Thus we have the estimates:
\begin{align*}
 & |u_n' -u_{n-1}'|   \leq 
 \tau[ B|u_{n-1}' - u_{n-2}'|  +  \\
& B' (\alpha_{n-2}+D)(|u_{n-1}-u_{n-2}|+|v_{n-1}-v_{n-2}|) ] \\
& + B(  |u_{n-2} -u_{n-3}|+  | v_{n-2} - v_{n-3}|) \\
& \leq \delta |u_{n-1}' - u_{n-2}'|  +\frac{C}{2^n} 
\end{align*}
for some constant $C$.
Then  we have the estimates:
\begin{align*}
 |u_n' & -u_{n-1}'|   \leq \delta |u_{n-1}' - u_{n-2}'|  + 
\frac{C}{2^n} 
 \leq \frac{1}{4} |u_{n-1}' - u_{n-2}'|  + 
\frac{C}{2^{n+1}} \\
& \leq \frac{1}{2^4} |u_{n-2}' - u_{n-3}'|  + 
\frac{C}{2^{n+2}} 
 \dots \leq  \frac{1}{2^{n-1}} |u_1'-u_0'| + \frac{C}{2^{2n-1}}.
\end{align*}

Similarly we have the estimates:
$$|\frac{\partial v_n}{\partial s}-\frac{\partial v_{n-1}}{\partial s}|
 \leq \frac{1}{2^{n-1}} |v_1'-v_0'| + \frac{C}{2^{2n-1}}$$
on $[ 0,  \tau] \times [0,  \infty)$.

So the convergence is  uniform.

\vspace{3mm}

{\bf Step 3:}
By step $1$ and $2$, we have the uniform estimates:
$$  |\frac{\partial u}{\partial x}|    \leq 2A+2D$$
on $[0,   \infty) \times [0,  \tau]$, and:
$$|\frac{\partial v}{\partial s}|  \leq 2A+2D$$
on $[ 0,   \tau] \times [0,  \infty)$
by letting $n \to \infty$.

Let us repeat the same process of extensions of the solutions,
on $[ 0,   \infty) \times [\tau, 2 \tau]$ for $u$ and 
on $[ \tau, 2  \tau] \times [0,  \infty)$ for $v$.

Notice that the initial norms have to replace from $A$ by $2A+2D$.
Then successively we have the estimates for $N=0,1,2, \dots$:
$$\begin{cases}
  |\frac{\partial u}{\partial x}|    \leq  2^NA+(2^N-1) 2D  & \text{ on }
 [ 0,   \infty) \times [(N-1)\tau, N \tau]\\
  |\frac{\partial v}{\partial s}|    \leq 2^NA+(2^N-1)2D & \text{ on } 
 [(N-1)\tau, N \tau]\times [ 0,   \infty)
\end{cases}$$

So in total, we have the following estimates:
\begin{align*}
&   |\frac{\partial u}{\partial x}| (x,s)   \leq 2^{\tau^{-1}s } (A  + 2D) \\
&  |\frac{\partial v}{\partial s}| (x,s) \leq  2^{\tau^{-1}x}  (A+2D), \\
&   |\frac{\partial u}{\partial s}|, \ \  |\frac{\partial v}{\partial x}| \leq D.  
\end{align*}

This completes the proof.
\vspace{3mm} \\
{\em Example 5.1:} 
Let us consider the  hyperbolic systems  of the form:
$$u_s = \frac{au}{1+u}-u, \quad v_x= g_t(u,v) -v.$$

Suppose the initial conditions satisfy   $u(1,0)=a-1$.
Then along the half line $\{(1,s): s \geq 0\}$, 
the ODE  $u_s= \frac{au}{1+u}-u$ has the unique solution $u(1,s)\equiv a-1$.
Then by differentiating the first equality by $x$ variable, one obtains the equation
$u_{xs} =(a^{-1}-1)u_x $,
whose solutions are given by:
$$u_x(1,s) = \exp((a^{-1}-1)s) u_x(1,0).$$
They are uniformly bounded   if $a>1$ hold.

If we choose
 $a <1$, $u_x$ grow  exponentially, even though
  $u$ take negative values.

\vspace{3mm}

For our purpose of the estimates on the second derivatives,
 this may not be so useful.

Let us induce the uniform  energy estimates by assuming negagive coefficients
on derivatives as below.

Firstly we have the general estimates:
\begin{lem}
Let $w(s)$ satisfy the estimates:
$$- a w(s) + b \leq w'(s) \leq -c w(s) +d $$
for some positive $a,b,c,d >0$. 
Then the estimates hold:
$$\frac{b}{a} + (w(0) -\frac{b}{a})\exp(-as) \leq w(s) \leq
\frac{d}{c} + (w(0) -\frac{d}{c})\exp(-as) .$$
\end{lem}
{\em Proof:}
Let us rewrite the inequalities as:
$$-a(w(s)-  \frac{b}{a}) \leq (w(s) -  \frac{b}{a})'.$$
Then we obtain the left hand side estimates of the conclusions immediately.

The right hand side can be treated similarly.
This  completes the proof.
 \vspace{3mm}

Let us consider the hyperbolic system of PDE:
$$u_s= f_t(v,u)-u,  \quad v_x=g_t(v,u)-v.$$

\begin{prop}
Assume negativities:
$$-a   \ \ \leq (f_t)_u-1 , \quad   (g_t)_v -1 \leq  \  \ -c$$
for some $0 <a,b$. 
Then the uniform estimates hold:
\begin{align*}
& \frac{b}{a} + (u_x(x,0) -\frac{b}{a})\exp(-as) \leq u_x(x,s) \leq
\frac{b}{c} + (u_x(x,0) -\frac{b}{c})\exp(-as),  \\
& \frac{d}{a} + (v_s(0,s) -\frac{d}{a})\exp(-ax) \leq v_s(x,s) \leq
\frac{d}{c} + (v_s(0,s) -\frac{d}{c})\exp(-as)
\end{align*}
where:
$$b =  \sup | (f_t)_v (g_t-v)|, \quad d=  \sup |(g_t)_u(f_t-u)|.$$

In particular $|u_x|$ and $|v_s|$ are both uniformly bounded.
\end{prop}
{ \em Proof:}
By differentiations, 
let us consider the equations:
\begin{align*}
& u_{xs} = ((f_t)_u -1)u_x + (f_t)_v (g_t -v),  \\
& v_{xs} = ((g_t)_v-1)v_s + (g_t)_u(f_t-u).
\end{align*}
Then the conclusions follow by applying lemma $5.9$.

This completes the proof.
\vspace{3mm} \\
{\bf 5.D Refinement and the higher distorsions:}
 Let us consider a Mealy automaton:
$${\bf A} : \quad \psi: Q \times S \to S, \quad \phi: Q \times S \to Q$$
 with two alphabets $S=\{s_0,s_1\} =\{L,L+1\}$.
Our aim in $5.D$ is to verify the following:
\begin{prop}
Let ${\bf A}$ be a Mealy automaton with $2$ alphabets.

There is an refinement of ${\bf A}$ with the pair of functions $(\bar{\phi}, \bar{\psi})$
so that the corresponding relatively elementary functions $(\tilde{f}_t, \tilde{g}_t)$
satisfy  the estimates:
 \begin{align*}
& [ \ |  ( \tilde{f}_t(v,u) -u) ( (\tilde{f}_t)_u(v,u)-1)|+
  |( \tilde{f}_t)_v(u,v)||v_s| \ ]
(x,s+ \alpha)   <2u(x,s+1),  \\
& [ \ |  ( \tilde{g}_t(v,u) -v) ( (\tilde{g}_t)_v(v,u)-1)|
+   |( \tilde{g}_t)_u(u,v)||u_x| \ ]
(x+ \alpha,s)   <2v(x+1,s).
 \end{align*}
for any solutions $(u,v)$ and all $0  \leq \alpha  \leq 1$.
 \end{prop}
Proof requires constructions of $(\max, +)$-functions
by several steps and occupies  $5.D$.
We also need  some general estimates of rational functions with positive coefficients.
\vspace{3mm} \\
{ \bf 5.D.2 Prototype :}
Let us describe a prototype of rational functions,
which appear  by refinement.

Let $\xi: {\mathbb R}  \to {\mathbb R} $ be another relatively $(\max, +)$-function,
given by:
\begin{align*}
\xi(x) & = \max( \min(x+ \delta, L) , x-\delta) \\
& = \max( -\max(-(x+ \delta) , -L) , x-\delta)  .
\end{align*}

$\xi(x) $ is increasing for $x <L$ and decreasing for $x  <L$.
So for any $q >\delta >0$, $\xi$ gives the funnctions as:
$$\xi : [L- q, L+q] \to [L-q, L+q].$$

Let $f_t$ be the corresponding relatively elementary functions to $\xi$:
$$f_t(z) = t^{-\delta}z +  \frac{1}{N_0} \ \frac{t^L}{  t^{L- \delta } + z}z$$
where $N_0 =N_0(t^{ \delta}  ) \geq 1$ are chosen so that the estimates hold:
$$ \mu \equiv 1- t^{-\delta} - \frac{t^{\delta}}{N_0} >0.$$

Let us consider the equalities:
 \begin{align*}
& (f_t(z) -z) '   =   t^{ -\delta} -1 + 
 \frac{1}{N_0} \frac{t^{2L - \delta}}{  (t^{L - \delta} + z)^2} ,  \\
& f_t(z) -z = (t^{- \delta}-1 + \frac{1}{N_0}
\frac{t^L}{t^{L-  \delta} +z} ) \ z.
 \end{align*}
Then the estimates hold:
$$     \frac{t^{2L-\delta}} {(t^{L - \delta} + z)^2} , \quad
\frac{t^L}{t^{L-  \delta} +z}  \leq t^{  \delta} .$$
Thus we have the estimates:
\begin{align*}
& -1
< t^{- \delta}-1 <  (f_t(z)-z)' \leq t^{-\delta} -1 + \frac{t^{\delta}}{N_0} =- \mu <0, \\
& -z <
(t^{-\delta}-1)z < f_t(z) -z  \leq  (t^{- \delta}-1 + \frac{t^{ \delta}}{N_0}) \ z = -  \mu z.
 \end{align*}

So in total we have the estimates:
$$|(f_t(z) -z)'| |f_t(z)-z| \leq (1- t^{- \delta} ) \ z.$$
This will be one of the required estimates for our asymptotic estimates of solutions to PDE systems.
 \vspace{3mm} \\
{\bf 5.D.3 Construction of exit functions:}
Let us construct admissible pairs concretely.
Such functions arise from  {\em almost diagonal functions} as below.
Let:
$$\xi(x)  = \max( \min(x+ \delta, L) , x-\delta) $$
be the relative $(\max, +)$-functions in $5.D$.
Given the initial value $x_0=0$, let us iterate it as
$x_{n+1}= \xi(x_n)$. It is easy to see: 
$$x_n  \equiv L  
\qquad  \text{ if } x_0=L, L+1  \ \text{ and } n \geq \delta^{-1}.$$

Let us have  another relative $(\max, +)$-functions with
 `two-step stairs' by:
$$\xi_2(x) = \max(  \min( \xi(x), L+1),    x-3\delta)  .$$

This satisfies the properties:
$$\xi_2(x) 
\begin{cases}
>x  &  \quad x < L  \\
=L & \quad L- \delta  \leq x \leq L +  \delta \\
 <x & \quad x > L
\end{cases}$$
The iterations satisfy the same properties,
$x_n \equiv L$ for all $n \geq  \delta^{-1}$ for  $x_0= L, L+1$.
On the other hand
let us translate $\xi_2$. Then its properties change as:
$$\xi_2(x+ 2\delta) 
\begin{cases}
>x  &  \quad x < L+1   \\
=L+1 & \quad L+1- \delta  \leq x \leq L + 1+ \delta \\
 <x & \quad x > L+1
\end{cases}$$
In particular the orbits behave differently.
The iterations satisfy:
 $$x_n \equiv L+1
\qquad  \text{ if } x_0=L , L+1 \ \text{ and } n \geq \delta^{-1}.$$

Inductively let us have relative $(\max, +)$-functions with
 $n$-step stairs by:
$$\xi_n(x) = \max(  \min( \xi_{n-1}(x), L+n-1),    x-(2n-1)\delta)  .$$
To represent exit functions, 
one can use $\xi_l$ when   we have $l$ alphabets.
In fact if we want an action which exchanges $L+a$ and $L+a+1$ with $a<l$, then 
the translations: 
\begin{align*}
 & \xi_l(\quad + 2a\delta): L+a, L +a +1\to L+a, \\
& \xi_l(\quad + 2(a+1)\delta): L+a, L +a+1 \to L+a+1
\end{align*}
can play such roles.
Notice that all $\xi_n$ are $1$-Lipschitz functions.

\vspace{3mm}

Let us fix $n$, and choose a large number $N_0=N_0(n)$
which will be determined later.

Let $f_t$ be the rational functions corresponding to $\xi$ above:
$$f_t(z) = t^{- \delta} z+ \frac{1}{N_0} \frac{t^Lz}{z+ t^{L- \delta}}.$$
For $0\leq m \leq n$
let us put rational functions inductively by:
$$f_t^m (z)  = t^{-(2m-1)\delta} z+  
\frac{1}{N_0}\frac{t^{L+m-1} f_t^{m-1}(z)}{t^{L+m-1}+ f_t^{m-1}(z)} $$
which  correspond to $\xi_m$ defined above.

\begin{lem} let us choose $N_0 \geq 1$ so that 
the estimates:
$$t^{-\delta} +\frac{t^{\delta}}{N_0}  \equiv 1- \mu <1$$
holds. Then
 the following estimates hold:
\begin{align*}
& -1 + t^{-(2n-1)\delta} < (f^n_t(z)-z)' \leq t^{-\delta} -1 + \frac{t^{\delta}}{N_0} =- \mu <0, \\
& (-1+ t^{-(2n-1) \delta})z < f^n_t(z) -z  \leq  (t^{- \delta}-1 + \frac{t^{ \delta}}{N_0}) \ z = -  \mu z.
 \end{align*}
In particular the estimates:
$$|(f_t(z) -z)'| |f_t(z)-z| \leq (-1+ t^{-(2n-1)\delta})^2 \ z.$$
 \end{lem}
{ \em Proof:}
We have already verifed the conclusions for $m=1$. 

Suppose  $f_t^{m-1}$ satisfy the conclusions for $m \leq n$.

Let us consider:
\begin{align*}
& |(f_t^m)' (z)| 
 = | t^{-(2m-1)\delta} + (\frac{ t^{L+m-1}}{t^{L+m-1}+ f_t^{m-1}(z)})^2 \frac{ (f_t^{m-1})'(z)}{N_0}|\\
& \leq t^{-(2m-1)\delta} +\frac{ |(f_t^{m-1})'(z)|}{N_0} \leq t^{-(2m-1)\delta} + \frac{1}{N_0} 
\leq t^{-\delta} + \frac{t^{\delta}}{N_0} < 1- \mu.
\end{align*}
So we have the estimates $-1 < (f_t^n(z)-z)' \leq -\mu<0$.

Next we have the estimates:
\begin{align*}
f_t^m (z)  &= t^{-(2m-1)\delta} z+\frac{1}{N_0}
 \frac{ t^{L+m-1}}{t^{L+m-1}+ f_t^{m-1}(z)} f_t^{m-1}(z)\\
& \leq t^{-(2m-1)\delta} z + \frac{1}{N_0} f_t^{m-1}(z) 
\leq  ( t^{-(2m-1)\delta} + \frac{1-\mu}{N_0} )z \\
& \leq (t^{-\delta} + \frac{t^{\delta}}{N_0})z <(1- \mu)z.
\end{align*}
Thus we have verified the conclusions for $m$.
This completes the proof.

 \vspace{3mm}

\begin{cor}
By choosing $N_0=N_0(t^{\delta})$ as above, 
the solutions to the equation:
$$u_s  = f_t^n(u) -u$$ satisfy the following estimates:
$$  |  (f_t^n(u) -u) ( (f^n_t)_u(u)-1)|(s+\alpha)   <\frac{\tau^2}{1-\tau^2} u(s+1)$$
for  $\tau = 1- t^{-(2n-1)\delta} $.
 \end{cor}
{ \em Proof:}
Let us put $\tau = 1- t^{-(2n-1)\delta}$.
By lemma $5.12$, the estimates
$|f_t^n(z)-z| < \tau^2z$ hold.  

Let us choose $0\leq \alpha_0 \leq 1$ so that 
$\sup_{0\leq s \leq 1} u(s)=  u( \alpha_0)$ holds.
Then we have 
the inequalities:
$$u(s+\alpha_0) \leq u(s+1) + \int_{s+\alpha_0}^{s+1} |u_s| <
u(s+1) + \tau^2 \int_s^{s+1} u(a) da. $$
By the mean value theorem, we have the estimates:
$$\int_s^{s+1} u(a)da = u(s+ \beta) \leq u(s+\alpha_0).$$
Thus combining with these estimates, we obtain the following:
$$(1- \tau^2) u(s+ \alpha_0) \leq  u(s+1).$$

Now finally we have the desired estimates:
\begin{align*}
&  |  (f_t^n(u) -u) ( (f^n_t)_u(u)-1)|(s+\alpha)  <
\tau^2 u(s+ \alpha)  \\
& \leq \tau^2 u(s+ \alpha_0) \leq \frac{\tau^2}{1- \tau^2} u(s+1).
\end{align*}
This completes the proof.
\vspace{3mm} \\
{\bf 5.D.4 Transition functions with $2$ alphabets:}
Now let us consider a Mealy automaton:
$${\bf A} : \quad \psi: Q \times S \to S, \quad \phi: Q \times S \to Q.$$
We will construct the transition functions by relatively $(\max, +)$-functions
 with two alphebets $S= \{s_0,s_1\} $.
The reason for this restriction is just for simplicity of the notations.
The general case can be constructed similarly.

Let us embed as $s_0 =L$ and $s_1=L+1$ for some large $L >>1$.
Let $Q=\{q^0, \dots, q^l\}$ be the set of the states.
We change the set, prepare twice of  the number of them as
$ \{q^0,\bar{q}^0, \dots, q^l,\bar{q}^l\}$, and embed them as:
$$\bar{Q} \equiv  \{q^0,\bar{q}^0,  \dots, q^l,\bar{q}^l\} \subset \delta  \ {\mathbb Z}$$
so that:
 $$q^j = \bar{q}^j+4\delta, \ \ q^{j+1}=  q^j + 8\delta$$
hold for  $j=0, \dots, l$.

\begin{lem}
There exists a refinement by
relatively $(\max,+)$-functions 
$\bar{\psi},\bar{\phi}: {\mathbb R} \times {\mathbb R} \to {\mathbb R}$.
\end{lem}
{\em Proof:}
{\bf Construction of $\bar{\psi}$ :}
Let us put the function:
 \begin{align*}
&  \mu(y) = \min[2\delta,  \max(  \tau(y-q^0), \dots,  \tau(y-q^l))],  \\
& \tau (y)=\min[  \max(0, y-\delta), \max(0, -y +7\delta))].
 \end{align*}
For $q \in Q$, these functions satisfy the following properties:
 $$\mu(y) =
\begin{cases}
0 & |y-q^j|  \leq  \delta, \\
2\delta & |\bar{q}^j -y|\leq \delta, \\
0 & y \geq q^l+ 7\delta  \text{ or } y \leq q^0 + \delta
\end{cases}$$

Now we put the desired functions by:
 $$\bar{\psi}(y,x) = \xi_2(x + \mu(y)).$$
{\bf Construction of $\bar{\phi}$ :}
Denote $\bar{*}^j = \bar{q}^j$
for  $* =q^j$ or $\bar{q}^j$, and
let $\epsilon_0$ be the permutation between $L$ and $L+1$. 
Let:
 $$\bar{\phi}: {\mathbb R} \times {\mathbb R} \to [q^0, q^l + 7  \delta] $$
be the relatively $(\max, +)$-function which satisfy the following properties:

\vspace{2mm}

(1)
$\bar{\phi}(*, \delta l ) =*$  if $\delta l \ne {\mathbb Z}$
 for  $l \in {\mathbb Z}$ and  $* \in \bar{Q}$.

\vspace{2mm}

(2) $|\bar{\phi}(y,x) -y |
\leq k\delta $ for all $x \in {\mathbb R}$ and 
$y \in [q^0, q^l+7\delta]$ for  some $k$.

\vspace{3mm}

(3) 
Suppose $\phi(q',s)=q$ and  let $\epsilon_0$ be the 
 permutation between $s_0 = L$ and $s_1= L+1$.
Then  for $*' = q' $ or $\bar{q}'$:
$$
\bar{\phi}(*', L) =
\begin{cases}
q & \text{ if }  \psi(q, \quad) = \epsilon_0, \\
\bar{q} &  \text{ if } \psi(q, \quad) = id
\end{cases}$$
$$
\bar{\phi}(*', L+1) =
\begin{cases}
\bar{q} & \text{ if } \psi(q, \quad) = \epsilon_0, \\
 q & \text{ if }  \psi(q, \quad) = id
\end{cases}$$

If we choose $\delta ^{-1}$ equal to some integer $l_0$, 
then these constructions give the refinement.
This completes the proof.
\vspace{3mm} \\
{\em Example:}
Let us consider the one state case $Q=\{q\}$. Then
we have representations of $\bar{\phi}$ by:
$$\bar{\phi}(y,x) =
\begin{cases}
 \min[  \bar{ q},  \max(y,  y+ x- (L+1- \delta)] & \text{ if } \psi(q, \quad) = \epsilon_0, \\
 \min[  \bar{q},   \max(y,y-x+  L+\delta)] &
\text{ if } \psi(q, \quad) = \text{ id }
 \end{cases}$$
\quad
\vspace{3mm} \\
{\bf 5.D.5 Estimates of relatively elementary functions:}
Let us give $C^0$ and $C^1$ estimates for  relatively 
elementary functions in terms of the corresponding $(\max, +)$-functions.
Let $f_t(z_0, \dots,z_{n-1})$ 
 correspond to  $\varphi(x_0, \dots,x_{n-1})$, and
fix $\bar{y}^0 \equiv (y_1, \dots, y_{n-1}) \in {\mathbb R}^{n-1}$. Then we consider: 
$$\varphi \equiv \varphi( \quad, y_1, \dots,y_{n-1}): {\mathbb R} \to {\mathbb R}.$$
Assume that $\varphi$ has the form:
 $$\varphi(x)= L_0+ \alpha (x-x_0)$$
 as a function 
  on some neighbourhood  $I_0=(x_0 - \tau, x_0+\tau)  \in {\mathbb R}$.

For $\bar{w}^0=t^{\bar{y}^0} \in {\mathbb R}^{n-1}_{>0}$ and $z_0 = t^{x_0}$,
 let us put :
$$f_t \equiv f_t(\quad, \bar{w}^0): (0, \infty) \to (0, \infty).$$

\begin{lem}
Let $M$ be the number of the components of $f_t$. 

(1) For any 
$t^{x_0 -\tau} \leq z \leq t^{x_0+\tau}$,  the estimates hold:
$$( \frac{ f_t(z)}{ t^{L_0}z_0^{- \alpha}z^{\alpha}})^{\pm 1}  \leq M.$$

(2)
Suppose $\varphi$ takes bounded values from both sides:
$$L \leq \varphi \leq L+a.$$
Then there exists some constant $C$ so that the estimates:
$$|\frac{\partial f_t}{\partial z}|(z)  \leq C t^a  \frac{t^L }{z}$$
hold for all $t>1$, where 
$C=C(\text{deg } h_t, \text{deg } k_t,M)$ with $f_t=t^L  \frac{h_t}{k_t}$.

(3) Suppose $a>0$. Then there exists an integer $N $ 
  so that
the estimates:
$$t^{L-a}  \leq N f_t(z) \leq t^{L+3a}$$
hold for  all $t^{L-a} \leq z \leq t^{L+3a}$ and  all sufficiently large $t >>1$.
\end{lem}
{\em Proof:}
(1) 
Let us put $x=\log_t z$ and $\varphi_t(x)=\log_t f_t(z)$.
By lemma $2.1$, the estimates hold for all $x \in I_0$:
\begin{align*}
|\varphi_t(x)-\varphi(x)|  & = |\log_t f_t(z) - \log_t (t^{L_0+\alpha (x-x_0)})| \\
& = \log_t ( \frac{ f_t(z)}{ t^{L_0}z_0^{- \alpha}z^{\alpha}})^{\pm 1} 
 \leq \log_tM
 \end{align*}
By removing $\log_t$ from both sides, one obtains the desired estimates.

(2) Notice that $f_t$ is a parametrized rational functions.
For each $t>1$, $f_t$ takes bounded values from both sides, 
since:
$$M^{-1}  t^{|\varphi(x)|} \leq  t^{|\varphi(x)|}t^{-|\varphi_t(x)-\varphi(x)|} \leq
|f_t(z)| \leq t^{|\varphi(x)|}t^{|\varphi_t(x)-\varphi(x)|} \leq M t^{|\varphi(x)|}.$$
In particular the degree of 
 $f_t$  must be equal to $0$ with respect to $z$.
So the derivative of $f_t$ with respect to $z$  has negative degree.

Let us denote  $f_t= t^L \frac{h_t}{k_t}$ by polynomials. Then by the above estimates, 
we obtain  the uniform bounds:
$$M^{-1} \leq  \frac{h_t}{k_t} \leq Mt^a$$ 
 where the coefficients of both $h_t$ and $k_t$ 
are  rational  in $t$.
Notice  that all the coefficients take positive values,
which is a characteristic in tropical geometry.

Suppose $h_t$ has degree $n_0  \geq 1$ in $z$. Then its derivative satisfies the estimates:
$$0 \  \leq \ z \ \frac{\partial h_t}{\partial z}(z) \ \leq \ c(n_0) h_t$$
for some constant $c(n_0)$ which is independent of $t>1$.

Now suppose deg $h_t=$ deg $k_t = n_0$, and denote 
$h_t' = \frac{\partial h_t}{\partial z}$.
Then the estimates hold:
$$ |\frac{\partial f_t}{\partial z}|(z) = t^L | \{ \frac{h_t'}{k_t} - \frac{h_t k_t'}{k_t^2} \}| \\
  \leq 2  t^L Mt^a \frac{c(n_0)}{z}= 2  c(n_0) Mt^a \frac{t^L }{z}.  $$
  
  (3) 
  Firstly we show that there exists an integer $N $  and some $L \leq L' \leq L+a$
  so that the estimates:
$$t^{L'} \leq N f_t(t^L) \leq t^{L+a}$$
hold for all sufficiently large $t >>1$.

  There is some $L'$ so that
  the  estimates $( \frac{ f_t(t^L)}{ t^{L'}})^{\pm 1}  \leq M$  hold by (1)
  for some $L \leq L' \leq L+a$ and all $t>1$.
  Since $ \frac{ f_t(t^L)}{ t^{L'}}$ is a rational function in $t>1$, 
  there is a limit:
   $$M^{-1} \leq \lim_{t \to \infty}  t^{-L'} f_t(t^L) \equiv  \kappa \leq M.$$
  One may assume the estimates $N M \leq t^a$   with $N=[\kappa^{\pm 1}]+1$ 
  for all  large $t >>1$, where $[\kappa^{\pm 1}] = \max([\kappa], [\kappa^{-1}]) \geq 0$.

  Then one obtains the estimates
   $1 \leq N t^{- L'}f_t(t^L) \leq NM \leq t^a$, and so:
     $$ t^L \leq t^{L'}  \leq N f_t(t^L)  \leq t^{L'+a}\leq t^{L+2a}$$
since $L \leq L' \leq L+a$ hold.  It follows   from (2) that  the estimates:
  $$N M^{-1}t^L  \leq  N f_t(t^x) \leq Nf_t(t^L) + CNt^{L+a}\log t^{4a}$$
  hold for all $L -a \leq x \leq L+3a$.
  
  One may assume the estimates for all sufficiently large $t >>1$:
  $$N^{-1}M \leq t^a, \quad CNt^{L+a}\log t^{4a} \leq  t^{L+2a}.$$
  Then combining with these estimates:
  $$t^{L-a} \leq Nf_t(t^x) \leq 2t^{L+2a} \leq t^{L+3a}$$
    hold for all $L -a \leq x \leq L+3a$.
This completes the proof.
\vspace{3mm} \\
{\bf 5.D.6 Proof of proposition $5.12$:}
\vspace{3mm} \\
{\bf Step 1:}
Let $\tau(y)$ be in $5.D.4$ and $l_t$ correspond to $\tau$:
$$l_t(w)= [ (1+t^{-\delta}w)^{-1}+ (1+ t^{7 \delta}w^{-1})^{-1}]^{-1}.$$
The inequalities hold:
$$ (1+t^{6\delta})^{-1} \leq l_t(w)^{-1} \leq 2, \quad
 |l_t(w)|' \leq 4(1+ t^{-\delta}).$$

Let:
$$h_t(w) =K^{-1} [t^{-2\delta} +( l_t(t^{-q^0}w) + \dots + l_t(t^{-q^l} w))^{-1}]^{-1}$$
 be corresponding to $\mu$, where $K\in {\mathbb N}$ are chosen so that 
the estimates $h_t(w) \leq 1$ hold.
Notice that the estimates hold:
$$[(l+1)(1+t^{6 \delta})]^{-1} + t^{-2\delta} \leq (Kh_t(w))^{-1} \leq \frac{2}{l+1} + t^{-2\delta}.$$
By inceasing the number of $l$ if necessarily, one may assume 
the lower bound $1- \chi \leq h_t(w)$ for sufficiently small $0< \chi <1$.

We have the estimates:
\begin{align*}
|h_t(w)'| & =\frac{ t^{-q^0+2\delta} }{K}
|\frac{\Sigma_{j=0}^l t^{-(q^j-q^0)} l_t'(t^{-q^j}w)}
{(t^{2\delta} + l_t(t^{-q^0}w) + \dots + l_t(t^{-q^l} w))^2}| \\
& \leq \frac{t^{2\delta}}{t^{q^0}}  \frac{4}{K(l+1)} \sup_{w \in {\mathbb R}}  |l_t'(w)| 
\leq C K^{-1}t^{-q^0}
\end{align*}
where $C$ is independent of $t>1$.

Now let $\bar{\psi}$ be in $5.D.3$, and recall
$\xi$ with corresponding  $f_t$.
Then corresponding to $\bar{\psi}$ is the following:
$$f_t(w,z) = f_t(z h_t(w)).$$

By lemma $5.13$, both the estimates:
\begin{align*}
& (-1+ t^{-3 \delta})z < f_t(w,z) -h_t(w) z <- \mu z,  \\
&   -1+ t^{-3 \delta} < (f_t)_z (w,z) -h_t(w)   <-  \mu .
\end{align*}
In particular we have the estimates:
\begin{align*}
& -z < (-2+ t^{-3 \delta}  + h_t(w)) z < f_t(w,z) -z <- \mu z,  \\
&   -1 < -2+ t^{-3 \delta}+h_t(w)  < (f_t)_z (w,z) -1   <-  \mu z,  \\
& | (f_t)_w(w,z)  |= |h_t'(w)|z |f_t'(zh_t(w))|   \leq C K^{-1}t^{-q^0}  z.
\end{align*}
\quad
\vspace{3mm}  \\
{\bf Step 2:}
Let $g_t(w,z)$ correspond to $\bar{\phi}$, and 
 estimate:
$$ |( f_t)_w(g_t-w)|,  \quad |(g_t)_z (f_t-z)|.$$

 Let us modify $g_t(w,z)$ so that  there are $C$ and  some $a$ independent of $t>1$,
 and the estimates hold:
\begin{align*}
& |(g_t)_z|,  \quad | (g_t)_w| \leq Ct^{a  \delta},   \\
& - w <  g_t(w,z)-w <   -\mu w,  \\
& -1  < (g_t)_w-1 < - \mu  
\end{align*}
for all $ t^{-a \delta} < w  \leq t^{a\delta}$.

Let us denote $\tilde{\phi}(y,x) = \bar{\phi}(y,x) -y$.
Then there is some $a$ with the estimates 
$|\tilde{\phi}(y,x)| \leq a\delta$.
So we have the equality
$\tilde{\phi}(y,x) =\max(\tilde{\phi}(y,x), - a\delta)$.

Let $\tilde{g}_t(w,z)$ correspond to $\tilde{\phi}$.
Then
 $\tilde{g}_t'(w,z)=\frac{ \tilde{g}_t }{N_0}+ t^{-a \delta}$
are tropically equivalent to $\tilde{g}_t$ for $N_0 \in {\mathbb N}$.
By lemma $5.16$, 
$\tilde{g}_t$ is uniformly bounded and so 
one may assume the estimates:
$$1- \chi  < \tilde{g}_t'(w,z) < 1+  \chi$$
for sufficiently small $0 <  \chi$.
Moreover: 
$$|\frac{\partial \tilde{g}_t}{\partial w}| \leq C t^{2a \delta}t^{-q^0}$$
hold.
Now we have the equality:
$$\frac{\partial \tilde{g}'_t}{\partial w} -1=
  \tilde{g}'_t-1 +   \frac{w }{N_0} \frac{\partial \tilde{g}_t}{\partial w}$$
where  the estimates
$| \frac{w }{N_0} \frac{\partial \tilde{g}_t}{\partial w}| 
\leq    Ct^{2a\delta} \frac{w}{N_0 t^{q^0} } \leq \frac{C t^{3a \delta}}{N_0} <<1$ hold, and so
we obtain the desired estimates.

Since
$$\bar{g}_t'(w,z) = \frac{w}{2} \ \tilde{g}_t'(w,z)$$
are tropically equivalent to $\bar{g}_t(w,z)$,
we obtain the desired rational functions.

So we obtain the followings:
$$ |( f_t)_w(g_t-w)| \leq CK^{-1}t^{a  \delta}  z   ,  \quad 
 |(g_t)_z (f_t-z)|\leq C t^{a\delta} z.$$
{\bf Step 3:}
By proposition $5.11$ and step 2, we have the  following estimates:
\begin{align*}
|(f_t)_v(v,u)| |v_s| & = |h_t'(w)|z |f_t'(zh_t(w))| |v_s| \\
& \leq C K^{-1}  t^{-q^0}  z t^{d \delta} \ z = CK^{-1} t^{-q^0+ a  \delta} z^2 \\
 |(g_t)_u(v,u)| |u_x| 
& \leq C t^{2a \delta}  z.
\end{align*}
If we take $q^0 \geq L+1$ and large $t>>1$, then we obtain the estimates:
$$
|(\tilde{f}_t)_v(v,u)| |v_s|(x,s) 
 \leq   \mu' u(x,s), \quad 
 |(\tilde{g}_t)_u(v,u)| |u_x| (x, s)
 \leq  \mu'  v(x,s)$$
for sufficiently small $ 0<\mu' <<1$. 

Now combining with corollary $5.14$, 
we obtain the estimates:
\begin{align*}
|(\tilde{f}_t)_v(v,u)| |v_s|(x,s+ \alpha) 
& \leq  \mu''  u(x,s+1), \\
 |(\tilde{g}_t)_u(v,u)| |u_x| (x+\alpha, s)
& \leq  \mu'' v(x+ 1,s)
\end{align*}
for some $0<  \mu'' <2$. 
This completes the proof of proposition $5.15$.

 \vspace{3mm}

Basically the constructions of the refnement  are quite general, 
and it would be reasonable to expect the following:
\vspace{3mm} \\
{\em Conjecture 5.1:}
For any pairs of relatively $(\max,+)$-functions $(\phi, \psi)$,
there are refinements $(\bar{\phi}, \bar{\psi})$ by $1$-Lipschitz functions
so that the corresponding 
relatively elementary functions $(\bar{f}_t, \bar{g}_t)$ are admissible.

\vspace{1cm}

Tsuyoshi Kato

Department of Mathematics

Faculty of Science, Kyoto University

Kyoto 606-8502 Japan

 \end{document}